\documentclass[a4paper,12pt]{article}

\usepackage[T1]{fontenc}
\usepackage{lmodern}          
    
\usepackage{mathtools}
\usepackage{amssymb}
\usepackage{amsthm}
\usepackage{thmtools}
\usepackage{thm-restate}
\usepackage{microtype}

\usepackage[margin=1in,includefoot,footskip=30pt]{geometry}
\usepackage{graphicx}
\usepackage{tikz-cd}
\usepackage{xcolor}
\usepackage{enumitem}
\setlist{nosep,leftmargin=*}
\usepackage{comment}
\usepackage[toc]{appendix}      
\usepackage[round]{natbib}      % author–year; \citet, \citep, \citeyear, \citealp
\usepackage{hyperref}

\usepackage[capitalise,noabbrev]{cleveref}
\newtheorem{theorem}{Theorem}[subsection]
\newtheorem{lemma}[theorem]{Lemma}
\newtheorem{proposition}[theorem]{Proposition}
\newtheorem{corollary}[theorem]{Corollary}
\theoremstyle{definition}
\newtheorem{definition}[theorem]{Definition}
\newtheorem{example}[theorem]{Example}
\theoremstyle{remark}
\newtheorem{remark}[theorem]{Remark}

\numberwithin{equation}{subsection}

\newcommand{\pw}{\text{pw}}
\newcommand{\Mod}{\mathrm{Mod}}
\newcommand{\Th}{\mathrm{Th}}
\newcommand{\Thpw}{\Th^{\pw}}

\newcommand{\eqv}{\mathrel{\dashv\!\;\vdash}}

\newcommand{\id}{\mathrm{id}}
\hypersetup{
  hidelinks,
  pdfauthor={Benjamin Zayton},
  pdftitle={Piecewise Interpretations in Many-Sorted Logic}
}

\title{Piecewise Interpretations in Many-Sorted Logic}
\author{Benjamin Zayton\thanks{  Department of Philosophy, University of Vienna.
  E-mail: \href{mailto:benjamin.zayton@univie.ac.at}
  {\texttt{benjamin.zayton@univie.ac.at}}.}}
\date{}
\begin{document}
\maketitle
\begin{abstract}
\noindent
Piecewise interpretations, introduced by Visser, generalize ordinary interpretations between first-order theories by allowing formulae to be interpreted using disjoint unions of definable sets. This makes it possible to capture natural constructions, such as the projective completion of an affine plane, but also leads to a notion of equivalence of theories, namely piecewise bi-interpretability. We develop the theory of piecewise interpretations systematically in the many-sorted setting. Each piecewise interpretation induces a semantic dual, a functor between categories of models. In addition, we introduce a syntactic notion of morphism between piecewise interpretations, which likewise has a semantic dual. Using these semantic duals, we collect theories, piecewise interpretations, and morphisms between them into a bicategory, in particular defining a well-behaved notion of composition of interpretations, associative up to coherent isomorphism. As an application, we consider a familiar idea, championed by Quine, namely that one can replace any many-sorted theory with a corresponding single-sorted one by replacing sorts with sortal predicates. We show that this leads to a biequivalence between the developed bicategories of theories with finitely many sorts and single-sorted theories, thereby providing an argument for the dispensability of sorts.
\end{abstract}

\section{Introduction}
Interpretations between first-order theories provide a fundamental tool
for understanding relationships between mathematical structures and theories. In the time since Tarski's introduction of his relative interpretations in \textit{Undecidable Theories} \citeyearpar{tarski1953undecidable}, various generalizations of the classical notion have arisen \citep{alscher2016theorien}. These notions, and their philosophical applications, have received significant attention and are relatively well understood.\footnote{For instance, \citet{button_philosophy_2018} investigate bi-interpretability as a standard of sameness of theory, \citet{niebergall2000logic} defends interpretability as a notion of reducibility, and \citet{engler2026bedeutung} investigates the relationship between interpretability and synonymy.} The goal of this paper is to systematically study the far less thoroughly investigated \emph{piecewise interpretations}, introduced by \citet{visser2004categories, visser2021extension} in a single-sorted setting. As opposed to ordinary interpretations, these allow one to interpret formulae using disjoint unions of definable sets. The theory developed in the paper is quite general, encompassing ordinary interpretations. It also introduces a natural notion of equivalence of theories, namely piecewise bi-interpretability, and remedies certain defects of the ordinary notion of interpretation, to do with the reduction of many-sorted to single-sorted logic.

\medskip
\noindent
The main focus will be the development of the mathematical theory of piecewise interpretations in the many-sorted setting. In \cref{sec:def}, we define many-sorted piecewise translations and interpretations, discuss the projective completion of an affine plane as an example, develop semantic duals of translations, and introduce morphisms between piecewise translations ($i$-morphisms) and their semantics. In \cref{sec:bicategory}, we show that theories, piecewise interpretations, and $i$-morphisms form a bicategory $\Thpw$ by semantic means, using the induced semantic duals of interpretations. In particular, we show that the developed notion of composition of piecewise interpretations is associative (up to coherent $i$-isomorphism), thereby answering a small question due to Visser. In \cref{sec:application}, we present a (pseudo)functorial reduction of many-sorted to single-sorted logic, based on the familiar idea that one can replace sorts with sortal predicates. This reduction will take the form of a biequivalence between the developed bicategory of theories (with finitely many sorts) and piecewise interpretations and the single-sorted analogue of that bicategory. 

\medskip\noindent
Our decision to place the study of interpretations against the background of category theory follows \citet{visser2004categories}, while our approach to interpretations in many-sorted logic is inspired by \cite{halvorson2019logic}. \cite{meadows2024beyond} pursues a structurally similar but more general and semantic project, studying functorial comparisons between theories mediated through their categories of models, organized through $2$-category theory. Taking a mostly syntactic route, \citet{d2020bicategories} systematically develop a bicategory of interpretations in coherent logic, establishing a tight relationship to categorical logic. Beyond Visser, piecewise interpretations have been studied independently by \citet{li2025theoretical} in his master's thesis: he presents a definition extending Visser's against a many-sorted background, connecting piecewise interpretations to Meadows' Morita interpretations and providing a notion of composition. The definition we have chosen differs in presentation from his and Visser's, taking a cue from categorical logic, and we develop the theory further by placing our notion of composition into a bicategorical framework.

\subsection{Notation and Preliminaries}

Given a set $X$, we denote by $X^*$ the set of finite sequences of elements of that set. Elements of $X^*$ will be denoted by $\bar{x} = (x_1, \ldots, x_n)$. We also write $\bar{x} = x_1 \times x_2 \times \cdots \times x_n$, identifying $x\in X$ with the sequence $(x)$, where $\times$ denotes concatenation of sequences. Concatenation with the empty sequence $\varepsilon$ is the identity. For a sequence $\bar{x}$, we denote by $|\bar{x}|\in \mathbb{N}$ its length and by $x_i$ its $i$-th component. For $n\in\mathbb{N}$ we write $[n]:=\{1,\ldots,n\}$, with $[0]:=\varnothing$. Identifying a sequence $\bar x=(x_1,\ldots,x_n)$ over a set $X$ with the function $[n]\to X$ sending $j$ to
$x_j$ and given a function $s:[m]\to[n]$, we write $\bar x^{s}$ for the composite $\bar x\circ s$, componentwise described as $\bar x^{s}:=(x_{s(1)},\ldots,x_{s(m)})$.
This is functorial: $\bar x^{\id}=\bar x$, and $(\bar x^{s})^{t}=\bar x^{s\circ t}$ for $t:[l]\to[m]$, by associativity of composition.

\medskip
\noindent
A \emph{finite family} over a set $X$ is a function $p: I_p \to X$ where $I_p$ is a finite, non-empty set, called the index set of the family. We denote by $ \mathcal{F}(X) := \{p : I_p \to X \mid I_p \text{ is finite and non-empty}\}$ the collection of all finite families over $X$. We will frequently need to consider $\mathcal F(X^*)$, the collection of finite families of finite sequences over $X$. Given a finite family of sets $p:I_p\to X$, we denote by $\bigsqcup_{k\in I_p}p(k)$ its disjoint union, set-theoretically constructed as the set of ordered pairs $(k,a)$ such that $a\in p(k)$.

\medskip
\noindent
Next, we fix our logical vocabulary. A relational many-sorted first-order signature $\Sigma=(\mathcal{S},\mathcal{R},\mathcal{V})$ consists of a non-empty set $\mathcal{S}$ of sort symbols, a set $\mathcal{R}$ of sorted relation symbols, and a family $\mathcal{V}=(\mathcal{V_\sigma})_{\sigma\in \mathcal{S}}$ of sets of sorted variables.\footnote{We work without function and constant symbols, analysing them as (families of) definite descriptions in the Russellian manner.} To ensure we have enough variables available, we fix an infinite cardinal $\kappa$ and restrict attention to signatures with $|\mathcal{R}\cup\mathcal{S}|\leq\kappa$ and $|\mathcal{V}_\sigma|=\kappa$ for all $\sigma\in\mathcal{S}$. Each relation symbol $R\in \mathcal{R}$ has an arity $\bar{\sigma}\in \mathcal{S}^*$, written $R:\bar{\sigma}$. Among the logical vocabulary, we count verum $\top$ and falsum $\bot$, taking their arity to be the empty sequence $\varepsilon$. For sentential operators, we have $\{\neg,\land,\lor,\to\}$. Moreover, each variable $x \in \mathcal{V}$ (resp. sequence $\bar{x}\in \mathcal{V}^*$ of variables) comes with an arity ${\sigma}\in \mathcal S$ (resp. sequence of sorts $\bar{\sigma}$ of the same length), expressed by $x:\sigma$ (resp. $\bar{x}:\bar \sigma$). Quantifiers are sorted as well, with one quantifier $\forall_\sigma$ for each sort $\sigma\in \mathcal{S}$, although we will only explicitly represent this sorting with a subscript where helpful. For a sequence of variables $\bar y=(y_1,\ldots,y_n):\bar \sigma$ we define $\forall \bar y\,\psi:=\forall_{\sigma_1}y_1\ldots\forall_{\sigma_n}y_n\psi$. In particular, $\forall\varepsilon\psi=\psi$. Officially, we will treat the existential quantifier $\exists$ as the dual of the universal quantifier. The signature $\Sigma$ generates a sorted first-order language $\mathcal{L}_\Sigma$ in the standard inductive fashion: given a relation symbol $R:\bar \sigma$ and a sequence of variables $\bar x:\bar \sigma$, one can form an atomic formula $R(\bar x)$, and general formulae are generated from these via the connectives and the quantifiers, taking sorting into account.\footnote{A consequence of our treatment of variables is that each such language $\mathcal{L}_\Sigma$ has cardinality $\kappa$.} 

\medskip\noindent For a syntactic calculus for many-sorted logic, we direct the reader to \cite{halvorson2019logic}; starting from a syntactic calculus for single-sorted logic, one only needs to modify the rules for the quantifiers so as to take the sorting into account. We will denote the resulting consequence relation by $\vdash$. A \emph{theory} $T$ in a signature $\Sigma$ is a deductively closed set of sentences in the language $\mathcal{L}_\Sigma$.
 
\medskip
\noindent
We will also work with formulae-in-context, that is, pairs of a formula $\varphi\in \mathcal{L}_\Sigma$ and a context $\bar x$ for $\varphi$. A context for a formula is a duplicate-free list $\bar{x}$ of variables which contains all the free variables of $\varphi$. The notation $(\bar{x}\mid \varphi)$ denotes the formula-in-context consisting of $\varphi$ and the context $\bar{x}$ for $\varphi$. Of particular importance will be the \textit{canonical context}\footnote{See Appendix~\ref{app:cc_and_subst} for an inductive definition.} $\bar{x}(\varphi)$ of a formula $\varphi$, which contains exactly the free variables of $\varphi$ in order of first occurrence. The sort sequence $\bar{\sigma}(\varphi)$ of $\bar{x}(\varphi)$ is called the \textit{arity} of $\varphi$; we write $\varphi:\bar \sigma$. 

\medskip
\noindent
Given a context $\bar{y} = (y_1,\ldots,y_m) : \bar{\tau}$ for a formula $\varphi$ with canonical context $\bar{x}(\varphi) = (x_1,\ldots,x_n) : \bar{\sigma}$, there is a unique injection $c : [n] \to [m]$ such that $x_j = y_{c(j)}$ for $1 \leq j \leq n$. We call $c$ the \emph{context injection} for $\bar{y}$ relative to $\varphi$. Given a formula $\varphi$ and a sequence of variables $\bar{y}:\bar{\sigma}(\varphi)$, we write $\varphi(\bar y)$ for the result of the capture-avoiding simultaneous substitution of $\bar y$ for $\bar{x}(\varphi)$ in $\varphi$.\footnote{See Appendix~\ref{app:cc_and_subst} for an inductive definition.}

\begin{restatable}[Canonical Context under Substitution]{lemma}{ccsubst}\label{lem:canon-subst}
Let $\varphi$ be a formula with canonical context $\bar{x} = (x_1,\ldots,x_n) : \bar{\sigma}$, and let $\bar{y} = (y_1,\ldots,y_n) : \bar{\sigma}$. Then the canonical context of $\varphi(\bar{y})$ is $\bar{y}^\circ$, the sequence resulting from $\bar y$ by deleting repeated members from right to left.
\end{restatable}
\begin{proof}
See Appendix~\ref{app:cc_and_subst}.
\end{proof}
\noindent
Let us record some other useful abbreviations: given a sequence of variables $\bar x:\bar \sigma$, we write
$\top_{\bar{x}} := \bigwedge_{i=1}^{|\bar{x}|} (x_i =_{\sigma_i} x_i)$ and $\bot_{\bar{x}} := \neg\top_{\bar{x}}$, where the empty conjunction is the nullary $\top$. We write $\top_{\bar \sigma}$ and $\bot_{\bar \sigma}$ if the choice of duplicate-free variable sequence does not matter. Moreover, given another sequence $\bar y:\bar \sigma$, we define identity for sequences of the same arity as $\bar x=_{\bar \sigma}\bar y:=\top_{\bar x\times\bar y}\land\bigwedge_{i=1}^{|\bar x|}(x_i=_{\sigma_i}y_i)$, with the first conjunct serving as padding to get the arities right. 
 
\section{Piecewise Interpretations: Basic Definitions and Semantics}\label{sec:def}
In this section, we will define piecewise translations and interpretations between theories, study their semantic duals, present some basic examples, and introduce the notion of an $i$-morphism between piecewise translations.
\subsection{Piecewise Translations and Interpretations}
The following definition is inspired by Visser's definition of piecewise translation and extends it to the many-sorted setting, using Halvorson's approach to many-sorted translations as a basis \citep{visser2004categories, visser2021extension, halvorson2019logic}. The main idea behind the definition is that sorts can be taken to finite families of sequences of sorts, indexed by the \emph{set of pieces} for the sort, and relation symbols are taken to finite families of accordingly sorted formulae.  
\begin{definition}[Piecewise Reconstrual]\label{def:reconstrual}
    Let $\Sigma_1$ and $\Sigma_2$ be two first-order signatures. A \emph{piecewise reconstrual} $f:\Sigma_1\to \Sigma_2$ is a tuple $(f_\mathcal{S}, f_\mathcal{V}, f_\mathcal{D}, f_\mathcal{R})$, consisting of the following data: The \emph{sort function} $f_\mathcal{S}$ assigns each sort $\sigma\in \mathcal{S}_1$ a finite family of finite sequences of sorts $(\sigma_f^k)_{k\in P_\sigma}$ of $\mathcal{S}_2$ with hereditarily finite index set $P_\sigma$, called \emph{the set of pieces for $\sigma$}.\footnote{The assumption that index sets be hereditarily finite allows one to globally well-order these, and will ensure that certain size conditions can be satisfied. Through composition of the translations induced by reconstruals, sets of natural numbers as index sets would lead one to sets of these sets as index sets, and hence ultimately to hereditarily finite sets.} The \emph{variable function} $f_\mathcal{V}$ assigns to each variable $x:\sigma$ in $\mathcal{V}_1$ a finite family of sequences of variables $(x_f^k)_{k\in P_\sigma}$ in $\mathcal{V}_2$ such that $x_f^k:\sigma_f^k$ is duplicate-free and appropriately sorted, and such that $x_f^k$ and $y_f^l$ are pairwise disjoint if $x\neq y$ (no-overlap).
    The \emph{domain function} $f_\mathcal{D}$ assigns to each sort $\sigma\in \mathcal{S}_1$ a finite family of formulae $(D_\sigma^k)_{k\in P_\sigma}$ of $\mathcal{L}_{\Sigma_2}$, such that $D_\sigma^k:\sigma_f^k$.
    The \emph{relation function} $f_\mathcal{R}$ assigns to each relation symbol $R:\bar\sigma=\sigma_1\times\cdots\times\sigma_n$ of $\Sigma_1$ a family of $\mathcal L_{\Sigma_2}$-formulae $(R_f^{\bar k})_{\bar k\in P_{\bar\sigma}}$, where $P_{\bar\sigma}:=P_{\sigma_1}\times\cdots\times P_{\sigma_n}$ and $R_f^{\bar k}:(\sigma_1)_f^{k_1}\times\cdots\times(\sigma_n)_f^{k_n}$. Note that we take $m$-fold Cartesian products to consist of $m$-tuples.

    \medskip\noindent
We extend the data of a reconstrual to sequences as follows: For a sequence $\bar \sigma=(\sigma_1,\ldots,\sigma_n)\in \mathcal{S}_1^*$, we set $\bar{\sigma}_f^{\bar{k}} := \sigma_{1,f}^{k_1}\times \ldots\times \sigma_{n,f}^{k_n}$ for $\bar{k} = (k_1,\ldots,k_n) \in P_{\bar{\sigma}}$.
 For the empty sequence $\varepsilon$, we set $P_\varepsilon=\{\varepsilon\}$. Given a sequence $\bar{x} = (x_1,\ldots,x_n)$ with $\bar{x}:\bar{\sigma}$ and $\bar{k} \in P_{\bar{\sigma}}$, we set $(\bar x_f^{\bar k}):=(x_1)_f^{k_1}\times\cdots\times(x_n)_f^{k_n}$. In particular, if $\sigma_f^k=\varepsilon$, then $x_f^k:\varepsilon$. To extend the domain function to sequences of sorts, given $\bar\sigma=(\sigma_1,\ldots,\sigma_n)$ we take some duplicate-free sequence of variables $\bar x:\bar \sigma$ and set $D_{\bar\sigma}^{\bar k}:=\bigwedge_{i=1}^n D_{\sigma_i}^{k_i}((x_i)_f^{k_i})$, to ensure that there is no accidental sharing of free variables. For the empty sequence, we set $D_\varepsilon^\varepsilon:=\top$. 
\end{definition}
\begin{remark}[Convention for Choice of Variables]
From now on, whenever a statement about translations or, later, $i$-morphisms involves the choice of some sequences of variables, we tacitly assume that these sequences are chosen to be fresh, pairwise distinct/disjoint, duplicate-free and such that their images under all pertinent translations and pieces are pairwise disjoint.
\end{remark}
Where necessary, we will make the dependence on $f$ of the data above explicit by using sub-/superscripts ($P^f_\sigma, D^f_\sigma$). Two comments on this definition: Due to the requirement on arities, replacing a translated formula with a logically equivalent one in a reconstrual will not necessarily result in a reconstrual.
Moreover, in practice, one does not always work with a globally fixed variable function, but instead just locally chooses fresh variables when translating any given formula. To partially alleviate the strictness of the definition and dependence on the globally chosen variable function, we will make use of the following equivalence relation on reconstruals:
\begin{definition}[Provable Equivalence of Reconstruals]\label{def:prov_equiv_reconstrual}
Let $f,g:\Sigma_1\to \Sigma_2$ be two piecewise reconstruals and let $T$ be a $\Sigma_2$-theory. Then $f$ and $g$ are $T$-provably equivalent, denoted by $f\equiv_Tg$, iff $\sigma_f=\sigma_g$ for all $\sigma\in\mathcal{S}_1$, and for all $\sigma\in\mathcal{S}_1$, $k\in P_\sigma^f$, and $\bar u:\sigma_f^k$, $T\vdash D_\sigma^{f,k}(\bar u)\leftrightarrow D_\sigma^{g,k}(\bar u)$, and, for all relation symbols, including identity, $R:\bar \sigma$, pieces $\bar k\in P_{\bar \sigma}^{f}$ and $\bar u:\bar \sigma_f^{\bar k}$, $T\vdash D_{\bar\sigma}^{f,\bar k}(\bar u)\to
            \bigl(R_f^{\bar k}(\bar u)\leftrightarrow R_g^{\bar k}(\bar u)\bigr)$.
\end{definition}
 \begin{remark}[Presentation of Reconstruals]\label{rem:presentation_reconstruals}
     The assumption that the sizes of signatures are bounded by an infinite cardinal $\kappa$ and that each sort always has $\kappa$-many variables has the consequence that, given a sort function, domain function, and relation function as in \cref{def:reconstrual}, there always is a variable function such that these data comprise a reconstrual. Moreover, the choice of variable function does not impact $\equiv_T$-equivalence. Therefore, to define a $\equiv_T$-equivalence class of a reconstrual, it is enough to write down a sort function, and the domain and relation function, up to substitution of variables in the target language. 
\end{remark}
\noindent
To define a proper notion of translation between the full languages $\mathcal{L}_{1}$ and $\mathcal{L}_{2}$ generated by two signatures $\Sigma_1$ and $\Sigma_2$, we recursively extend any given reconstrual. As an intermediate step, this will require a definition of the translation for formulae-in-context, which then restricts to the desired translation via the canonical context. The recursion reduces the translation of a formula in an arbitrary context to its translation in canonical context and the translation of complex formulae in canonical context to the translation of their components in general contexts; the atomic constituents in canonical context are then handled by the underlying reconstrual.
\begin{definition}[Piecewise Translation]\label{def:translation_extension}
Let $f:\Sigma_1\to \Sigma_2$ be a piecewise reconstrual. It extends to a function $f$ from formulae-in-context of $\mathcal{L}_{1}$ to families of formulae of $\mathcal{L}_{2}$ by the following recursive clauses. If the context $\bar x$ is the canonical context $\bar x(\varphi)$, we write $\varphi$ for $(\bar x(\varphi)\mid \varphi)$.

\begin{enumerate}[label=(\roman*)]

\item Truth Values: $\top_f^{\varepsilon}=\top$ and $\bot_f^{\varepsilon}=\bot$.

\item Contexts: Given a formula $\varphi$ with canonical context $\bar{x}(\varphi)=(x_1,\ldots,x_n):\bar\sigma$ and a context $\bar{y}=(y_1,\ldots,y_m):\bar \tau$ for $\varphi$ with $\bar x\neq \bar y$, consider the unique injection $c:[n]\to[m]$ such that $x_i=y_{c(i)}$ for all $1\leq i\leq n$. Via $c$, a piece $\bar k\in P_{\bar \tau}$ induces a piece $\bar k^c\in P_{\bar \sigma}$ by $(\bar k^c)_i:=k_{c(i)}$ for $1\leq i\leq n$. We can hence set $(\bar y\mid \varphi)_f^{\bar k} := \varphi_f^{\bar k ^c}$. 
\item Atomic: Given an atomic formula $R(z_1,\ldots,z_m)$ with $\bar z:\bar\tau$ and canonical context $\bar x=(x_1,\ldots,x_n):\bar \sigma$, note that there is a unique function $d:[m]\to[n]$ such that $x_{d(i)}=z_i$. Via $d$, a piece $\bar k\in P_{\bar \sigma}$ induces a piece $\bar k^d\in P_{\bar \tau}$ by $(\bar k^d)_i=k_{d(i)}$ for $1 \leq i\leq m$. We can hence set: $
R(\bar z)_f^{\bar k}:=R_f^{\bar k^d}(\bar z_f^{\bar k^d})$.
\item Negation: Given a negated formula $\neg \varphi:\bar \sigma$ and a piece $\bar k\in P_{\bar \sigma}$, we set $(\neg\varphi)_f^{\bar k} := \neg (\varphi_f^{\bar k})$

\item Binary: Given a complex formula $\varphi = \varphi_1 \star \varphi_2:\bar\sigma$ where $\star \in\{\land,\lor,\to\}$ with canonical context $\bar x=\bar x(\varphi_1\star \varphi_2)$ and a piece $\bar k \in P_{\bar \sigma}$, we set $\varphi_f^{\bar k} := (\bar x\mid\varphi_1)_f^{\bar k} \star (\bar x\mid \varphi_2)_f^{\bar k}$.

\item Quantifier: Given a quantified formula $\varphi = \forall_\tau y\, \psi$ with $y:\tau$ and canonical context $\bar x:=\bar x(\varphi):\bar \sigma$, and a piece $\bar k\in P_{\bar \sigma}$, note that $\bar{x}\times y:\bar \sigma\times\tau$ is a context for $\psi$. We can now set:
\begin{equation}\label{quantifier_translation}
    \varphi_f^{\bar k} :=\bigwedge_{h \in P_{\tau}} \forall y_f^{h} \left( D_\tau^{h}(y_f^h) \to (\bar x\times y\mid \psi)_f^{\bar k\times h} \right) 
\end{equation}
When forming a conjunction indexed by a set of pieces, we will generally assume that the conjuncts are ordered in accordance with the well-order on the hereditarily finite sets.
\end{enumerate}
In general, a (piecewise) translation  $f:\mathcal{L}_{1}\to\mathcal{L}_{2}$ is a function from $\mathcal{L}_{1}$ to $\mathcal{F}(\mathcal{L}_{2})$ induced in this way by a specific reconstrual from $\Sigma_1$ to $\Sigma_2$. Two piecewise translations $f,g$ are called provably equivalent (with respect to an $\mathcal{L}_2$-theory $T$) if the underlying reconstruals are provably equivalent (with respect to $T$).\footnote{Our approach to piecewise translations is presentationally different from Visser's even in the single-sorted setting, as he assigns pieces directly to the free variables of a formula, roughly speaking. His approach avoids canonical contexts and is more amenable to recursive extension, but the approach via canonical contexts allows for easier control of substitutions and hence compositions. Independently, \citet{li2025theoretical} extends Visser's approach to the many-sorted setting. We think that, ultimately, his definition should turn out to be equivalent to the one chosen here, but we have not systematically investigated the matter.}
\end{definition}
\noindent
Equivalently, a translation can be seen as taking a formula $\varphi:\bar \sigma \in \mathcal{L}_{1}$ and a piece $\bar k \in P_{\bar \sigma}$ to a formula $\varphi_f^{\bar k} \in \mathcal{L}_{2}$. The behaviour of this translation with respect to free variables and arities is as expected:
\begin{restatable}[Canonical Context of Translations]{lemma}{cctrans}\label{lem:variables_translations}
Let $f:\Sigma_1\to\Sigma_2$ be a piecewise reconstrual and let $\varphi$ be an $\mathcal{L}_{1}$-formula with canonical context $\bar{x}:\bar \sigma$. Then, for any piece $\bar k \in P_{\bar \sigma}$, the canonical context of $\varphi_f^{\bar k}$ is $\bar x_f^{\bar k}$. In particular, the formula $\varphi_f^{\bar k}$ has exactly the free variables $\bar{x}_f^{\bar k}$. 
\end{restatable}
\begin{proof}
See Appendix~\ref{app:cc_and_subst}.
\end{proof}
\noindent
Given an $\mathcal{L}_1$-sentence $\varphi$, its canonical context is the empty sequence $\varepsilon$ with set of pieces $P_{\varepsilon}=\{\varepsilon\}$. By the preceding lemma, the translation $\varphi_f^{\varepsilon}$ of $\varphi$ at the unique piece $\varepsilon$ has canonical context $\varepsilon$, and is hence a sentence. More generally, for any context $\bar y:\bar \sigma$ and $\bar k\in P_{\bar \sigma}$, $(\bar y\mid\varphi)_f^{\bar k}=\varphi_f^\varepsilon$. Given a translation $f:\mathcal{L}_1\to\mathcal{L}_2$ we hence obtain a well-defined function between sets of sentences:
\begin{align} \label{eq:sentence_preserv}
f:\mathcal{L}_{1}^{\textsc{Sent}} &\to \mathcal{L}_{2}^{\textsc{Sent}}\\
\varphi&\mapsto \varphi_f:=\varphi_f^{\varepsilon}
\end{align}

\noindent
This provides a first indication that piecewise translations can serve as the basis of an appropriate notion of interpretation between theories. 
In analogy to the ordinary theory of interpretations and relying on piecewise translations mapping sentences to sentences, we define:
\begin{definition}[Piecewise Interpretation]
    Let $S$ be a $\Sigma_1$-theory and $T$ be a $\Sigma_2$-theory. Then a piecewise translation $f:\mathcal{L}_1\to \mathcal{L}_2$ is an interpretation of $S$ in $T$ iff 
    for all sentences $\varphi\in\mathcal{L}_1$:
    \begin{equation}
        S\vdash\varphi\Rightarrow T\vdash \varphi_f. \label{theoremhood_preserv}
    \end{equation}
    In that case, we write $f:S\to T$, and say that $T$ interprets $S$ (piecewise) via $f$.
\end{definition}
If one specializes to the case where the set of pieces has exactly one piece for each sort, one recovers the ordinary notion of interpretation between many-sorted theories, as developed in \cite{halvorson2019logic}, in turn generalizing the usual notion between single-sorted theories.
\subsection{Identity under a Piecewise Interpretation}\label{sec:identity}
An interesting feature of the ordinary notion of interpretation is that it allows for the interpretation of identity by more general equivalence relations. To make room for this feature, we follow Halvorson in treating identity as a non-logical symbol, and hence tacitly think of each $\mathcal{L}$-theory $S$ as proving the usual axiom schemata for identity:
\begin{align}
 S&\vdash \forall_\sigma x\, (x=_\sigma x)&&&&&
    \tag{$=$ - Intro}\label{id_intro}\\
  S&\vdash \forall\bar x\forall\bar y\,(\varphi(\bar x)\land (\bar x=_{\bar\sigma} \bar y)
          \rightarrow \varphi(\bar y))&&&&&
    \tag{$=$ - Elim}\label{id_elim}
\end{align} for all variables $x:\sigma$, $\bar x,\bar y:\bar \sigma$, and $\mathcal{L}$-formulae $\varphi:\bar \sigma.$
Now, a piecewise interpretation $f:S\to T$ takes the identity $=_\sigma$ on $\sigma\in\mathcal{S}_1$ to a family $(=_\sigma)_f:P_\sigma\times P_\sigma\to\mathcal{L}_{2}$ with $(=_\sigma)_f^{k,l}:\sigma_f^k\times\sigma_f^l$. Since $f$ preserves theoremhood and the identity axioms are theorems of $S$, the translations of identity introduction and identity elimination are theorems of $T$: 
\begin{align}
 T&\vdash \bigwedge_{k\in P_\sigma^f}\forall x_f^k\, (D_x^k \to (x=_\sigma x)^{k}_f)&&&&&
    \tag{$=_f$ - Intro}\label{f-identity-intro}\\
  T&\vdash \bigwedge_{\bar k\in P_{\bar{\sigma}}^f} \bigwedge_{\bar l\in P_{\bar{\sigma}}^f}(\forall \bar x_f^{\bar k}\forall \bar y_f^{\bar l}\,(D_{\bar x}^{\bar k}\land D_{\bar y}^{\bar l}\land \varphi(\bar x)_f^{\bar k} \land (\bar x=_{\bar\sigma} \bar y)^{\bar k\times\bar l}_f\rightarrow \varphi(\bar y)_f^{\bar l}))&&&&&
    \tag{$=_f$ - Elim}\label{eq: f-identity-elim}
\end{align} 
where we set $D_{\bar x}^{\bar k}:=D_{\bar\sigma}^{\bar k}(\bar x_f^{\bar k})$ to facilitate readability. In particular, this implies, for all $k,l,m\in P_\sigma^f$:
\begin{align}
  T\vdash\; & \forall x_f^k\,(D_{x}^k\to (x=_\sigma x)_f^{k}) 
    \tag{Reflexivity}\label{Reflexivity}\\
  T\vdash\; &\forall x_f^k\forall y_f^l\, (D_{x,y}^{k,l}\to\left( (x=_\sigma y)_f^{k,l}\to (y=_\sigma x)_f^{l,k}\right))\tag{Symmetry}\label{Symmetry}\\
  T\vdash\; &\forall x_f^k\forall y_f^l\forall z_f^m\, (D_{x,y,z}^{k,l,m}\to \left((x=_\sigma y)_f^{k,l}
    \land (y=_\sigma z)_f^{l,m}\to (x=_\sigma z)_f^{k,m}\right))
    \tag{Transitivity}\label{Transitivity}
\end{align}
Reflexivity follows from the translation of identity introduction; in its presence, symmetry and transitivity are jointly equivalent to the translation of  identity elimination for identity. Semantically, the family $(=_\sigma)_f$ defines an equivalence relation on a disjoint union of extensions of sequences of sorts, as we will see. The summands of this disjoint union will be given by the extensions of the domain formulae of the interpretation. For these, identity introduction implies that $S\vdash \exists x (x=_\sigma x)$ for all sorts $\sigma\in \mathcal{S}_1$, and the translation of this sentence follows from $T$. This translation can easily be seen to be equivalent to the following, stating joint non-emptiness of the domain formulae for $f$ for each $\sigma$:
 \begin{equation}
        T\vdash\bigvee_{k\in P_\sigma^f} \exists x_{f}^k D_\sigma^k(x_f^k)\tag{Domain}\label{eq:domain}
    \end{equation}
\medskip
\noindent
We now know how piecewise interpretations treat identity. It is worth singling out those translations which behave like interpretations in this respect:
\begin{definition}[Piecewise Translation Respecting Identity]
    Let $f:\mathcal{L}_1\to\mathcal{L}_2$ be a piecewise translation. We say that $f$ respects identity with respect to an $\mathcal{L}_2$-theory $T$ iff $f$ and $T$ satisfy the non-emptiness of domains condition \eqref{eq:domain} and the identity conditions \eqref{f-identity-intro} and \eqref{eq: f-identity-elim}. 
\end{definition}
Note that a translation respects identity with respect to a theory $T$ if it satisfies the domain condition \eqref{eq:domain}, \eqref{f-identity-intro} and \eqref{eq: f-identity-elim} for atomic formulae, by induction. It follows that if $f\equiv_Tg$ and $f$ respects identity with respect to $T$, then so does $g$. 

\medskip\noindent
Such translations are the ones for which we will develop a theory of semantic duals, from which the following preservation of provability result will easily follow, generalizing the result for ordinary interpretations, proved syntactically as Theorem 5.4.7 in \citet{halvorson2019logic}: 
\begin{theorem}[Preservation of Provability]\label{thm:preserv_of_prov_act}
     Let $f:\mathcal{L}_1\to\mathcal{L}_2$ be a piecewise translation which respects identity with respect to an $\mathcal{L}_2$-theory $T$. For a set of $\mathcal{L}_1$-sentences $\Delta$, denote its image under $f$ by $\Delta_f:=\{\delta_f\mid\delta\in \Delta\}$, and let $\varphi$ be an $\mathcal{L}_1$-sentence. Then:
     \begin{equation}
         \Delta\vdash\varphi  \Rightarrow T\cup \Delta_f\vdash\varphi_f
     \end{equation}
\end{theorem}
\noindent

In particular, this entails that one can test whether a translation is an interpretation using a set of axioms for the source theory in question:
\begin{corollary}[Syntactic Criterion for Interpretationhood]\label{cor:syntactic criterion}
    Let $S$ be an $\mathcal{L}_1$-theory with a set of axioms $\Delta$ and $T$ an $\mathcal{L}_2$-theory. Then a piecewise translation $f:\mathcal{L}_1\to\mathcal{L}_2$ is an interpretation of $S$ in $T$ iff $f$ respects identity with respect to $T$ and $T\vdash\delta_f$ for each $\delta\in\Delta$.  
\end{corollary}
\noindent
In particular, this tells us that a translation $f:\mathcal{L}_1\to\mathcal{L}_2$ which respects identity with respect to an $\mathcal{L}_2$-theory $T$ can be seen as an interpretation of the theory of identity $L(\Sigma_1)$ in the signature $\Sigma_1$, axiomatized by \eqref{id_intro} and \eqref{id_elim}, into $T$. 
\subsection{Two Basic Examples: \texorpdfstring{$n$}{n} Elements and Projective Completions}
We are now ready for two concrete examples of piecewise interpretations. The first elementary example shows that piecewise interpretations induce a coarser interpretability relation than ordinary interpretations.
\begin{example}[$n$-Element Theories]\label{ex:n-element}
Let $n \geq 1$ be a natural number and let $\mathrm{E}_n$ be the theory in a signature with a single sort $\sigma$ and no relation symbols, except for the identity, saying that there are exactly $n$ elements of sort $\sigma$. For $n \geq 2$, there is no ordinary interpretation $\mathrm{E}_n \to \mathrm{E}_1$, as such an interpretation would need to define a quotient structure with $n$ elements inside a structure with a single element.

\medskip\noindent
There is, however, a piecewise interpretation $f_n : \mathrm{E}_n \to
\mathrm{E}_1$. Set $P_\sigma := [n]$ and, for $k,l\in [n]$,
$\sigma_{f_n}^k := \varepsilon$ and $D_\sigma^{f_n,k} := \top$, with identity being translated as
\begin{equation} 
(x=_\sigma y)_{f_n}^{k,l}:=
    \begin{cases}\top & \text{if }k=l,\\ \bot & \text{if }k\neq l.\end{cases}
\end{equation}
It is easy to check that this is indeed an interpretation by noting that it respects identity and computing the translation of the single axiom of $\mathrm{E}_n$. Indeed, $f_n$ can be used to provide an interpretation of $\mathrm{E}_n$ into any theory.
\end{example}
In the following example, we capture the construction of the projective completion of an affine plane as a piecewise interpretation, following Hartshorne's presentation of the geometric situation \citep{hartshorne1967foundations}.\footnote{We thank Albert Visser for suggesting this example to us.}  
\begin{example}[Projective Completion of an Affine Plane]\label{ex:projective_completion}
Consider the two-sorted signature $\Sigma = (\{\mathsf{p}, \mathsf{l}\},
\{I\})$ with a single non-identity relation symbol $I : \mathsf{p} \times
\mathsf{l}$ for incidence. We abbreviate parallelism by
$\ell_1 \parallel \ell_2 :=\ell_1 =_\mathsf{l} \ell_2 \lor
\neg\exists p\,(I(p,\ell_1) \land I(p,\ell_2))$.

\medskip
\noindent
The \emph{theory of affine planes} $T_{\mathrm{aff}}$ over $\Sigma$ has axioms stating that any two distinct points lie on a unique line, that for every point $p$ and line $\ell$, there is a unique line through $p$ parallel to $\ell$ and that there exist three non-collinear points. One can show that $T_{\mathrm{aff}}$ proves parallelism to be an equivalence relation. The \emph{theory of projective planes} $T_{\mathrm{proj}}$ over the same signature exchanges the axiom about unique parallels with an axiom saying that any two lines meet, and adds that every line contains at least three points.

\medskip
\noindent
We construct a piecewise interpretation $f : T_{\mathrm{proj}} \to T_{\mathrm{aff}}$.
The sort function assigns two pieces to each sort:
\begin{align*}
    P_\mathsf{p} &= \{1, 2\}, &
    \mathsf{p}_f^1 &= (\mathsf{p}), &
    \mathsf{p}_f^2 &= (\mathsf{l}), \\
    P_\mathsf{l} &= \{1, 2\}, &
    \mathsf{l}_f^1 &= (\mathsf{l}), &
    \mathsf{l}_f^2 &= \varepsilon.
\end{align*}
Piece~$1$ of $\mathsf{p}$ represents the affine points, and piece~$2$
represents the points at infinity, each given by a line up to parallelism.
Piece~$1$ of $\mathsf{l}$ represents the affine lines, while piece~$2$
represents the line at infinity.

\medskip
\noindent
All domain formulae give the full domain of the respective sorts:
$D_\mathsf{p}^1 = \top_\mathsf{p},\quad
D_\mathsf{p}^2 = \top_\mathsf{l},\quad
D_\mathsf{l}^1 = \top_\mathsf{l},\quad
D_\mathsf{l}^2 = \top_\varepsilon = \top$.

\medskip
\noindent
The translation of identity on $\mathsf{p}$ determines when two projective
points (possibly from different pieces) are identified:
\begin{alignat*}{2}
    (=_\mathsf{p})_f^{1,1}(x_1, x_2) &= (x_1 =_\mathsf{p} x_2), &\qquad
    (=_\mathsf{p})_f^{1,2} &= \bot_{\mathsf{p}\times\mathsf{l}}, \\
    (=_\mathsf{p})_f^{2,1} &= \bot_{\mathsf{l}\times\mathsf{p}}, &\qquad
    (=_\mathsf{p})_f^{2,2}(\ell_1, \ell_2) &= (\ell_1 \parallel \ell_2).
\end{alignat*}
That is, affine points are identified via ordinary equality, points at infinity are identified when their representing lines are parallel, and no affine point is ever identified with a point at infinity. The required reflexivity, symmetry, and transitivity conditions follow from the fact that $T_{\mathrm{aff}}$ proves parallelism to be an equivalence relation.

\medskip
\noindent
The identity on $\mathsf{l}$ is translated as follows:
\begin{alignat*}{2}
    (=_\mathsf{l})_f^{1,1}(\ell_1, \ell_2) &= (\ell_1 =_\mathsf{l} \ell_2), &\qquad
    (=_\mathsf{l})_f^{1,2} &= \bot_{\mathsf{l}}, \\
    (=_\mathsf{l})_f^{2,1} &= \bot_{\mathsf{l}}, &\qquad
    (=_\mathsf{l})_f^{2,2} &= \top.
\end{alignat*}
\noindent
Note that $(=_\mathsf{l})_f^{1,2}$ and $(=_\mathsf{l})_f^{2,1}$ have arity
$(\mathsf{l})$, while $(=_\mathsf{l})_f^{2,2}$ has arity $\varepsilon$ and is hence
a sentence; it expresses that the line at infinity is equal to itself.
Finally, the incidence relation $I : \mathsf{p} \times \mathsf{l}$ is
translated on the set of pieces $P_{\mathsf{p}\times\mathsf{l}} = \{1,2\}
\times \{1,2\}$:
\begin{alignat*}{2}
    I_f^{1,1}(p, \ell) &= I(p, \ell), &\qquad
    I_f^{1,2} &= \bot_{\mathsf{p}}, \\
    I_f^{2,1}(a, \ell) &= (a \parallel \ell), &\qquad
    I_f^{2,2} &= \top_{\mathsf{l}}.
\end{alignat*}
In prose: an affine point lies on an affine line according to ordinary incidence; no affine point lies on the line at infinity; a point at infinity $[\ell_1]_\|$ lies on an affine line $\ell_2$ iff $\ell_1 \parallel \ell_2$; and every point at infinity lies on the line at infinity. 

\medskip\noindent
To see that this translation respects identity, note that it clearly satisfies the domain condition and \eqref{f-identity-intro} and \eqref{eq: f-identity-elim} for the identities. For \eqref{eq: f-identity-elim} for incidence, the only interesting case is for $I_f^{2,1}$, where it follows from symmetry and transitivity of parallelism. To see that it is indeed an interpretation into $T_{\mathrm{aff}}$, one could now go through the axioms, or, instead, use the semantic means introduced in the next subsection.
\end{example}
\noindent
Three notes on this example: Firstly, one could also produce somewhat more involved interpretations $T_{\mathrm{proj}}\to T_{\mathrm{aff}}$ without pieces. Secondly, nothing hinges on allowing the empty sequence as the value of the sort function here; without it, the interpretations would just look a bit clunkier. Thirdly, there is also an interpretation $g:T_{\mathrm{aff}}\to T_{\mathrm{proj}}$, which, however, uses parameters to single out a line.
\subsection{The Semantic Dual of a Piecewise Interpretation}
In this subsection, we will see that each piecewise interpretation has an induced semantic dual: if $f:S\to T$ is a piecewise interpretation, there is a corresponding functor $f^*:\Mod(T)\to \Mod(S)$.\footnote{In Meadows' \citeyearpar{meadows2024beyond} terminology, the dual is a mod-functor compatible with the original translation.} 
The right morphisms between models of classical theories for our purposes will be the elementary embeddings between them. 
\begin{definition}[Category of Models of a Theory]
    Let $T$ be a classical first-order theory. Then $\Mod(T)$ denotes the category with objects models $M$ of $T$ and morphisms elementary embeddings between models.\footnote{A small remark on size: this category is only locally small, but by downward Löwenheim--Skolem, restricting to models with carriers of cardinality up to $\kappa$ would allow one to retain completeness while obtaining an essentially small category. For a more sophisticated discussion of size issues and categories of models, we direct the reader to \citet[Section~2.2.2]{meadows2024beyond}.} 
\end{definition}
At this stage, we need to be a bit more precise about our model theory for many-sorted first-order logic. To define the semantic dual, we will introduce a notion of extension not only for atomic, but also for complex formulae. Let us fix some standard definitions:
\begin{definition}[$\Sigma$-Structures and Assignments]
    A $\Sigma$-structure $M$ assigns to each sort $\sigma\in \mathcal{S}$ a non-empty set $M(\sigma)$ and to each relation symbol $R:\bar\sigma=\sigma_1\times\cdots\times\sigma_n$ a subset $M(R)$ of $M(\bar\sigma):=M(\sigma_1)\times\cdots\times M(\sigma_n)$. Accordingly, the verum $\top:\varepsilon$ is assigned the subset $\{\varepsilon\}\subseteq M(\varepsilon)=\{\varepsilon\}$ and the falsum $\bot:\varepsilon$ is assigned the subset $\varnothing\subset M(\varepsilon)$. We will assume that the identity $=_\sigma$ is assigned the extension $M(=_\sigma)=\{(x,y)\in M(\sigma)\times M(\sigma)\mid x=y\}$. 

    \medskip\noindent
    Given a $\Sigma$-structure $M$, a variable assignment $\alpha:\mathcal{V}_\Sigma\to M$ for $M$ is a function assigning to each variable $x:\sigma$ an element $\alpha(x)\in M(\sigma)$. The set of all variable assignments for $M$ is denoted by $\mathrm{A}(M)$. For a sequence $\bar x:\sigma_1\times\cdots\times\sigma_n$, we extend any given assignment $\alpha$ by $\alpha(\bar x)=(\alpha(x_1),\ldots,\alpha(x_n))$. For a duplicate-free $\bar x:\bar\sigma$, a sequence of elements $\bar a\in M(\bar\sigma)$, and an assignment $\alpha$, we define a variant $\alpha[\bar x\mapsto\bar a]$ by $\alpha[\bar x\mapsto\bar a](x_i)=a_i$ and $\alpha$ otherwise.  
\end{definition}
Satisfaction of formulae with respect to a structure and an assignment, denoted by $M,\alpha\vDash\varphi$, is then defined recursively in the standard Tarskian fashion. We get the following expected lemma on substitution and satisfaction, whose proof is an easy induction on complexity: 
\begin{lemma}[Substitution and Satisfaction]\label{lemma:subst_and_sat}
    Let $M$ be a $\Sigma$-structure, $\varphi:\bar \sigma$ be a $\Sigma$-formula with canonical context $\bar x$, and $\bar y:\bar \sigma$. Then, for any variable assignment $\alpha$ for $M$, the following holds:
    \begin{equation}
        M,\alpha\vDash\varphi(\bar y)\;\iff\;M,\alpha[\bar x\mapsto\alpha(\bar y)]\vDash \varphi
    \end{equation}
\end{lemma}

\noindent
Satisfaction of sentences does not depend on the variable assignment chosen, and we say that a $\Sigma$-structure is a model of a theory $T$, denoted by $M\vDash T$ iff $M\vDash \varphi$ for all $\varphi\in T$. Using assignments, we can define a notion of extension for general formulae:
\begin{definition}[Extension of a Formula in a Structure]
    Given a $\Sigma$-structure $M$ and a $\Sigma$-formula $\varphi:\bar \sigma$ with canonical context $\bar x:\bar\sigma$, the extension of $\varphi$ in $M$ is denoted by $M(\varphi)$ and defined by
    \begin{equation}
        M(\varphi)=\{\alpha(\bar x)\mid \alpha\in \mathrm{A}(M) \text{ s.t. }M,\alpha\vDash\varphi\}\subseteq M(\bar \sigma). 
    \end{equation}
\end{definition}
Intuitively, the extension of a formula consists of all those elements of the given structure satisfying the formula, generalizing from the case of relation symbols. Note that the extension of a sentence is $M(\varepsilon)$ if it is true in $M$ and $\varnothing$ otherwise. With the notion of extension at hand, we can define the semantic dual of a piecewise translation: 
\begin{definition}[Semantic Dual of a Piecewise Translation]\label{def:semantic_dual}
    Let $f:\mathcal{L}_1\to \mathcal{L}_2$ be a piecewise translation respecting identity with respect to an $\mathcal{L}_2$-theory $T$. Given a model $M$ of $T$, we define a $\Sigma_1$-structure $f^*M$ as follows: For a sort $\sigma$ of $\Sigma_1$, the extension of $\sigma$ in $f^*M$ is 
    \begin{equation}
        f^*M(\sigma):=\Bigl(\bigsqcup_{k\in P_\sigma ^f} M(D_\sigma^{f,k})\Bigr)\Big/\sim_\sigma^f   
    \end{equation} where $(k,\bar a)\sim_\sigma^f (l,\bar b)$ iff $\bar a\times\bar b\in M((=_\sigma)_f^{k,l})$. Note that $\sim_\sigma^f$ is an equivalence relation on a non-empty set by the considerations in \cref{sec:identity}, and denote by $\pi_\sigma^f:\sqcup_{k\in P_\sigma^f}M(D_\sigma^{f,k})\to f^*M(\sigma)$ the associated projection. To extend this to sequences of sorts, note that for $\bar\sigma=(\sigma_1,\ldots,\sigma_n)$, 
    \begin{equation}\label{eq:domain_sequence_extension}
        M(D^{f,\bar k}_{\bar{\sigma}}) = \{\bar a_1\times\cdots\times\bar a_n\mid \bar a_i\in M(D_{\sigma_i}^{f,k_i}) \text{ for }1\leq i\leq n\},
    \end{equation} whence
    \begin{align} \label{eq:projection_map}
\pi_{\bar \sigma}^f:\bigsqcup_{\bar k\in P^f_{\bar\sigma}} M(D_{\bar\sigma}^{f,\bar k})&\to \prod_{i=1}^{|\bar\sigma|}f^*M(\sigma_i)\\
(\bar k,\bar a)&\mapsto (\pi_{\sigma_1} ^f(k_1,\bar a_1),\ldots,\pi_{\sigma_n}^f(k_n,\bar a_n))
\end{align}
 is well-defined and surjective, where $\bar a=\bar a_1\times\cdots\times\bar a_n$ is such that $| \bar a_i|=| (\sigma_i)_f^{k_i}|$ for $1\leq i\leq n$. Given a relation symbol $R:\bar \sigma$ of $\Sigma_1$, its extension is defined by
     \begin{equation}
        f^*M(R):=\pi_{\bar \sigma}^f\Bigl(\bigsqcup_{\bar k\in P_{\bar\sigma} ^f} (M(R_f^{\bar k})\cap M(D_{\bar \sigma}^{f,\bar k}))\Bigr)\label{translation_extension_atomic}
    \end{equation}
\end{definition}
The structures defined via the duals of piecewise interpretations are not interpretable structures in the usual sense, as these only allow for the taking of products, subsets, and quotients under equivalence relations, but not for disjoint unions. However, the formation of disjoint unions does not seem to be more problematic than these in any obvious way. In particular, it does not rely on problematic infinitary or nonconstructive resources.

\medskip\noindent
Schematically, equation \eqref{translation_extension_atomic}, suitably understood, also applies to general formulae. To make this precise, we need a notion of variable assignment for a translation for a structure:
\begin{definition}[$f$-Assignment for a Structure $M$]\label{def:f_assignment}
    Let $f:\mathcal{L}_1\to \mathcal{L}_2$ be a piecewise translation and let $M$ be a $\Sigma_2$-structure. An $f$-assignment for $M$ is a tuple $\langle\gamma, \alpha\rangle$, where $\gamma$ is a function assigning to each $\Sigma_1$-variable $x:\sigma$ a piece $\gamma(x)\in P_\sigma^f$, and $\alpha$ is an $M$-assignment satisfying $\alpha(x_f^{\gamma(x)})\in M(D_\sigma^{f,\gamma(x)})$ for every $\Sigma_1$-variable $x:\sigma$. If $f$ respects identity with respect to a $\Sigma_2$-theory $T$ and $M\vDash T$, we define an $f^*M$-assignment $\pi\langle\gamma,\alpha\rangle$ by $x\mapsto\pi_\sigma^f(\gamma(x),\alpha(x_f^{\gamma(x)}))\in f^*M(\sigma)$.
\end{definition}
Three small notes on this definition: Firstly, if a translation $f:\mathcal{L}_1\to \mathcal{L}_2$ respects identity with respect to an $\mathcal{L}_2$-theory $T$ and $M\vDash T$, then there is an $f$-assignment for $M$. Secondly, surjectivity of $\pi_\sigma^f$ implies that every assignment for $f^*M$ is of the form $\pi\langle\gamma,\alpha\rangle$ for an $f$-assignment $\langle\gamma,\alpha\rangle$ for $M$. Thirdly, given a duplicate-free $\bar x:\bar \sigma$, a piece $\bar k\in P_{\bar \sigma}^f$, and a sequence of elements $\bar a\in M(D_{\bar \sigma}^{f,\bar k})$, there is always an $f$-assignment $\langle\gamma,\alpha\rangle$ such that $\gamma(\bar x)=\bar k$ and $\alpha(\bar x_f^{\bar k})=\bar a$, since we can extend $\gamma$ and $\alpha$ arbitrarily. 

\medskip\noindent
More substantially, we can now formulate a proposition on how satisfaction and extension in the semantic dual are related to satisfaction and extension in the original model. From a technical point of view, this is perhaps the most important proposition in this paper. It also justifies the name \emph{semantic dual} for $f^*$:
\begin{proposition}[Satisfaction of Translated Formulae]\label{prop:satisfaction_under_translation}
    Let $f:\mathcal{L}_1\to \mathcal{L}_2$ be a piecewise translation respecting identity with respect to an $\mathcal{L}_2$-theory $T$, and let $M$ be a model of $T$. Then, for every $\mathcal{L}_1$-formula $\varphi$ with canonical context $\bar x:\bar \sigma$ and every $f$-assignment $\langle\gamma,\alpha\rangle$ for $M$,
    \begin{equation}\label{eq:assignment_translation}
        M,\alpha\vDash\varphi_f^{\gamma(\bar x)}\;\iff\;f^*M, \pi\langle\gamma,\alpha\rangle\vDash\varphi.
    \end{equation}
In particular, for every $\mathcal{L}_1$-sentence $\varphi$,
\begin{equation}
    M\vDash\varphi_f\,\iff\;f^*M\vDash\varphi.
\end{equation}
\end{proposition}
\begin{proof}
    We proceed by induction on complexity of $\varphi$. A preliminary remark regarding translations in context is expedient: given a $\Sigma_2$-structure $M$ and an $f$-assignment $\langle\gamma,\alpha\rangle$ for $M$, let $\bar y:\bar \tau$ be a context for $\varphi:\bar \sigma$ and denote by $\bar x$ the canonical context of $\varphi$. Then, by the context clause for translations, 
    \begin{equation}\label{eq:gamma_context}
        (\bar y\mid\varphi)_f^{\gamma(\bar y)}=\varphi_f^{\gamma(\bar x)}.
    \end{equation}
    
    \medskip\noindent
    The cases for the truth values and negation are immediate. The case for the binary connectives follows from the Tarskian clauses, the binary clauses for translations, and the observation \eqref{eq:gamma_context}.

    \medskip\noindent
    Turning to atomic formulae, let $\varphi=R(\bar z)$ with $\bar z=(z_1,\ldots,z_m):\bar\tau$ and canonical
context $\bar x=(x_1,\ldots,x_n)=\bar z^\circ:\bar \sigma$, and let $d:[m]\to [n]$ be such that $x_{d(i)}=z_i$, as in the atomic clause for translations. Then the atomic clause for translation together with $(\gamma (\bar x))^d=\gamma(\bar x^d)=\gamma(\bar z)$ imply $\varphi_f^{\gamma(\bar x)}=R_f^{\gamma(\bar z)}(\bar z_f^{\gamma(\bar z)})$, whence
\cref{lemma:subst_and_sat} yields
\begin{equation}\label{eq:atomic_left}
    M,\alpha\vDash\varphi_f^{\gamma(\bar x)}\iff
    \alpha(\bar z_f^{\gamma(\bar z)})\in M(R_f^{\gamma(\bar z)}).
\end{equation}
On the right-hand side, by definition, we get
\begin{equation}
f^*M, \pi\langle\gamma,\alpha\rangle\vDash R(\bar z) \iff
\bigl(\pi\langle\gamma,\alpha\rangle(z_1),\ldots,\pi\langle\gamma,\alpha\rangle(z_m)\bigr)\in f^*M(R).
\end{equation}

\noindent
Now, from left to right, $\alpha(\bar z_f^{\gamma(\bar z)})\in M(R_f^{\gamma(\bar z)})$ by assumption and $\alpha(\bar z_f^{\gamma(\bar z)})\in M(D_{\bar\tau}^{f,\gamma\bar z})$ by the definition of an $f$-assignment,
together with \eqref{eq:domain_sequence_extension}. Therefore, $\pi\langle\gamma,\alpha\rangle(\bar z)=\pi_{\bar \tau}^f(\gamma(\bar z),\alpha(\bar z_f^{\gamma(\bar z)}))\in f^*M(R)$, by the definition of the semantic dual \eqref{translation_extension_atomic}. From right to left, choose
$\bar l\in P_{\bar\tau}^f$ and $\bar a=\bar a_1\times\cdots\times\bar a_m\in M(R_f^{\bar l})\cap M(D_{\bar\tau}^{f,\bar l})$ such that
$\pi_{\bar\tau}^f(\bar l,\bar a)=\bigl(\pi\langle\gamma,\alpha\rangle(z_1),\ldots,\pi\langle\gamma,\alpha\rangle(z_m)\bigr)$; by construction,
$\bar a_i\times\alpha(z_{i,f}^{\gamma (z_i)})\in M\bigl((=_{\tau_i})_f^{l_i,\gamma (z_i)}\bigr)$, so instantiating \eqref{eq: f-identity-elim} for $R$ at $\bar a$ and $\alpha(\bar z_f^{\gamma(\bar z)})$ yields the right-hand side of \eqref{eq:atomic_left}.

\medskip\noindent
For the universal quantifier, let $\varphi=\forall_\tau y\,\psi$ with canonical context $\bar x$. Applying the clause for the universal quantifier and reducing contexts, we get
\begin{equation}\label{eq:quantifier_gamma}
    \varphi_f^{\gamma(\bar x)}=\bigwedge_{h\in P_\tau^f}\forall y_f^h
    \bigl(D_\tau^{f,h}(y_f^h)\to\psi_f^{\gamma[y\mapsto h](\bar x(\psi))}\bigr).
\end{equation}
Therefore,
\begin{equation}
    M,\alpha\vDash\varphi_f^{\gamma(\bar x)}\iff \forall h\in P_\tau^f\forall \bar a\in M(D_\tau^{f,h}): M,\alpha[y_f^h\mapsto\bar a]\vDash \psi_f^{\gamma[y\mapsto h](\bar x(\psi))}.
\end{equation}
Applying the induction hypothesis with $f$-assignment $\langle \gamma',\alpha'\rangle:= \langle\gamma[y\mapsto h], \alpha[y_f^h\mapsto\bar a]\rangle$, together with surjectivity of $\pi_\tau^f$, we see that this is equivalent to 
\begin{equation}
    \forall b\in f^*M(\tau): f^*M,\pi\langle\gamma,\alpha\rangle[y\mapsto b]\vDash\psi
\end{equation} which is equivalent to the desired right-hand side by the Tarskian clause for the quantifier.
\end{proof}
\noindent
We can now fulfil the promise to prove \cref{thm:preserv_of_prov_act}, on preservation of provability, by semantic means: 
\begin{proof}[Proof of \cref{thm:preserv_of_prov_act}]
    Let $f:\mathcal{L}_1\to \mathcal{L}_2$ be a piecewise translation respecting identity with respect to an $\mathcal{L}_2$-theory $T$. Moreover, let $\Delta$ be a set of $\mathcal{L}_1$-sentences, and let $\varphi$ be an $\mathcal{L}_1$-sentence. Assume that $\Delta\vdash \varphi$. Let $M\vDash T$ such that $M\vDash\Delta_f$. Then, by \cref{prop:satisfaction_under_translation}, $f^*M\vDash\Delta$, whence, by soundness, $f^*M\vDash \varphi$. By \cref{prop:satisfaction_under_translation} again, $M\vDash \varphi_f$, whence, since $M\vDash T$ was arbitrary, completeness allows us to conclude that $T\cup\Delta_f\vdash\varphi_f$.
\end{proof}
From the observations immediately following \cref{def:f_assignment}, we get an essentially equivalent reformulation of \cref{prop:satisfaction_under_translation} in terms of extensions:
\begin{corollary}[Extensions of Translated Formulae]\label{cor:extensions_elementwise}
   In the setting of \cref{prop:satisfaction_under_translation}, the following holds: for every $\mathcal{L}_1$-formula $\varphi:\bar\sigma$, every piece $\bar k\in P_{\bar\sigma}^f$ and every tuple $\bar a\in M(D_{\bar\sigma}^{f,\bar k})$
    \begin{equation}\label{eq:extension_elementwise}
        \bar a\in M(\varphi_f^{\bar k})\;\iff\;\pi_{\bar\sigma}^f(\bar k,\bar a)\in f^*M(\varphi).
    \end{equation}
    In particular, taking images under $\pi^f_{\bar\sigma}$,
    \begin{equation}\label{eq:extension_projection}
        f^*M(\varphi)=\pi_{\bar\sigma}^f\Bigl(\bigsqcup_{\bar k\in P_{\bar \sigma}^f}
            M(\varphi_f^{\bar k}) \cap M(D_{\bar \sigma}^{f,\bar k})\Bigr).
    \end{equation}
\end{corollary}

\noindent
Putting this together with adequacy of the semantics allows us to state a semantic characterization of interpretationhood, paving the way for us to understand $f^*$ as a functor:
\begin{corollary}[Characterization of Interpretations via Model Map]\label{cor:model_map}
    Let $f:\mathcal{L}_1\to \mathcal{L}_2$ be a piecewise translation respecting identity with respect to an $\mathcal{L}_2$-theory $T$. The translation $f$ is a piecewise interpretation $f:S\to T$ iff for any model $M$ of $T$, the $\mathcal{L}_1$-structure $f^*M$ is a model of $S$. In particular, for a piecewise interpretation, $f^*$ induces a function $f^*:\Mod(T)\to \Mod(S)$.
\end{corollary}
\begin{proof}
    By adequacy of the semantics, $f$ is an interpretation $f:S\to T$ iff, for every $\mathcal{L}_1$-sentence $\varphi$ made true in every model of $S$, $\varphi_f$ is made true in every model of $T$. By \cref{prop:satisfaction_under_translation}, this is the case iff for every such sentence $\varphi$ and model $M$ of $T$, $f^*M\vDash\varphi$. But the latter just says that $f^*M$ is a model of $S$. 
\end{proof}
\noindent
Moreover, these elementary properties of the semantic dual are enough to also characterize provable equivalence semantically, thereby establishing invariance of interpretationhood under it.
\begin{proposition}[Provable Equivalence and Model Map]\label{prop:prov_equiv_char}
    Let $f,g:\mathcal{L}_1\to\mathcal{L}_2$ be piecewise translations respecting identity with respect to an $\mathcal{L}_2$-theory $T$, and assume that $\sigma_f=\sigma_g$ for every
    $\sigma\in\mathcal{S}_1$. Then the following are equivalent.
    \begin{enumerate}[label=(\roman*)]
        \item $f\equiv_T g$.
        \item $f^*M=g^*M$ for every model $M$ of $T$.
        \item The domain functions of $f$ and $g$ are provably equivalent, and for every
        $\mathcal{L}_1$-formula $\varphi:\bar\sigma$, every $\bar k\in P_{\bar\sigma}^f$ and
        every $\bar u:\bar\sigma_f^{\bar k}$,
        \begin{equation}\label{eq:prov_equiv_formulae}
            T\vdash D_{\bar\sigma}^{f,\bar k}(\bar u)\to
            \bigl(\varphi_f^{\bar k}(\bar u)\leftrightarrow\varphi_g^{\bar k}(\bar u)\bigr).
        \end{equation}
    \end{enumerate}
    Therefore, if $f:S\to T$ and $f\equiv_Tg$, then $g:S\to T$.
\end{proposition}
\begin{proof}
(i) to (ii) follows by inspection of \cref{def:semantic_dual}, soundness and \cref{lemma:subst_and_sat}.
For (ii) to (iii), let $M\vDash T$ and $\sigma\in\mathcal S_1$. The elements of $f^*M(\sigma)$ are the $\sim^f_\sigma$-classes, which partition $\bigsqcup_{k\in P^f_\sigma}M(D^{f,k}_\sigma)$; hence $\bigcup f^*M(\sigma)=\bigsqcup_{k\in P^f_\sigma}M(D^{f,k}_\sigma)$, and likewise for $g$. Since $P^f_\sigma=P^g_\sigma$, the equality $f^*M(\sigma)=g^*M(\sigma)$ thus yields $\bigsqcup_{k}M(D^{f,k}_\sigma)=\bigsqcup_{k}M(D^{g,k}_\sigma)$, and by comparison of the fibres over each $k\in P^f_\sigma$, $M(D^{f,k}_\sigma)=M(D^{g,k}_\sigma)$. Since $M\vDash T$ was arbitrary, completeness yields the provable equivalence of the domain formulae. Moreover, $f^*M(\sigma)$ and $g^*M(\sigma)$ are the same partition of the same set, so $\sim^f_\sigma={\sim^g_\sigma}$ and $\pi^f_{\bar\sigma}=\pi^g_{\bar\sigma}$. Claim~\eqref{eq:prov_equiv_formulae} now follows from \eqref{eq:extension_elementwise} for $f$ and $g$, together with completeness.
(iii) to (i) merely requires specializing to atomic formulae. The last claim then follows from \cref{cor:model_map}.
\end{proof}
\noindent

\begin{example}[$n$-Element Theories, continued]\label{ex:n-element-cont}
    Consider the piecewise interpretation $f_n:\mathrm{E}_n\to \mathrm{E}_1$ from \cref{ex:n-element}. In a model $M$ of $\mathrm{E}_1$, the model $f_n^*M$ is constructed using the set $M(\top)=\{\varepsilon\}$, since $\sigma_{f_n}^k=\varepsilon$ and $D_\sigma^{f_n,k}=\top$. Set $c_k:=\pi_\sigma^{f_n}(k,\varepsilon)$. Then
    \begin{equation}
        f_n^*M(\sigma)=\{c_k\mid k\in [n]\} \text{ with } f_n^*M(=_\sigma)=\{(c_k,c_k)\mid k\in [n]\}.
    \end{equation}
    This is clearly a model of $\mathrm{E}_n$.
\end{example}
\begin{example}[Projective Completion, Continued]
\label{ex:projective_completion_dual}
Consider the piecewise translation $f$ from \cref{ex:projective_completion}, which respects identity with respect to $T_{\mathrm{aff}}$. Given a model $M\vDash T_{\mathrm{aff}}$ and going through the definition of the semantic dual, we see that we can identify the carriers of $f^*M$ with
\[
    f^*M(\mathsf p)
    = M(\mathsf p)\sqcup\bigl(M(\mathsf l)/{\parallel_M}\bigr),
    \qquad
    f^*M(\mathsf l)
    = M(\mathsf l)\sqcup\{\ell_\infty\},
\]
where $\parallel_M$ denotes parallelism in $M$. Here an affine point $p$ corresponds to
$\pi_{\mathsf p}^f(1,p)$, an equivalence class under parallelism $[a]_{\parallel_M}$ to $\pi_{\mathsf p}^f(2,a)$, an affine line $\ell$ to $\pi_{\mathsf l}^f(1,\ell)$, and $\ell_\infty$ to $\pi_{\mathsf l}^f(2,\varepsilon)$.

\medskip\noindent
Under these identifications, the interpreted identity relations
are ordinary equality, and incidence is given by
\begin{align*}
    f^*M(I)
    ={}&
    \bigl\{(p,\ell)\mid (p,\ell)\in M(I)\bigr\}\\
    &{}\cup
    \bigl\{([a]_{\parallel_M},\ell)
       \mid a,\ell\in M(\mathsf l),\ a\parallel_M\ell\bigr\}\\
    &{}\cup
    \bigl\{([a]_{\parallel_M},\ell_\infty)
       \mid a\in M(\mathsf l)\bigr\}.
\end{align*}
Thus, in particular, each affine line $\ell$ has a corresponding point at infinity $[\ell]_{\parallel_M}$, and the additional line $\ell_\infty$ contains exactly the points at infinity. This is the usual projective completion of $M$, which is a projective plane \citep{hartshorne1967foundations}. As a consequence, $f^*M\vDash T_{\mathrm{proj}}$ for every $M\vDash T_{\mathrm{aff}}$. \cref{cor:model_map} then entails that $f:T_{\mathrm{proj}}\to T_{\mathrm{aff}}$ is indeed a piecewise interpretation. 
\end{example}

\subsection{The Semantic Dual as a Functor between Categories of Models}
Given a piecewise interpretation $f:S\to T$, the function $f^*:\Mod (T)\to \Mod(S)$ can be seen as the object-part of a functor.\footnote{For ordinary interpretations, this is detailed clearly in \citet{halvorson2019logic}.}
The definition of $f^*$ on morphisms, that is, on elementary embeddings between models of the target theory, proceeds in the following fashion:  Recall that an elementary embedding between $\Sigma$-structures $j:M\to N$ is a function from the carrier of $M$ to that of $N$, that is, a family of functions $(j_\sigma:M(\sigma)\to N(\sigma))_{\sigma\in \mathcal{S}}$,  such that, if $\varphi:\bar \sigma=\sigma_1\times\cdots\times\sigma_n$ is an $\mathcal{L}_\Sigma$-formula and $\bar a\in M(\bar\sigma)$, 
\begin{equation}
    (a_1,\ldots,a_n)\in M(\varphi) \Leftrightarrow (j_{\sigma_1}(a_1),\ldots,j_{\sigma_n}(a_n))\in N(\varphi). \label{elementary-embedding_atomic}
\end{equation}
We will suppress these subscripts from now on, as the respective sort will always be clear from context.
\begin{definition}[Semantic Dual on Elementary Embeddings]\label{def:semantic_dual_ee}
    Let $f:\mathcal{L}_1\to \mathcal{L}_2$ be a piecewise translation respecting identity with respect to $T$, and let $j:M\to N$ be an
    elementary embedding between models of $T$. Define
    $f^*j:f^*M\to f^*N$ by
    \begin{equation}\label{eq:dual-ee}
      f^*j\bigl(\pi_{M,\sigma}^f(k,\bar a)\bigr)
      \;:=\;
      \pi_{N,\sigma}^f(k,\,j(\bar a)).
    \end{equation}
    Well-definedness follows from the fact that $j$ is an elementary embedding.
\end{definition}
\begin{proposition}[Semantic Dual is a Functor]
    Let $f:S\to T$ be a piecewise interpretation. Then there is a corresponding functor $f^*:\Mod(T)\to\Mod(S)$, sending a model $M\vDash T$ to $f^*M\vDash S$ and an elementary embedding $j:M\to N$ to the elementary embedding $f^*j:f^*M\to f^*N$.
\end{proposition}
\begin{proof}
Given a model $M\vDash T$, \cref{cor:model_map} tells us that $f^*M$ is a model of $S$. To see that $f^*j$ is an elementary embedding if $j: M\to N$ is, let $\bar{a}\in\prod_i f^*M(\sigma_i)$ with $a_i=\pi_{\sigma_i}^f(k_i,\bar{b}_i)$, and set $\bar{k}=(k_1,\ldots,k_n)$.
Then for any $\mathcal{L}_S$-formula $\varphi:\bar{\sigma}$,
\cref{cor:extensions_elementwise} yields
\begin{align}
  \bar{a}\in f^*M(\varphi)
  &\iff \bar{b}_1\times\cdots\times\bar{b}_n
        \in M(\varphi_f^{\bar{k}})\cap M( D_{\bar{\sigma}}^{f,\bar{k}})\\
  &\iff j(\bar{b}_1)\times\cdots\times j(\bar{b}_n)
        \in N(\varphi_f^{\bar{k}})\cap N(D_{\bar{\sigma}}^{f,\bar{k}})\\
  &\iff f^*j(\bar{a})\in f^*N(\varphi),
\end{align}
where the middle step applies elementarity of $j$ and the outer steps
use \cref{cor:extensions_elementwise} for $M$ and $N$
respectively. 

\medskip\noindent
Functoriality is easy: $f^*j$ acts on
representatives by applying $j$ componentwise, so
$f^*(j_2 \circ j_1) = f^*j_2 \circ f^*j_1$ and
$f^*(\id_M) = \id_{f^*M}$.
\end{proof}
\noindent
Note that if $f,g:S\to T$ satisfy $f\equiv_Tg$, then $f^*=g^*$, by \cref{prop:prov_equiv_char} and \cref{def:semantic_dual_ee}.
\subsection{Morphisms between Piecewise Translations}
While the assignment $f\mapsto f^*$, sending an interpretation to its semantic dual, is clearly not injective, it becomes injective (on interpretations with the same sort functions) after descending to $\equiv_T$-equivalence classes by \cref{prop:prov_equiv_char}. However, provable equivalence does not capture other natural relations between the semantic duals of translations: for instance, one might be interested in naturally isomorphic duals or duals which are the same up to a reindexing of the sets of pieces, both of which provide other natural candidates for equivalence of translations.
 
\medskip
\noindent
In this subsection, we develop a notion of equivalence for translations which does not unduly depend on the sets of pieces used, on the basis of a syntactically expressible notion of morphism between piecewise translations. Semantically, such a morphism will correspond to a natural transformation between the semantic duals of the interpretations. 
\begin{definition}[$i$-Morphism between Piecewise Translations]\label{def:i-morphism}
Let $f,g:\mathcal{L}_1\to \mathcal{L}_2$ be two piecewise translations. A pre-$i$-morphism $\eta$ from $f$ to $g$ is a collection $(\eta_\sigma:P_\sigma^f\times P_\sigma^g\to \mathcal{L}_2)_{\sigma\in \mathcal{S}_1}$ of finite families of $\mathcal{L}_2$-formulae for each sort of $\mathcal{L}_1$ such that $\eta_\sigma^{k,l}:\sigma_f^k\times\sigma_g^l$. We extend to sequences $\bar\sigma\in \mathcal{S}_1^*$ and pieces $\bar k\in P_{\bar \sigma}^f$ and $\bar l\in P_{\bar \sigma}^g$ by choosing $\bar x=\bar x_1\times\cdots\times\bar x_n:\bar\sigma_f^{\bar k}$ and $\bar y=\bar y_1\times\cdots\times\bar y_n:\bar\sigma_g^{\bar l}$ and setting $\eta_{\bar \sigma}^{\bar k,\bar l}:=\top_{\bar x\times\bar y}\land\bigwedge_{i=1}^{|\bar \sigma|}\eta_{\sigma_i}^{k_i,l_i}(\bar x_i,\bar y_i)$. It is an $i$-morphism with respect to an $\mathcal{L}_2$-theory $T$, expressed by $\eta:f\Rightarrow g$, iff the theory $T$ proves the universal closures of the following formulae for all appropriate pieces, writing $\eta_{x,y}^{k,l}$ for $\eta_\sigma^{k,l}(x_f^k, y_g^l)$:
\begin{align}
     \eta_{x,y}^{k,l}&\to \left(D^{f,k}_{x} \land D^{g,l}_{y}\right)  \tag{Dom/Range}\label{Dom/Range} \\
     \left(\eta^{k,l}_{x,y} \land D^{f,k'}_{x'} \land D^{g,l'}_{y'} \land (x =_\sigma x')^{k,k'}_f \land (y =_\sigma y')^{l,l'}_g\right) &\to \eta^{k',l'}_{x',y'} \tag{Well-Def}\label{well-def}\\
     D_{x}^{f,k}&\to \bigvee_{l\in P_\sigma^g}\exists y_g^l\eta_{x,y}^{k,l} \tag{Joint-Exist}\label{Joint-Exist}\\
     \left(\eta^{k,l_1}_{x,y}\land \eta_{x,z}^{k,l_2}\right)&\to (y=z)_g^{l_1,l_2} \tag{Unique}\label{Unique}    
\end{align}
Moreover, for each formula $\varphi:\bar\sigma\in \mathcal{L}_1$, pair of pieces $\bar k\in P_{\bar \sigma}^f$  and $\bar l\in P_{\bar \sigma}^g$, the theory $T$ must prove
\begin{equation}
\eta_{\bar x, \bar y}^{\bar k,\bar l}\to(\varphi(\bar x)_f^{\bar k}\leftrightarrow \varphi(\bar y)_g^{\bar l}) \tag{Element-Embed}\label{Element-Embed}
\end{equation} where $\eta_{\bar x,\bar y}^{\bar k, \bar l}:=\eta_{\bar\sigma}^{\bar k,\bar l}(\bar x_f^{\bar k},\bar y_g^{\bar l})$. Here $\bar\sigma=\varepsilon$ is permitted: in that case the antecedent is the empty conjunction $\top$, and since $P_\varepsilon=\{\varepsilon\}$ and $\varphi_f^{\varepsilon}=\varphi_f$, instantiation yields $T\vdash\varphi_f\leftrightarrow\varphi_g$.
\end{definition}
For $i$-morphisms from a translation respecting identity, \eqref{Unique} already follows from \eqref{Element-Embed} applied to the identity, together with \eqref{Dom/Range} and \eqref{Reflexivity} for $f$, so one does not need to check it in these cases. Typically, we will encounter situations in which we have a (pre-)$i$-morphism $\eta:f\Rightarrow g$, where we antecedently know that $f$ is an interpretation $S\to T'$. In such cases, we will assume $T=T'$. 
\medskip
\noindent
Note that if we have translations $f\equiv_T f'$ and $g\equiv_Tg'$, then a collection of formulae defining a pre-$i$-morphism/an $i$-morphism from $f$ to $g$ also defines a pre-$i$-morphism/an $i$-morphism from $f'$ to $g'$. We will only ever consider (pre-)$i$-morphisms up to provable equivalence in a target theory. 
\begin{definition}[Provable Equivalence of $i$-Morphisms]\label{def:i-morphisms_prov_equiv}Two parallel pre-$i$-morphisms $\eta$ and $\zeta$ between translations $f$ and $g$ are provably equivalent with respect to a theory $T$, denoted by $\eta\equiv_T\zeta$, iff, for every pair of pieces $(k,l)\in P_\sigma^f\times P_\sigma^g$ and all $x,y:\sigma$, the following holds:
    \begin{equation}
        T\vdash \eta_{x,y}^{k,l}\leftrightarrow\zeta_{x,y}^{k,l}\label{provable equivalence for i-morphisms}
    \end{equation} 
\end{definition}
    Taking the sentential case of condition \eqref{Element-Embed} as sufficient for equivalence of translations, we obtain a weaker but still instrumentally useful and interesting notion of equivalence for translations:
\begin{definition}[Sentential Equivalence]\label{def:sent_equiv}
    Two parallel translations $f,g:\mathcal{L}_1\to \mathcal{L}_2$ are sententially equivalent with respect to an $\mathcal{L}_2$-theory $T_2$ iff the target theory $T_2$ proves $\varphi_f\leftrightarrow \varphi_g$ for any $\mathcal{L}_1$-sentence $\varphi$.
\end{definition} 
Sentential equivalence is clearly entailed by provable equivalence. The following proposition provides a useful tool for establishing interpretationhood of translations via $i$-morphisms, established by inspecting the sentential instances of \eqref{Element-Embed}.
\begin{proposition}[Criterion for Interpretationhood]\label{prop:criterion_interpretationhood}
    Let $f,g:\mathcal{L}_1\to \mathcal{L}_2$ be piecewise translations and $T_2$ an $\mathcal{L}_2$-theory. If there is an $i$-morphism $\eta:f\Rightarrow g$ with respect to $T_2$, then $f$ and $g$ are sententially equivalent with respect to $T_2$. Therefore, if $f$ is an interpretation $f:T_1\to T_2$, then so is $g$.
\end{proposition}
          \noindent
For $i$-morphisms between interpretations, we have a sensible notion of semantic dual, viz., piecewise definable collections of elementary embeddings between the induced models. 
\begin{proposition}[Dual of an $i$-Morphism]\label{Prop:semantic_dual_i-morphism}
    Let $f,g:S\to T$ be piecewise interpretations and $\eta:f\Rightarrow g$ an $i$-morphism. Then there is a unique natural transformation $\eta^*:f^*\Rightarrow g^*$, the \emph{semantic dual} of $\eta$, such that for all $M\vDash T$, sorts $\sigma$ of $S$, pieces $k\in P_\sigma^f$ and $l\in P_\sigma^g$, and tuples $\bar a\in M(D_\sigma^{f,k})$ and $\bar b\in M(D_\sigma^{g,l})$:
    \begin{equation}
        \bar a\times\bar b\in M(\eta_\sigma^{k,l})
        \quad\Longleftrightarrow\quad
        \eta^*_{M,\sigma}\bigl(\pi_\sigma^f(k,\bar a)\bigr)=\pi_\sigma^g(l,\bar b)
        \tag{Dual}\label{eq:dual}
    \end{equation}
\end{proposition}
\begin{proof}
    \emph{Existence, \eqref{eq:dual}, and Uniqueness:} Let $M\vDash T$. Every element of $f^*M(\sigma)$ is of the form $\pi_\sigma^f(k,\bar a)$ with $\bar a\in M(D_\sigma^{f,k})$. By \eqref{Joint-Exist} and \eqref{Dom/Range} there are $l\in P_\sigma^g$ and $\bar b\in M(D_\sigma^{g,l})$ with $\bar a\times\bar b\in M(\eta_\sigma^{k,l})$, and we set $\eta^*_{M,\sigma}(\pi_\sigma^f(k,\bar a)):=\pi_\sigma^g(l,\bar b)$. By condition \eqref{Unique} this is independent of $(l,\bar b)$, and \eqref{well-def} ensures independence of the representative $(k,\bar a)$. The left-to-right direction of \eqref{eq:dual} holds by this definition together with \eqref{Unique}. For the right-to-left direction, suppose $\eta^*_{M,\sigma}(\pi_\sigma^f(k,\bar a))=\pi_\sigma^g(l,\bar b)$ with $\bar b\in M(D_\sigma^{g,l})$, and let $(l',\bar b')$ be such that $\bar a\times\bar b'\in M(\eta_\sigma^{k,l'})$ and $\bar b'\times\bar b\in M((=_\sigma)_g^{l',l})$. Since $\bar a\times\bar a\in M((=_\sigma)_f^{k,k})$ by \eqref{f-identity-intro} for $f$, condition \eqref{well-def} implies $\bar a\times\bar b\in M(\eta_\sigma^{k,l})$. Uniqueness follows immediately from surjectivity of $\pi_\sigma^f$ and \eqref{eq:dual}.
    
    \smallskip
    \noindent
    \emph{Elementarity of the components:} Fix $\varphi:\bar\sigma\in\mathcal{L}_S$ with $|\bar \sigma|\geq 1$. Then:
\begin{align}
  &\forall \bar k\in P_{\bar\sigma}^f, \forall \bar l\in P_{\bar\sigma}^g:  T\vdash \eta^{\bar k,\bar l}_{\bar x,\bar y}\to
    (\varphi(\bar x)_f^{\bar k}\leftrightarrow\varphi(\bar y)_g^{\bar l})\\
  \Leftrightarrow\quad
  &\forall M\vDash T,\; \forall \bar k\in P_{\bar\sigma}^f, \forall \bar l\in P_{\bar\sigma}^g
    \forall\bar a\in M(\bar\sigma_f^{\bar k}),\ \forall \bar b\in M(\bar\sigma_g^{\bar l}):\\
    &\bar a\times\bar b\in M(\eta^{\bar k,\bar l})\Rightarrow
    \bigl(\bar a\in M(\varphi_f^{\bar k})
      \Leftrightarrow\bar b\in M(\varphi_g^{\bar l})\bigr)\\
  \Leftrightarrow\quad
  &\forall M\vDash T,\;
    \forall\bar a \in f^*M(\bar \sigma):\;
    \bar a\in f^*M(\varphi)
      \iff \eta^*_M(\bar a)\in g^*M(\varphi).
\end{align}
    The first equivalence is an instance of soundness and completeness. The second uses \eqref{Dom/Range}, \eqref{Joint-Exist}, \eqref{eq:dual}, and \cref{cor:extensions_elementwise}. Quantifying over all appropriate $\varphi$, the top line ranges over the instances of \eqref{Element-Embed} at positive arities, and the bottom line states that every $\eta^*_M$ is an elementary embedding.

\smallskip
\noindent
    \emph{Naturality:} For naturality, we need to show that if $j:M\to N$ is an elementary embedding, we have $\eta^*_N\circ f^*j=g^*j\circ \eta^*_M$. So let $j:M\to N$ be elementary and $\pi_\sigma^f(k,\bar a)\in f^*M(\sigma)$. Pick $(l,\bar b)$ with $\eta^*_{M,\sigma}(\pi_\sigma^f(k,\bar a))=\pi_\sigma^g(l,\bar b)$. Then $\bar a\times\bar b\in M(\eta_\sigma^{k,l})$ by \eqref{eq:dual}, so $j(\bar a)\times j(\bar b)\in N(\eta_\sigma^{k,l})$ with $j(\bar a)\in N(D_\sigma^{f,k})$ and $j(\bar b)\in N(D_\sigma^{g,l})$,  by elementarity of $j$. Applying \eqref{eq:dual} in $N$ yields $\eta^*_{N,\sigma}(\pi_\sigma^f(k,j(\bar a)))=\pi_\sigma^g(l,j(\bar b))$, which is the required equation, since $f^*j$ and $g^*j$ act on representatives by $j$.
\end{proof}
\noindent
As a translation $f:\mathcal{L}_1\to \mathcal{L}_2$ which respects identity with respect to an $\mathcal{L}_2$-theory $T$ can be seen as an interpretation $L(\Sigma_1)\to T$, the proposition above also applies to them. We have a partial converse, allowing us to check semantically whether a pre-$i$-morphism is an $i$-morphism.
\begin{proposition}[Semantic Recognition of $i$-Morphisms]\label{prop:semantic-recognition}
    Let $f,g:S\to T$ be piecewise interpretations, let $\theta$ be a pre-$i$-morphism from $f$ to $g$ satisfying \eqref{Dom/Range} with respect to $T$, and suppose there is a family of elementary embeddings $(F_M:f^*M\to g^*M)_{M\vDash T}$ satisfying \eqref{eq:dual} with $F$ in place of $\eta^*$. Then $\theta:f\Rightarrow g$ is an $i$-morphism and $\theta^*=F$.
\end{proposition}
\begin{proof}
    We rely on completeness throughout. For \eqref{well-def}, suppose $\bar a\times\bar b\in M(\theta_\sigma^{k,l})$, and that $(k',\bar a')\in P_\sigma ^f\times M(D_\sigma^{f,k'})$ and $(l',\bar b')\in P_\sigma^g\times  M(D_\sigma^{g,l'})$ represent the same respective classes. Then $F_{M,\sigma}(\pi_\sigma^f(k',\bar a'))=F_{M,\sigma}(\pi_\sigma^f(k,\bar a))=\pi_\sigma^g(l,\bar b)=\pi_\sigma^g(l',\bar b')$ by the left-to-right direction of \eqref{eq:dual}, so the right-to-left direction yields $\bar a'\times \bar b'\in M(\theta_\sigma^{k',l'})$. For \eqref{Joint-Exist}, given $\bar a\in M(D_\sigma^{f,k})$, choose a representative $\pi_\sigma^g(l,\bar b)$ of $F_{M,\sigma}(\pi_\sigma^f(k,\bar a))$ and apply \eqref{eq:dual}. For \eqref{Element-Embed}, elementarity of the components of $F$ allows us to run the elementarity part of the proof of \cref{Prop:semantic_dual_i-morphism} backwards, with the sentential instances of \eqref{Element-Embed} following from the fact that $f^*M\equiv g^*M$ and completeness. Finally, $\theta^*=F$ by the uniqueness part of \cref{Prop:semantic_dual_i-morphism}.
\end{proof}
\noindent
Together with the following corollary, this proposition allows us to show that a given pre-$i$-morphism is equivalent to an $i$-morphism and itself an $i$-morphism at once. The corollary itself is an easy consequence of \eqref{eq:dual} and completeness:
\begin{corollary}[Local Faithfulness]\label{cor:local_faithfulness}
    Let $\eta,\zeta:f\Rightarrow g$ be $i$-morphisms between piecewise interpretations $f,g:S\to T$. Then $\eta\equiv_T\zeta$ iff $\eta^*=\zeta^*$.
\end{corollary}
\noindent
The promised notion of equivalence between translations is realized by the notion of an $i$-isomorphism. Recall that an elementary embedding is automatically injective, so an elementary embedding is an isomorphism iff it is surjective, and that every isomorphism is elementary.
  \begin{definition}[$i$-Isomorphism]\label{def:i-iso}
    An $i$-morphism $\eta:f\Rightarrow g$ between piecewise translations with respect to $T$ is an \emph{$i$-isomorphism} iff $T$ proves    $D_{y}^{g,l} \to
        \bigvee_{k \in P_\sigma^f} \exists z_f^k\, \eta_{z,y}^{k,l}$. If $f$ and $g$ respect identity with respect to $T$, this is equivalent to $\eta^*_M$ being an isomorphism at every model $M\vDash T$. 
    If $\eta:f\Rightarrow g$ is an $i$-isomorphism, we denote by $\eta^{-1}:g\Rightarrow f$ the $i$-isomorphism defined by $(\eta^{-1})_{\sigma}^{l,k}(\bar x,\bar y)=\top_{\bar x\times\bar y}\land\eta_\sigma^{k,l}(\bar y,\bar x)$. 
\end{definition} 
Restricting to $i$-morphisms between interpretations, we can apply \cref{prop:semantic-recognition} with $F:=(\eta^*)^{-1}$ to see that the inverse of an $i$-isomorphism defined above is indeed again an $i$-isomorphism. 
A particularly important class of $i$-isomorphisms in the following will be those between translations agreeing up to a family of bijections of their sets of pieces. Semantically, the interpretations related by such an $i$-isomorphism will merely differ in their choice of index sets for the respective disjoint unions.
\begin{proposition}[$i$-Isomorphism induced by a Bijection of Pieces]\label{prop:piece_bijection_morphism}
    Let $f,g:T_1\to T_2$ be piecewise interpretations and let $(b_\sigma:P_\sigma^f\to P_\sigma^g)_{\sigma\in\mathcal{S}_1}$
    be a family of
    bijections between the families of sets of pieces for each sort, satisfying the following conditions for all
    $\sigma\in \mathcal{S}_1$, $\bar\sigma\in \mathcal{S}_1^*$, $k\in P_\sigma^f$, $\bar k\in P_{\bar \sigma}^f$, $R:\bar \sigma\in \mathcal{R}_1$, including the identity symbols $=_\sigma$, and $\bar z:\sigma_f^k\in \mathcal{V}_2^*$ (resp. $\bar u:\bar \sigma_f^{\bar k}\in \mathcal{V}_2^*$):
   \begin{align}
     &\sigma_f^k = \sigma_g^{b_\sigma(k)} 
       \tag{Sort Condition}\label{Sort_Condition} \\
     T_2\vdash & D_\sigma^{f,k}(\bar z)
       \leftrightarrow D_\sigma^{g,b_\sigma(k)}(\bar z)
       \tag{Domain Condition}\label{Domain_Condition}\\
     T_2\vdash &D_{\bar\sigma}^{f,\bar k}(\bar u)\to\bigl(R_f^{\bar k}(\bar u)
  \leftrightarrow R_g^{b(\bar k)}(\bar u)\bigr)\tag{Relation Condition} \label{Relation_Condition}    
\end{align} where $b(\bar k):=(b_{\sigma_1}(k_1),\ldots,b_{\sigma_n}(k_n))$. Then the family of formulae
    $(\eta_b)_\sigma : P_\sigma^f \times P_\sigma^g \to \mathcal{L}_2$
    defined by
    \begin{equation}
       (\eta_b)_\sigma^{k,l}(\bar z,\bar w) \;:=\;
        D_\sigma^{f,k}(\bar z) \land D_\sigma^{g,l}(\bar w)
        \land (=_\sigma)_g^{b_\sigma(k),\,l}(\bar z,\bar w),
    \end{equation} where $\bar z:\sigma_f^{k}$ and $\bar w:\sigma_g^{ l}$, is an $i$-isomorphism $\eta_b:f\Rightarrow g$ with respect to $T_2$. Furthermore, the dual $\eta_b^*$ acts via $(\eta_b^*)_M(\pi_\sigma^f(k,\bar a))=\pi_\sigma^g(b_\sigma(k),\bar a)$.
\end{proposition}
\begin{proof}
  Apply \cref{prop:semantic-recognition} with $F_{M,\sigma}(\pi_\sigma^f(k,\bar a)):=\pi_\sigma^g(b_\sigma(k),\bar a)$; well-definedness of $F$ uses the \eqref{Sort_Condition}, the \eqref{Domain_Condition}, and the \eqref{Relation_Condition} for $=_\sigma$. From the latter, we also get injectivity, while surjectivity follows from the \eqref{Domain_Condition} and bijectivity of $b_\sigma$. Finally, that $F_M$ is an isomorphism and hence elementary follows from the \eqref{Relation_Condition} for all $R\in \mathcal{R}_1$.
\end{proof}
\noindent
This proposition also supplies us with an easy way to see that if $f\equiv_Tg$, then $f$ and $g$ are $i$-isomorphic, by taking $b=\id$ in the above proposition. Moreover, given two pairs of translations $(f,g)$ and $(f',g')$ such that $f\equiv_T f'$ and $g\equiv_T g'$, and a family $(b_\sigma:P_\sigma^f\to P_\sigma^g)$ satisfying the conditions above for $(f,g)$, it also satisfies them for $(f',g')$, and the resulting $i$-isomorphisms $\eta_b$ and $\eta_b'$ are provably equivalent.
\section{The Bicategory of Piecewise Interpretations}\label{sec:bicategory}
From a category-theoretic point of view, we have yet to establish that piecewise interpretations make for a good notion of morphism between theories, because we have not explained how they compose. In this section, we will show that piecewise interpretations and $i$-morphisms between them form a bicategory by developing appropriate notions of composition for them. That they form a bicategory means, in particular, that composition of piecewise interpretations is associative only up to a coherent $i$-isomorphism.\footnote{\citet{d2020bicategories}, working in coherent logic, also obtain a non-strict bicategory of interpretations, but that is because they allow terms beyond variables in their languages, leading to non-trivial unitors. We expect that this approach to terms in interpretations could be combined neatly with the theory here.} The overall strategy will be to avoid tedious syntactic verifications of the bicategory properties by employing the developed semantic duals. Let us start by recalling the definition of the notion of bicategory we will be using:
\begin{definition}[Bicategory]\label{def:bicategory}
A \emph{bicategory} $\mathcal{B}$ consists of: \emph{a class of objects} $\mathrm{ob}(\mathcal{B})$; for all $X,Y\in\mathrm{ob}(\mathcal{B})$ a \emph{hom-category} $\bigl(\mathcal{B}(X,Y),\bullet,\id\bigr)$, whose objects are called \emph{$1$-morphisms} $f: X\to Y$ and whose morphisms are called \emph{$2$-morphisms} $\eta: f\Rightarrow f': X\to Y$ of $\mathcal{B}$; composition in $\mathcal{B}(X,Y)$ is called \emph{vertical composition}; for each object $X$ an \emph{identity $1$-morphism} $\id_X:X\to X$, and for all $X,Y,Z$ a composition functor
        \begin{equation}
             c_{X,Y,Z}:\mathcal{B}(Y,Z)\times\mathcal{B}(X,Y)\longrightarrow\mathcal{B}(X,Z);
        \end{equation}
 and, finally, invertible $2$-morphisms, namely the associators and left and right unitors,
 \begin{equation}
          a_{h,g,f}: (h\circ g)\circ f \Rightarrow h\circ (g\circ f): W\to Z
          \end{equation}
          \begin{equation}
              l^{f}: \id_X\circ f \Rightarrow f: W\to X \qquad
          r^{f}: f\circ \id_W \Rightarrow f: W\to X,
          \end{equation}
        for all $f: W\to X$, $g: X\to Y$, $h: Y\to Z$, natural in
        each argument and satisfying the \emph{pentagon} and \emph{triangle} identities (see \cref{rem:prebicategory_to_bicategory}). A bicategory in which $a$, $l$ and $r$ are all identities is called a \emph{$2$-category}.
\end{definition}
On a pair of $1$-morphisms $(g,f)$, the action of the composition functor is denoted by $g\circ f$. On a pair of $2$-morphisms $(\eta,\zeta)$ it is denoted by $\eta\star \zeta$ and called \emph{horizontal composition}. Two special cases of horizontal composition are worth singling out. For a $1$-morphism $f: X\to Y$ and a $2$-morphism $\zeta: h\Rightarrow h': W\to X$, the horizontal composite $\id_f \star \zeta: f\circ h\Rightarrow f\circ h'$
is denoted by $f\zeta$ and called the \emph{left whiskering of $\zeta$ by $f$}; dually, for $\eta: g\Rightarrow g': Y\to Z$ the composite $\eta \star \id_f: g\circ f\Rightarrow g'\circ f$ is denoted by $\eta f$ and called the \emph{right whiskering of $\eta$ by $f$}.

\medskip\noindent
In the case of piecewise interpretations, the objects of our bicategory $\Thpw$ are given by classical first-order theories. Given two theories $S$ and $T$, the collection of $1$-morphisms $\Thpw(S,T)$ between them is given by the set of (equivalence classes under provable equivalence of) piecewise interpretations $f:S\to T$. For any theory $T$, the identity interpretation $\id_T:T\to T$ is defined as follows: For each sort $\sigma$ of $T$, the set of pieces is $P_\sigma=\{0\}$, with $\sigma_{\id_T}^{0}:=\sigma$, and $x_{\id_T}^0=x$ for each variable $x:\sigma$. The domain formula $D_\sigma^{\id_T,0}$ is given by the verum $\top_\sigma$. For an atomic formula $\varphi$, we have $\varphi_{\id_T}^{\bar 0}:=\varphi$. It follows by induction that for any $\mathcal{L}$-formula $\varphi$, $\bar{x}(\varphi_{\id_T}^{\bar 0})=\bar{x}(\varphi)$ and $\vdash \varphi_{\id_T}^{\bar 0}\leftrightarrow \varphi$. The semantic dual of $\id_T$ is naturally isomorphic to $\id_{\Mod(T)}$ via the natural isomorphism 
\begin{equation}\label{ref:mod_unitor}
    \alpha_{T}: \id_T^*\overset{\sim}{\Rightarrow}\id_{\Mod(T)} \qquad \pi_\sigma^{\id_T}(0,a)\mapsto a.
\end{equation}

\medskip\noindent
The theory of semantic duals we have developed will allow us to define a locally faithful contravariant pseudofunctor from the bicategory $\Thpw$ into the $2$-category of categories $\mathbf{CAT}$, whose objects are categories, $1$-morphisms are functors, and $2$-morphisms are natural transformations between functors, with composition the usual compositions of functors and natural transformations. We will develop the bicategory $\Thpw$ and this pseudofunctor concurrently, using it to show that we are indeed defining a bicategory. Here is the general definition of a pseudofunctor:

\begin{definition}[Pseudofunctor]\label{def:pseudofunctor}
A \emph{pseudofunctor} $F:\mathcal{B}\to\mathcal{C}$ between bicategories
consists of a function $F:\mathrm{ob}(\mathcal{B})\to\mathrm{ob}(\mathcal{C})$, a functor $F_{X,Y}:\mathcal{B}(X,Y)\to\mathcal{C}(FX,FY)$ for all $X,Y\in\mathcal{B}$, an invertible $2$-morphism $\phi_X: \id_{FX}\Rightarrow F(\id_X)$ for each $X\in\mathcal{B}$ (the \emph{unitor}), an invertible $2$-morphism $\phi_{g,f}: Fg\circ Ff\Rightarrow F(g\circ f)$ for all $f: X\to Y$ and $g: Y\to Z$, natural in $f$ and $g$ (the \emph{compositor}), subject to two unit identities and an associativity identity:
\begin{align}
  F(l^{f})\bullet\phi_{\id_X,f}\bullet\bigl(\phi_{X}\,Ff\bigr)
    &= l^{Ff},
  \tag{PF-Left}\label{ax:pf-left}\\
  F(r^{f})\bullet\phi_{f,\id_W}\bullet\bigl(Ff\,\phi_{W}\bigr)
    &= r^{Ff},
  \tag{PF-Right}\label{ax:pf-right}\\
  \phi_{h,\,g\circ f}\bullet\bigl(Fh\,\phi_{g,f}\bigr)\bullet a_{Fh,Fg,Ff}
    &= F(a_{h,g,f})\bullet\phi_{h\circ g,\,f}\bullet\bigl(\phi_{h,g}\,Ff\bigr),
  \tag{PF-Assoc}\label{ax:pf-assoc}
\end{align}
for all $f: W\to X$, $g: X\to Y$ and $h: Y\to Z$. 
A contravariant pseudofunctor $G:\mathcal{B}\to \mathcal{C}$ is a pseudofunctor $G:\mathcal{B}^{op}\to \mathcal{C}$, where $\mathcal{B}^{op}$ is the bicategory defined by setting $\mathcal{B}^{op}(X,Y):=\mathcal{B}(Y,X)$, using the inverse associator as the new associator, and exchanging left and right unitor. Its compositor and unitor are invertible $2$-morphisms $Gf\circ Gg\Rightarrow G(g\circ f)$ and $\id_{GX}\Rightarrow G(\id_X)$.
\end{definition} 

For the semantic duals, we will present these structural morphisms by their inverses. In this presentation, the axioms are spelled out more explicitly in Appendix~\ref{app:bicategorical_results}, \cref{def:pseudofunctor_from_prebicategory}. For our purposes, the following basic property of pseudofunctors will be important:
\begin{definition}[Local Faithfulness]\label{def:locally-faithful}
A (contravariant) pseudofunctor $F:\mathcal{B}\to\mathcal{C}$ is
\emph{locally faithful} if it is injective on $2$-morphisms, that is, iff $F(\eta)=F(\zeta)$ implies $\eta=\zeta$ for all parallel $2$-morphisms $\eta,\zeta:f\Rightarrow g$.
\end{definition}
\subsection{Composition of Piecewise Interpretations}
As a next step, we will define the composition of piecewise interpretations. That is, we will define the object part $\circ$ of the composition functor one needs to obtain a bicategory.
\begin{definition}[Composition of Translations]\label{def:composition}
    Let $f:\mathcal{L}_1\to \mathcal{L}_2$ and $g:\mathcal{L}_2\to \mathcal{L}_3$ be piecewise translations. The reconstrual underlying their composition $g\circ f$ is given by the following data: The \emph{sort function} $(g\circ f)_{\mathcal{S}}$ assigns to each sort the family $\sigma_{g\circ f}^{(k,\bar l)}:=(\sigma_f^k)_g^{\bar l}$ with index set $P_\sigma^{g\circ f}:=\{(k,\bar l)\mid k\in P_\sigma^f,\;\bar l\in P_{\sigma_f^k}^g\}.$ The variable and domain functions below are indexed by this index set. The \emph{variable function} $(g\circ f)_{\mathcal{V}}$ sends a variable $x:\sigma$ to the family $x_{g\circ f}^{(k,\bar l)}:=(x_f^k)_g^{\bar l}$. The \emph{domain function} $(g\circ f)_{\mathcal{D}}$ assigns to each sort $\sigma\in\mathcal{S}_1$ the family $D_\sigma^{g\circ f, (k,\bar l)}:=(D_\sigma^{f,k})_g^{\bar l}\land D^{g,\bar l}_{\bar x(D_\sigma^{f,k})}$. The composite \emph{relation function} $(g\circ f)_{\mathcal{R}}$ sends a relation symbol $R\in\mathcal R_1$ of arity $\bar\sigma=(\sigma_1,\ldots,\sigma_n)$ to the family of formulae $R_{g\circ f}^{((k_1,\bar l_1),\ldots,(k_n, \bar l_n))}:=(R_f^{\bar k})_g^{\bar l}$,
        where $\bar k:=(k_1,\ldots,k_n)$ and $\bar l:= \bar l_1\times \cdots\times\bar l_n$.
\end{definition}

One can easily check that these data satisfy the no-overlap condition, duplicate-freeness and the sortal conditions of \cref{def:reconstrual}. In the following, we will need that the induced composed translation schematically handles general formulae just as it does relation symbols. To that end, one first needs to establish in a technical lemma that substitution and translation commute in an appropriate sense, which is done in \cref{lem:subst-translation} in Appendix~\ref{app:cc_and_subst}. 
\begin{restatable}[Composition for Complex Formulae]{proposition}{compositioncomplex}
\label{composition_complex}
Let $f:\mathcal L_1\to\mathcal L_2$ and $g:\mathcal L_2\to\mathcal L_3$ be
piecewise translations, and let $\varphi:\bar\sigma$ be an $\mathcal L_1$-formula
with $|\bar\sigma|=n$. For any composite piece $\bar m=((k_1,\bar l_1),\ldots,(k_n,\bar l_n))\in P_{\bar\sigma}^{g\circ f}$, we set
\begin{equation}
    \bar m_f:=(k_1,\ldots,k_n)\in P_{\bar\sigma}^f,
\qquad
\bar m_g:=\bar l_1\times\cdots\times\bar l_n\in P_{\bar\sigma_f^{\bar m_f}}^g.
\end{equation}
Then the formulae $\varphi_{g\circ f}^{\bar m}$ and $(\varphi_f^{\bar m_f})_g^{\bar m_g}$ have the same canonical context, and
\begin{equation}
\vdash_{L(\Sigma_3)} \varphi_{g\circ f}^{\bar m}
\leftrightarrow
(\varphi_f^{\bar m_f})_g^{\bar m_g}.
\end{equation}
\end{restatable}
\noindent
 We write $\varphi \mathrel{\dashv\!\;\vdash_T}\psi$ to denote that $\varphi$ and $\psi$ have the same canonical context and are provably equivalent with respecto to the theory $T$.
\begin{proof}
   See Appendix~\ref{app:cc_and_subst}.
\end{proof}
\begin{remark}[Composite Pieces]\label{rem:composite_pieces}
In the situation above, since $\bar\sigma_f^{\bar m_f}=\sigma_{1,f}^{k_1}\times\cdots\times\sigma_{n,f}^{k_n}$, the assignment
$\bar m\mapsto(\bar m_f,\bar m_g)$ is a bijection from $P^{g\circ f}_{\bar\sigma}$ to the set of pairs
$(\bar k,\bar l)$ with $\bar k\in P^f_{\bar\sigma}$ and $\bar l\in P^g_{\bar\sigma_f^{\bar k}}$.
\end{remark}
\noindent
For this composition operation to fit into the bicategory we are developing, the composition of two piecewise interpretations also needs to be a piecewise interpretation, which follows from the preceding proposition specialized to sentences:
\begin{corollary}[Composition of Interpretations is an Interpretation]\label{cor:composite_interpretations}
    Let $f:T_1\to T_2$ and $g:T_2\to T_3$ be piecewise interpretations. Then the composition $g\circ f$ is an interpretation $g\circ f:T_1\to T_3$.
\end{corollary}
\noindent
Moreover, \cref{composition_complex} allows us to explicitly describe the semantic dual of the composition of two piecewise interpretations.
\begin{lemma}[Extensions under Composition]\label{lem:composite-extensions}
Let $f:T_1\to T_2$ and $g:T_2\to T_3$ be piecewise interpretations. Let
$M\vDash T_3$, let $\bar\sigma\in\mathcal S_1^*$, and let $\bar m\in P^{g\circ f}_{\bar\sigma}$
with $\bar m_f$ and $\bar m_g$ as in \cref{composition_complex}. Then
\begin{equation}\label{eq:composite-domain}
M\bigl(D^{g\circ f,\bar m}_{\bar\sigma}\bigr)
=\Bigl\{\,\bar a\in M\bigl(D^{g,\bar m_g}_{\bar\sigma_f^{\bar m_f}}\bigr)
\;\Bigm|\;\pi^g_{\bar\sigma_f^{\bar m_f}}(\bar m_g,\bar a)
\in g^*M\bigl(D^{f,\bar m_f}_{\bar\sigma}\bigr)\Bigr\},
\end{equation}
and for every $\mathcal L_1$-formula $\varphi:\bar\sigma$ and every
$\bar a\in M(D^{g\circ f,\bar m}_{\bar\sigma})$,
\begin{equation}\label{eq:composite_formulae}
    \pi^{g\circ f}_{\bar\sigma}(\bar m,\bar a)\in (g\circ f)^*M(\varphi)
\iff \pi^f_{\bar\sigma}\bigl(\bar m_f,\pi^g_{\bar\sigma_f^{\bar m_f}}(\bar m_g,\bar a)\bigr)
      \in f^*(g^*M)(\varphi).
\end{equation}
\end{lemma}
\begin{proof}
    The domain equation \eqref{eq:composite-domain} follows from the definition of the composite domain formula and \eqref{eq:extension_elementwise} for $g$. For the equivalence for formulae \eqref{eq:composite_formulae}, let $\varphi:\bar \sigma$ be an $\mathcal{L}_1$-formula, $\bar m\in P_{\bar \sigma}^{g\circ f}$, and $\bar a\in M(D_{\bar \sigma}^{g\circ f,\bar m})$. We get the following chain of equivalences: 
    \begin{align}
    \pi^{g\circ f}_{\bar \sigma}(\bar m,\bar a)\in (g\circ f)^*M(\varphi) \overset{\text{\eqref{eq:extension_elementwise}}}{\Leftrightarrow}&\bar a\in M(\varphi_{g\circ f}^{\bar m})\overset{\cref{composition_complex}}{\Leftrightarrow} \bar a\in M((\varphi_f^{\bar m_f})_g^{\bar m_g})\\
    \overset{\eqref{eq:extension_elementwise}}{\Leftrightarrow}& \pi_{\bar{\sigma}_f^{\bar m_f}}^g(\bar m_g,\bar a)\in g^*M(\varphi_f^{\bar m_f})\\ \overset{\eqref{eq:extension_elementwise}}{\Leftrightarrow}&\pi_{\bar \sigma}^f(\bar m_f,\pi_{\bar{\sigma}_f^{\bar{m}_f}}^g(\bar m_g,\bar a))\in f^*(g^*M)(\varphi)
\end{align}
\end{proof}
\noindent An immediate corollary, using \cref{prop:prov_equiv_char}, is that composition descends to provable equivalence classes of interpretations. Moreover, this explicit description allows us to construct the compositor we will need to think of the semantic dual as being part of a pseudofunctor:
\begin{proposition}[Compositor for Semantic Duals]
\label{prop:alpha_gf}
Let $f:T_1\to T_2$ and $g: T_2\to T_3$ be piecewise
interpretations. Then there is a natural isomorphism
\begin{equation}
     \alpha_{g,f}: (g\circ f)^* \Rightarrow f^*g^*: \Mod(T_3)\to\Mod(T_1).
\end{equation}
\noindent
For a $T_3$-model $M$, the component
$\alpha_{g,f,M}: (g\circ f)^*M\to f^*(g^*M)$ is the
isomorphism whose $\sigma$-component is defined by
\begin{equation}
    \pi^{g\circ f}_\sigma((k,\bar l),\bar a)
  \longmapsto
  \pi^f_\sigma\bigl(k,\pi^g_{\sigma_f^{k}}(\bar l,\bar a)\bigr).
\end{equation}
\noindent
Here $k\in P^{ f}_{\sigma}, \bar l\in P_{\sigma_f^k}^{g}$ and
$\bar a\in M(D^{g\circ f, (k,\bar l)}_\sigma)$.
\end{proposition}

\begin{proof}
Recall that every $x\in (g\circ f)^*M(\sigma)$ is of the form $ x=\pi^{g\circ f}_{\sigma}((k,\bar l),\bar a)$ with $(k,\bar l)\in P_{\sigma}^{g\circ f}$ and $\bar a\in M(D_\sigma^{g\circ f,(k,\bar l)})$, whence $\alpha_{g,f}$ is total. Well-definedness and injectivity of $(\alpha_{g,f,M})_\sigma$ follow from \eqref{eq:composite_formulae} applied to $=_\sigma$. Surjectivity follows from construction and \eqref{eq:composite-domain}. That $\alpha_{g,f,M}$ is an elementary embedding follows immediately from the collection of equivalences given by \eqref{eq:composite_formulae} for all formulae $\varphi\in\mathcal{L}_1$. 

\medskip
\noindent
To demonstrate naturality, consider an elementary embedding $j: M\to N$ and compute:
\begin{align}
  &f^*(g^*j)\bigl(\pi^f_\sigma(k,\pi^g_{\sigma_f^{ k}}(\bar l,\bar a))\bigr)
  = \pi^f_\sigma\bigl(k,\pi^g_{\sigma_f^{ k}}( \bar l,j(\bar a))\bigr)\\
  &= \alpha_{g,f,N,\sigma}\bigl(\pi^{g\circ f}_\sigma((k, \bar l),j(\bar a))\bigr)
   = \alpha_{g,f,N,\sigma}\circ(g\circ f)^*j\bigl(\pi^{g\circ f}_\sigma(( k,\bar l),\bar a)\bigr).
\end{align}
\end{proof}
\noindent
Moving on to the identity of our bicategory, the identity interpretation indeed behaves as an identity with respect to composition,
up to canonical $i$-isomorphisms, the unitors of our bicategory:
\begin{definition}[Left and Right Unitors]\label{def:unitor}
    Let $f:T\to U$ be an interpretation. The left unitor $l^f:\id_U\circ f\Rightarrow f$ and right unitor $r^f:f\circ \id_T\Rightarrow f$ of $f$ are $i$-isomorphisms defined via bijections of pieces
according to \cref{prop:piece_bijection_morphism}:
\begin{align}
  l^f_\sigma: P_\sigma^{\id_U\circ f} &\to P_\sigma^f
  &
  r^f_\sigma: P_\sigma^{f\circ \id_T} &\to P_\sigma^f
  \\
  (k,\,0^{|\sigma_f^k|}) &\mapsto k
  &
  (0,\,k) &\mapsto k
\end{align}
Here $0^{|\sigma_f^k|}$ denotes the unique element of
$\prod_{i=1}^{|\sigma_f^k|}P_{(\sigma_f^k)_i}^{\id_U}$.
\end{definition}
That these bijections satisfy the conditions
of \cref{prop:piece_bijection_morphism} follows easily from the
definition and provable triviality of $\id_U$ and $\id_T$. Under the semantic dual, the unitors are appropriately coherent with the compositors and unitors of the semantic dual:
\begin{lemma}[Coherence of the Unitors]\label{lem:unitor-coherence}
Let $f:T\to U$ be a piecewise interpretation and recall that for a theory $S$, $\alpha_{S}:(\id_S)^*\overset{\sim}{\Rightarrow}\id_{\Mod(S)}$ is the natural isomorphism defined by $\pi_\sigma^{\id_S}(0,a)\mapsto a$. Then
\begin{equation}\label{eq:unitor_coherence}
    (l^f)^{*}=\bigl(f^{*}\alpha_{U}\bigr)\bullet\alpha_{\id_U,f}
    \qquad\text{and}\qquad
    (r^f)^{*}=\bigl(\alpha_{T}f^{*}\bigr)\bullet\alpha_{f,\id_T}.
\end{equation}
\end{lemma}

\begin{proof}
Let $M\vDash U$ and $\sigma\in\mathcal{S}_T$. We first consider the left unitor: Given an element $x$ of
$(\id_U\circ f)^{*}M(\sigma)$, $x$ is of the form
$\pi_\sigma^{\id_U\circ f}((k,0^{|\sigma_f^k|}),\bar a)$ with $(k, \bar 0, \bar a)\in P^f_\sigma\times \{0\}^{|\sigma_f^k|}\times M(D_\sigma^{f,k})$. By \cref{prop:piece_bijection_morphism}, $(l^f)^{*}_M\bigl(\pi_\sigma^{\id_U\circ f}((k,0^{|\sigma_f^k|}),\bar a)\bigr)
 =\pi_\sigma^{f}(k,\bar a)$. On the right-hand side of the left equation, by
\cref{prop:alpha_gf} and the definitions of whiskering and
$\alpha_{U}$,
\begin{equation}
    \bigl(f^{*}\alpha_{U}\bigr)_M\circ\alpha_{\id_U,f,M}\bigl(
    \pi_\sigma^{\id_U\circ f}\bigl((k,0^{|\sigma_f^k|}),\bar a\bigr)\bigr)
    = (f^{*}\alpha_{U})_M\big(
    \pi_\sigma^{f}\bigl(k,\,\pi^{\id_U}_{\sigma_f^k}(0^{|\sigma_f^k|},\bar a)\bigr)\bigr)
    =
    \pi_\sigma^{f}(k,\bar a).
\end{equation}
Moving on to the right unitor, given an element $x$ of $(f\circ\id_T)^{*}M(\sigma)$, $x$ is of the form
$\pi_\sigma^{f\circ\id_T}((0,k),\bar a)$ with $(k,\bar a)\in P^f_\sigma\times M(D_\sigma^{f,k})$, and $(r^f)^{*}_M(x)=\pi_\sigma^f(k,\bar a)$ by \cref{prop:piece_bijection_morphism}. On the right-hand side of the right equation, we see that
\begin{equation}
    \bigl(\alpha_{T}f^{*}\bigr)_M\bigl(\alpha_{f,\id_T,M}(
    \pi_\sigma^{f\circ\id_T}\bigl((0,k),\bar a\bigr))\bigr)
=\bigl(\alpha_{T}f^{*}\bigr)_M\bigl(\pi_\sigma^{\id_T}\bigl(0,\,\pi^{f}_{\sigma}(k,\bar a)\bigr)\bigr)
    =
    \pi_\sigma^{f}(k,\bar a).
\end{equation}
\end{proof}

\subsection{Associator for Piecewise Interpretations}
The next ingredient for our bicategory is an associator for each appropriate triple of interpretations, making composition associative up to a coherent $i$-isomorphism. We thereby answer a small question due to Visser, about the sense in which piecewise interpretations compose associatively (personal communication).
\begin{definition}[Associator for Interpretations]\label{def:associator}
    Let $f:T_1\to T_2, g:T_2\to T_3$ and $h:T_3\to T_4$ be piecewise interpretations. Given a sort $\sigma\in\mathcal S_1$, set $\tau_i=(\sigma_f^k)_i$ for each $k\in P_\sigma^f$ and $1\leq i\leq|\sigma_f^k|$. Then, iterating the definition of composition yields:
    \begin{align}
        &P^{(h\circ g)\circ f}_\sigma=\left\{\bigl(k,\,((l_1,\bar m_1),\ldots,(l_{| \sigma_f^k|},\bar m_{| \sigma_f^k|}))\bigr)\,\middle| \, \begin{aligned}
            &k\in P^f_\sigma, l_i\in P^g_{\tau_i}, \bar m_i\in P^h_{(\tau_i)^{l_i}_g}\\
            & \text{with }1\leq i\leq | \sigma_f^k|
        \end{aligned}
        \right\}\\
        \intertext{For the other bracketing, one computes:}
        &P^{h\circ (g\circ f)}_\sigma=\left\{\bigl((k,\bar l),\,\bar m_1\times\cdots\times\bar m_{| \sigma_f^k|}\bigr)\,\middle| \
       \begin{aligned}
            &k\in P^f_\sigma, l_i\in P^g_{\tau_i}, \bar m_i\in P^h_{(\tau_i)^{l_i}_g}\\
            & \text{with }1\leq i\leq | \sigma_f^k|
        \end{aligned}
        \right\}
    \end{align}
    Inspecting these two sets, we see that they can be matched up bijectively by regrouping and concatenating pieces. The associator $a_{h,g,f}$ is the $i$-isomorphism induced by the following family of bijections, as in \cref{prop:piece_bijection_morphism}:
    \begin{align}
        (a_{h,g,f})_\sigma:P^{(h\circ g)\circ f}_\sigma&\to P^{h\circ (g\circ f)}_\sigma\\
        \bigl(k,\,((l_1,\bar m_1),\ldots,(l_{| \sigma_f^k|},\bar m_{| \sigma_f^k|}))\bigr)&\mapsto \bigl((k, \bar l),\,\bar m_1\times\cdots\times\bar m_{| \sigma_f^k|}\bigr),
    \end{align}
    where $\bar l=(l_1,\ldots,l_{|\sigma^k_f|})$.
\end{definition}

\begin{proposition}[Associativity of Composition]\label{prop:comp_assoc}
    Let $f:T_1\to T_2, g:T_2\to T_3$ and $h:T_3\to T_4$ be piecewise interpretations. Then the associator $a_{h,g,f}:(h\circ g)\circ f \Rightarrow h\circ (g\circ f)$ is an $i$-isomorphism.
\end{proposition}
\begin{proof}
    We verify the instances of the conditions of \cref{prop:piece_bijection_morphism}
    for a bijection $(a_{h,g,f})_\sigma$ for an arbitrary sort $\sigma$. Write $p=(k,((l_1,\bar m_1),\ldots,(l_n,\bar m_n)))\in P_\sigma^{(h\circ g)\circ f}$
    with $n=|\sigma_f^k|$, and abbreviate $\tau_i:=(\sigma_f^k)_i$,
    $\bar l=(l_1,\ldots,l_n)$, and $\bar m=\bar m_1\times\cdots\times\bar m_n$.

   \medskip\noindent\emph{\eqref{Sort_Condition}:}
    Iterating the definition of the sort function of a composition on each side yields:
    \begin{align}
        \sigma_{(h\circ g)\circ f}^{p}
        &= (\tau_1)_{h\circ g}^{(l_1,\bar m_1)}\times\cdots\times
          (\tau_n)_{h\circ g}^{(l_n,\bar m_n)}
        = ((\tau_1)_g^{l_1})_h^{\bar m_1}\times\cdots\times
          ((\tau_n)_g^{l_n})_h^{\bar m_n}\\
          &=  ((\tau_1)_g^{l_1}\times\cdots\times
          (\tau_n)_g^{l_n})_h^{\bar m} = (\sigma_{g\circ f}^{(k,\bar l)})_h^{\bar m}
          \\
          &=
        \sigma_{h\circ(g\circ f)}^{a(p)}
    \end{align}
 \medskip\noindent\emph{\eqref{Domain_Condition}:} For the domain condition, we compute:
\begin{align}
     D_x^{(h\circ g)\circ f, p}\;&\eqv_{L(\Sigma_4)}\;
           ((D_x^{f,k})_g^{\bar l})_h^{\bar m}
           \;\land\; \bigwedge_{i=1}^n
           \Bigl[(D_{(x_f^{k})_i}^{g,l_i})_h^{\bar m_i}
           \;\land\; D_{((x_f^{k})_i)_g^{l_i}}^{h,\bar m_i}\Bigr]\\
           \;&\eqv_{L(\Sigma_4)}\; \bigl((D_x^{f,k})_g^{\bar l}
           \;\land\; \bigwedge_{i=1}^n D_{(x_f^{k})_i}^{g,l_i}\bigr)_h^{\bar m}
           \;\land\; \bigwedge_{i=1}^n D_{((x_f^{k})_i)_g^{l_i}}^{h,\bar m_i}\\
           &=\;D_x^{h\circ(g\circ f),a(p)},
\end{align} where the first $\eqv$ writes out the definition of the composite domain formula and applies \cref{composition_complex} at the intermediate stage for the left-most conjunct, the second distributes and applies the binary and the context clauses. 

    \medskip\noindent\emph{\eqref{Relation_Condition}:}
   Let $R:\bar\sigma$ be a relation symbol of $\mathcal{L}_1$ and let
   $\bar q\in P_{\bar\sigma}^{(h\circ g)\circ f}$. Applying the piece decomposition from
   \cref{composition_complex} first to $f$ and $h\circ g$, then to $g$ and $h$, set
   $\bar k:=\bar q_f$, $\bar l:=(\bar q_{h\circ g})_g$ and $\bar m:=(\bar q_{h\circ g})_h$.
   By the definition of the extension of $a$ on sequences, $a(\bar q)$ decomposes into the same $\bar k,\bar l,\bar m$. Iterative applications of \cref{composition_complex} yield
   \begin{equation}
       R_{(h\circ g)\circ f}^{\bar q}
       \;\eqv_{L(\Sigma_4)}\; \bigl((R_f^{\bar k})_g^{\bar l}\bigr)_h^{\bar m}
       \;\eqv_{L(\Sigma_4)}\; R_{h\circ(g\circ f)}^{a(\bar q)},
   \end{equation}
   which clearly entails the relation condition.

\end{proof}

 \medskip
 \noindent
 We can describe the semantic dual of the associator using the compositor from \cref{prop:alpha_gf}. Recall that if $F:C\to D$ and $G_1,G_2:D\to E$ are functors and $\eta:G_1\Rightarrow G_2$ is a natural transformation, then $\eta F:G_1F\Rightarrow G_2F$, the right whiskering of $\eta$ by $F$, denotes the natural transformation defined via $(\eta F)_X:=\eta_{F(X)}$.
 \begin{proposition}[Coherence of the Associator]\label{prop:associator_coherence}
     Let $f:T_1\to T_2, g:T_2\to T_3$ and $h:T_3\to T_4$ be piecewise interpretations. Then the following holds for the associator $a_{h,g,f}$:
     \begin{equation}\label{eq:associativity_coherence}
         a_{h,g,f}^*=\alpha_{h,g\circ f}^{-1} \bullet\alpha_{g,f}^{-1}h^*\bullet f^*\alpha_{h,g}\bullet\alpha_{h\circ g,f}
     \end{equation}
     
 \end{proposition}
 \begin{proof}
 
 Let $M\vDash T_4$, $\sigma$ be a sort of $\mathcal{L}_1$, and $x\in ((h\circ g)\circ f)^*M(\sigma)$. Then there is a piece $ p=(k,((l_1,\bar m_1),\ldots,(l_n,\bar m_n)))\in P_\sigma^{(h\circ g)\circ f}$ and a tuple $\bar a\in M(D_\sigma^{(h\circ g)\circ f,p})$ such that $x=\pi^{(h\circ g)\circ f}_\sigma( p,\bar a)$. Set $\tau_i:=(\sigma_f^k)_i$ and $n:=|\sigma_f^k|$. Then $\bar a$ can be written as $\bar a=\bar b_1\times\cdots\times \bar b_n$ with $\bar b_i\in M(D_{\tau_i}^{h\circ g,(l_i,\bar m_i)})$. Set $\bar l=(l_1,\ldots,l_n)$, $\bar m=\bar m_1\times\cdots\times\bar m_n$, and $\bar q=((l_1,\bar m_1),\ldots,(l_n,\bar m_n))$. For each $i$, define 
 \begin{equation}
     e_i:=\pi_{\tau_i}^{h\circ g}((l_i,\bar m_i),\bar b_i),  \quad \text{and} \quad  \bar c_i:=\pi_{(\tau_i)_g^{l_i}}^h(\bar m_i,\bar b_i).
 \end{equation}
 Set $ \bar c=\bar c_1\times\cdots\times \bar c_n$. We compute the right-hand side evaluated at $x$:\footnote{We omit the subscript $M$ after the first equality to reduce visual clutter, keeping in mind that right whiskering operates on it.} \begin{align} &\bigl( \alpha_{h,g\circ f}^{-1}\bullet(\alpha_{g,f}^{-1}h^*)\bullet (f^*\alpha_{h,g})\bullet\alpha_{h\circ g,f} \bigr)_M \bigl(\pi_\sigma^{(h\circ g)\circ f}(p,\bar a)\bigr)\\ 
 &= \alpha_{h,g\circ f}^{-1}(\alpha_{g,f}^{-1}h^*)(f^*\alpha_{h,g}) \bigl( \pi_\sigma^f(k,\pi_{\sigma_f^k}^{h\circ g}(\bar q,\bar a)) \bigr)\label{eq:assoc_two}\\
 &= \alpha_{h,g\circ f}^{-1}(\alpha_{g,f}^{-1}h^*)(f^*\alpha_{h,g}) \bigl( \pi_\sigma^f(k,(e_i)_{i=1}^n) \bigr)\\
 &= \alpha_{h,g\circ f}^{-1}(\alpha_{g,f}^{-1}h^*) \bigl( \pi_\sigma^f(k,((\alpha_{h,g})_M(e_i))_{i=1}^n) \bigr)= \alpha_{h,g\circ f}^{-1}(\alpha_{g,f}^{-1}h^*) \bigl( \pi_\sigma^f(k,(\pi_{\tau_i}^g(l_i,\bar c_i))_{i=1}^n) \bigr)\\
 &= \alpha_{h,g\circ f}^{-1}(\alpha_{g,f}^{-1}h^*) \bigl( \pi_\sigma^f(k,\pi_{\sigma_f^k}^g(\bar l,\bar c)) \bigr)= \alpha_{h,g\circ f}^{-1} \bigl( \pi_\sigma^{g\circ f}((k,\bar l),\bar c) \bigr)\label{eq:assoc_five}\\ 
 &= \alpha_{h,g\circ f}^{-1} \bigl( \pi_\sigma^{g\circ f}((k,\bar l), \pi_{\sigma_{g\circ f}^{(k,\bar l)}}^h(\bar m,\bar a)) \bigr)= \pi_\sigma^{h\circ(g\circ f)}(((k,\bar l),\bar m),\bar a)\\ 
 &= a_{h,g,f,M}^* \bigl(\pi_\sigma^{(h\circ g)\circ f}(p,\bar a)\bigr)\label{eq:assoc_seven}. 
 \end{align} Here, equalities \eqref{eq:assoc_two}, \eqref{eq:assoc_five}, and \eqref{eq:assoc_seven} use the componentwise structure of the projection maps, while the others apply the definitions of the involved maps.
 \end{proof}
\subsection{Vertical and Horizontal Composition of \texorpdfstring{$i$}{i}-Morphisms}
Having defined the operations we need on $1$-morphisms, that is, on piecewise interpretations, we can now move on to define the set of $2$-morphisms for each pair of parallel interpretations $f,g:T_1\to T_2$ by $\Thpw(f,g)=\{[\eta:f\Rightarrow g]_{\equiv_{T_2}}\}$. Such a $2$-morphism is just an equivalence class of $i$-morphisms under provable equivalence (\cref{def:i-morphisms_prov_equiv}). With these $2$-morphisms providing the required morphisms, the set of all interpretations between two theories can be equipped with the structure of a category: The identity $i$-morphism $\id_f:f\Rightarrow f$ is the collection $(\id_f)_\sigma:P_\sigma^f\times P_\sigma^f\to \mathcal{L}_{2}$ defined as follows, where $\bar u:\sigma_f^k$ and $\bar v:\sigma_f^l$:
\begin{equation}
            (\id_f)_\sigma^{k,l}(\bar u,\bar v) := D_{\sigma}^{f,k}(\bar u)\land D_\sigma^{f,l}(\bar v) \land(=_\sigma)_f^{k,l}(\bar u,\bar v).
\end{equation}
It follows from \cref{prop:piece_bijection_morphism} with $b$ the identity that this is indeed an $i$-morphism. Composition for the category $\Thpw(T_1,T_2)$ is given by the following definition:
\begin{definition}[Vertical Composition of $i$-Morphisms]\label{def:i-morphisms_vert_comp}
Let $\eta: g\Rightarrow h$ and $\zeta: f\Rightarrow g$ be $i$-morphisms, where $f,g,h: T_1\to T_2$ are interpretations. Their vertical composition $\eta\bullet\zeta: f\Rightarrow h$ is defined, for each sort $\sigma\in\mathcal{S}_1$ and pieces $k\in P_\sigma^f$, $l\in P_\sigma^h$, by
\begin{equation}
  (\eta\bullet\zeta)_\sigma^{k,l}(\bar{x},\bar{y})
  \;:=\;
  \bigvee_{m\in P_\sigma^g}\exists\bar{z}\,
  \bigl(\zeta_\sigma^{k,m}(\bar{x},\bar{z})
  \;\land\;
  \eta_\sigma^{m,l}(\bar{z},\bar{y})\bigr),
\end{equation}
with $\bar{x}:\sigma_f^k$, $\bar{y}:\sigma_h^l$, and $\bar{z}:\sigma_g^m$.
\end{definition} 
Vertical composition descends to equivalence classes of
$i$-morphisms, is associative, and has $\id_f$ as identity. All this
follows from local faithfulness (\cref{cor:local_faithfulness})
together with the following proposition, which shows that
$\eta\bullet\zeta$ is an $i$-morphism whose semantic dual is the
vertical composite of the duals:
\begin{proposition}[Semantic Dual of Vertical Composition]\label{prop:i-morphism_vert_comp}
    Let $\eta: g\Rightarrow h,\zeta:f\Rightarrow g$ be two $i$-morphisms, where $f,g,h:T_1\to T_2$ are interpretations. Then $\eta\bullet\zeta:f\Rightarrow h$ is an $i$-morphism satisfying $(\eta\bullet\zeta)^*=\eta^*\bullet\zeta^*$, where $\bullet$ denotes vertical composition of natural transformations. In particular, the component at each model is the usual composition of elementary embeddings.
\end{proposition}
\begin{proof}
    Condition \eqref{Dom/Range} for the pre-$i$-morphism $\eta\bullet\zeta$
    follows directly from the corresponding \eqref{Dom/Range} condition for $\zeta$ and $\eta$. By
    \cref{prop:semantic-recognition}, it therefore suffices to
    verify \eqref{eq:dual} for $\eta\bullet\zeta$ and the natural
    transformation $\eta^*\bullet\zeta^*:f^*\to h^*$. So let $M\vDash T_2$,
    $(k,\bar a)\in P_\sigma^f\times M(D_\sigma^{f,k})$, and
    $(l,\bar b)\in P_\sigma^h\times M(D_\sigma^{h,l})$. By definition, \eqref{Dom/Range} for $\eta$ (or $\zeta$), and \eqref{eq:dual} for $\eta$ and $\zeta$, $\bar a\times \bar b\in M((\eta\bullet\zeta)_\sigma^{k,l})$ iff there are $(m,\bar c)\in P_\sigma^g\times M(D_\sigma^{g,m})$ such that $\zeta^*_M(\pi_\sigma^f(k,\bar a))=\pi_\sigma^g(m,\bar c)$ and
    $\eta^*_M(\pi_\sigma^g(m,\bar c))=\pi_\sigma^h(l,\bar b)$. This is the case iff $\eta^*_M\bigl(\zeta^*_M(\pi_\sigma^f(k,\bar a))\bigr)=\pi_\sigma^h(l,\bar b)$, since $\zeta^*_M(\pi_\sigma^f(k,\bar a))=\pi_\sigma^g(m,\bar c)$ for some such $(m,\bar c)$ by surjectivity of $\pi_\sigma^g$.  
\end{proof}

\noindent
A straightforward consequence of this is that what we have called $i$-isomorphisms in \cref{def:i-iso} before are indeed exactly the $i$-morphisms invertible with respect to vertical composition. Summing up, for each pair of theories, there is a $1$-category $\Thpw(T_1,T_2)$ whose objects are piecewise interpretations and whose morphisms are $i$-morphisms between them, both quotiented with respect to provable equivalence. 

\medskip
\noindent
As a final ingredient, necessary to ensure that composition of piecewise interpretations can be seen as a functor $\Thpw(T_2,T_3)\times \Thpw(T_1,T_2)\to \Thpw(T_1,T_3)$, we need another notion of composition of $i$-morphisms, namely horizontal composition. Horizontal composition is a function $\star:\Thpw(f_2,g_2)\times \Thpw(f_1,g_1)\to \Thpw(f_2\circ f_1, g_2\circ g_1)$, where $f_1,g_1:T_1\to T_2$ and $f_2,g_2:T_2\to T_3$, satisfying certain coherence conditions. 
\begin{definition}[Horizontal Composition]\label{def:i-morphisms_hor_comp}
Let $\eta: f_2\Rightarrow g_2$ and $\zeta: f_1\Rightarrow g_1$ be
$i$-morphisms, where $f_1,g_1: T_1\to T_2$ and
$f_2,g_2: T_2\to T_3$.  Their horizontal composition
\begin{equation}
    \eta\star\zeta: f_2\circ f_1\Rightarrow g_2\circ g_1
\end{equation}
is defined, for each
sort $\sigma\in\mathcal{S}_1$ and pieces
$m_1=(k_1,\bar l_1)\in P_\sigma^{f_2\circ f_1}$,
$m_2=(k_2,\bar l_2)\in P_\sigma^{g_2\circ g_1}$, as follows.
\begin{align}
  &(\eta\star\zeta)_\sigma^{m_1, m_2}(\bar x,\bar y)
  \;:=\; D_\sigma^{f_2\circ f_1,m_1}(\bar x)
        \;\land\; D_\sigma^{g_2\circ g_1,m_2}(\bar y) \notag\\
  &\land\;\bigvee_{m_3\in P_\sigma^{g_2\circ g_1}}
    \exists\bar u\;\Bigg(
    \bigvee_{\bar k_4\in P_{\sigma_{g_1}^{k_3}}^{f_2}}
      \exists\bar v\; \bigl((\zeta_\sigma^{k_1,k_3})_{f_2}^{\bar l_1\times\bar k_4}(\bar x,\bar v)
        \;\land\;
        \eta_{\sigma_{g_1}^{k_3}}^{\bar k_4,\bar l_3}(\bar v,\bar u)
      \bigr)
  \land\;
    (=_\sigma)_{g_2\circ g_1}^{m_2,m_3}(\bar y,\bar u)
    \Bigg),
\end{align}
where $m_3=(k_3,\bar l_3)$ and $\bar x:\sigma_{f_2\circ f_1}^{m_1}$,
$\bar y:\sigma_{g_2\circ g_1}^{m_2}$, $\bar u:\sigma_{g_2\circ g_1}^{m_3}$, and $\bar v:\sigma_{f_2\circ g_1}^{(k_3,\bar k_4)}$. By construction and \cref{lem:variables_translations}, $(\eta\star\zeta)_\sigma^{m_1,m_2}:\sigma^{m_1}_{f_2\circ f_1}\times
\sigma^{m_2}_{g_2\circ g_1}$, as required of a pre-$i$-morphism.
\end{definition}
As with vertical composition, well-definedness is best established
through the semantic dual:
\begin{proposition}[Semantic Dual of Horizontal Composition]\label{prop:hor-comp-dual}
    Let $\eta: f_2\Rightarrow g_2$ and $\zeta: f_1\Rightarrow g_1$
    be $i$-morphisms, where $f_1,g_1: T_1\to T_2$ and
    $f_2,g_2: T_2\to T_3$ are interpretations. Then
    $\eta\star\zeta: f_2\circ f_1\Rightarrow g_2\circ g_1$ is an
    $i$-morphism satisfying
    \begin{equation}
        (\eta\star\zeta)^*
        =\alpha_{g_2,g_1}^{-1}\bullet(\zeta^*\star\eta^*)\bullet\alpha_{f_2,f_1},
    \end{equation}
    where $\zeta^*\star\eta^*$ denotes the horizontal composition of natural
    transformations. At a model $M\vDash T_3$, the latter can, by general category theory, be described as
    \begin{equation}\label{eq:hor-nat}
        (\zeta^*\star\eta^*)_M
        =\zeta^*_{g_2^*M}\circ f_1^*(\eta^*_M)
        =g_1^*(\eta^*_M)\circ\zeta^*_{f_2^*M}.
    \end{equation}
\end{proposition}

\begin{proof}
    Condition \eqref{Dom/Range} for the pre-$i$-morphism $\eta\star\zeta$
    holds by construction, and
    $F:=\alpha_{g_2,g_1}^{-1}\bullet(\zeta^*\star\eta^*)\bullet\alpha_{f_2,f_1}$
    is a composite of natural transformations, so by
    \cref{prop:semantic-recognition} it suffices to verify
    \eqref{eq:dual} for $\eta\star\zeta$ and $F$. So let $M\vDash T_3$,
    $(m_1,\bar a)\in P_\sigma^{f_2\circ f_1}\times
    M(D_\sigma^{f_2\circ f_1,m_1})$, and
    $(m_2,\bar b)\in P_\sigma^{g_2\circ g_1}\times
    M(D_\sigma^{g_2\circ g_1,m_2})$, and let
    $m_1=(k_1,\bar l_1)$.

    \smallskip\noindent
    First by definition, and two applications of \eqref{eq:dual}, which requires establishing the respective domain conditions via \eqref{Dom/Range} for $\eta$ at $M$, $\zeta$ at $f_2^*M$ via \cref{cor:extensions_elementwise}, and \eqref{Element-Embed} for $\eta$ for $D_\sigma^{g_1,k_3}$, we observe that $\bar a\times \bar b\in M((\eta\star\zeta)_\sigma^{m_1,m_2})$
    iff there are $m_3=(k_3,\bar l_3)\in P_\sigma^{g_2\circ g_1}$,
    $\bar k_4\in P_{\sigma_{g_1}^{k_3}}^{f_2}$,
    $\bar u\in M(D_\sigma^{g_2\circ g_1,m_3})$, and
    $\bar v\in M(D_\sigma^{f_2\circ g_1,(k_3,\bar k_4)})$ such that
    \begin{equation}\label{eq:hor-images}
      \zeta^*_{f_2^*M}\bigl(\pi_\sigma^{f_1}(k_1,\,
        \pi_{\sigma_{f_1}^{k_1}}^{f_2}(\bar l_1,\bar a))\bigr)
        =\pi_\sigma^{g_1}\bigl(k_3,\,
        \pi_{\sigma_{g_1}^{k_3}}^{f_2}(\bar k_4,\bar v)\bigr),
      \qquad
      \eta^*_M\bigl(\pi_{\sigma_{g_1}^{k_3}}^{f_2}(\bar k_4,\bar v)\bigr)
        =\pi_{\sigma_{g_1}^{k_3}}^{g_2}(\bar l_3,\bar u),
    \end{equation}
    and $\pi_\sigma^{g_2\circ g_1}(m_3,\bar u)
    =\pi_\sigma^{g_2\circ g_1}(m_2,\bar b)$. 

    \smallskip\noindent
    To see that this is the case iff 
    $F_M\bigl(\pi_\sigma^{f_2\circ f_1}(m_1,\bar a)\bigr)
    =\pi_\sigma^{g_2\circ g_1}(m_2,\bar b)$, for the forwards direction, one plugs these equalities into the computation of $F_M$. For the backwards direction, one chooses representatives of the intermediate steps, again observing that these satisfy the respective domain conditions using elementarity and \cref{cor:extensions_elementwise}.
\end{proof}
\noindent
Horizontal composition is not associative on the nose, but this is to be expected, since composition of piecewise interpretations is itself only associative up to the associator.

\subsection{Bicategory and Mod-Pseudofunctor}
Having defined the basic operations on piecewise interpretations and $i$-morphisms and verified that they are invariant under provable equivalence, we now have all the basic ingredients for the following central definition:
\begin{definition}[Bicategory of Theories and Piecewise Interpretations]\label{def:bicategory_pw_int}
    Let $\Thpw$ denote the bicategory whose objects are given by classical first-order theories, sets of $1$-morphisms $\Thpw(T_1,T_2)$ by equivalence classes of piecewise interpretations $f:T_1\to T_2$ under provable equivalence, and whose sets of $2$-morphisms $\Thpw(f,g)$ are given by equivalence classes of $i$-morphisms $\eta:f\Rightarrow g$ under provable equivalence. The identities are the identity interpretations and $i$-morphisms, composition of $1$-morphisms is given by composition of piecewise interpretations, and vertical and horizontal composition of $2$-morphisms is given by vertical and horizontal composition of $i$-morphisms. Unitors and associators are given by unitors and associators for piecewise interpretations.
\end{definition}
\begin{remark}[Size of $\Thpw$]
By design, the bicategory $\Thpw$ is locally small, although its objects form a proper class. Allowing theories in languages of larger cardinalities than the infinite cardinal $\kappa$ fixed in the introduction only adds more objects (provided one adequately enlarges variable supplies), without essentially changing the sets of $1$- or $2$-morphisms between them.
\end{remark}
At this stage, we have not verified the coherence conditions one needs for a bicategory: naturality of the associators and unitors, the interchange law, and so on. However, we do have all the data one needs for a bicategory, and we have semantic duals for theories, namely their categories of models, and for interpretations and $i$-morphisms. Moreover, we have appropriate unitors and compositors for these duals. So we have what one might call a pre-bicategory ($\Thpw$) and something which looks like a locally faithful contravariant pseudofunctor ($\Mod$) into a strict $2$-category ($\mathbf{CAT}$). In this situation, we can apply the following general theorem, whose proof makes use of ideas well known from coherence results for bicategories and can be found in Appendix~\ref{app:bicategorical_results}:
\begin{restatable}[Reflection Theorem]{theorem}{reflection}\label{thm:reflection}
Let $\mathcal{B}$ be a pre-bicategory and $\mathcal{C}$ be a strict $2$-category. Assume that the associators $a$ and the unitors $l$ and $r$ of $\mathcal{B}$ are isomorphisms and that there is a locally faithful contravariant pseudofunctor $F:\mathcal{B}\to \mathcal{C}$. Then $\mathcal{B}$ is a bicategory and $F:\mathcal{B}\to \mathcal{C}$ is a locally faithful contravariant pseudofunctor between bicategories.
\end{restatable}
\noindent
Applying this to our situation yields the following central theorem, justifying the preceding definition and encapsulating the work we have done on semantic duals.
\begin{theorem}[Mod-Functor for Piecewise Interpretations]\label{mod-functor}
    There is a locally faithful, contravariant pseudofunctor between bicategories $\Mod:\Thpw\to \mathbf{CAT}$, sending a theory $T$ to its category of models $\Mod(T)$, a piecewise interpretation $f:T_1\to T_2$ to its semantic dual $f^*:\Mod(T_2)\to \Mod(T_1)$, and an $i$-morphism $\eta:f\Rightarrow g$ to its semantic dual $\eta^*:f^*\Rightarrow g^*$. The compositor for $f$ and $g$ is given by the inverse of $\alpha_{g,f}$, and the unitor for $T$ is given by $\alpha_{T}$.
\end{theorem}
\begin{proof}
    That $\Thpw$ is a pre-bicategory with well-defined operations was established in the remarks above, and that the associators and unitors are isomorphisms is established by \cref{prop:comp_assoc} and the remarks under \cref{def:unitor}, respectively. To see that $\Mod$ indeed defines a locally faithful, contravariant pseudofunctor $\Mod:\Thpw\to \mathbf{CAT}$, note first that the required structure, in particular the unitor and compositor of the semantic dual, are well-defined on equivalence classes under provable equivalence, as a consequence of \cref{prop:prov_equiv_char}.  Functoriality on each hom-category (i.e., $(\eta\bullet\zeta)^*=\eta^*\bullet\zeta^*$ and $(\id_f)^*=\id_{f^*}$) was established in \cref{prop:i-morphism_vert_comp} and is immediate from the definition of $\id_f$, respectively. Continuing, the compositor and unitor are natural isomorphisms, cohering appropriately with the unitors and associators of $\Thpw$ (\cref{prop:alpha_gf,prop:hor-comp-dual,prop:associator_coherence,lem:unitor-coherence}). Finally, local faithfulness is the content of \cref{cor:local_faithfulness}, allowing us to apply \cref{thm:reflection}.
\end{proof}
\noindent
Let us now consider equivalence in this bicategory. Two theories $T_1$ and $T_2$ are equivalent in $\Thpw$ iff there is a pair of piecewise interpretations $f:T_1\to T_2$ and $g:T_2\to T_1$ whose compositions are isomorphic to the respective identities, that is, $g\circ f\cong \id_{T_1}$ and $f\circ g\cong \id_{T_2}$. Such theories will be called \textit{piecewise bi-interpretable}. 

\medskip
\noindent
Of course, if two theories are bi-interpretable in the ordinary sense, they are piecewise bi-interpretable. The converse does not hold in general, as we will see in \cref{sec:application}. Because pseudofunctors preserve equivalence, we also have the following immediate corollary of \cref{mod-functor}.  
\begin{corollary}[Piecewise Bi-Interpretability Implies Equivalence of Categories of Models]
    If two theories are piecewise bi-interpretable, they have equivalent categories of models.
\end{corollary}
\noindent
The converse, likewise, does not hold: A counterexample from \citet[Theorem~5.7]{barrett2016morita}, demonstrating that categorical equivalence does not imply Morita equivalence, can be repurposed to show this: essentially, one notes that piecewise bi-interpretable theories have isomorphic sentential algebras, and that one of the theories in their example has a sentential algebra with atoms, while the other does not.\footnote{See \cite{meadows2025halvorson} for a thorough discussion of these and related examples.} Moreover, we see that for $n\neq m$, the theories $\mathrm{E}_n$ and $\mathrm{E}_m$ from \cref{ex:n-element} are mutually piecewise interpretable, but not piecewise bi-interpretable, because their respective models have non-isomorphic automorphism groups, whence they cannot have equivalent categories of models. 

\section{A Functorial Reduction of (Finitely) Many-Sorted Logic to Single-Sorted Logic}\label{sec:application}
 The formalization of many-sorted first-order logic came in tandem with proposed reductions of many-sorted to single-sorted logic \citep{schmidt1938deduktive, schmidt1951zulassigkeit, wang_logic_1952}. Logicians quickly arrived at a procedure, hereafter called \textit{unification}, which takes a many-sorted theory and returns a corresponding single-sorted theory. The idea behind this procedure is quite simple: use predicate symbols instead of sorts and relativize quantifiers accordingly.\footnote{See \cite{oddsson2026constructive} for investigations into the metatheory of this procedure, and \cite{oberschelp1962untersuchungen} for a more abstract and general perspective on many-sorted logic in a single-sorted setting.} Quine was a proponent of the idea that this procedure shows that many-sorted logic is ultimately dispensable \citep[see, e.g.,][]{quine1956unification}. We will follow \citet{barrett2017quine} in their reconstruction of this idea. First, for any many-sorted language $\mathcal{L}$ with finitely many sorts $\sigma_1,\ldots,\sigma_n$, there is a single-sorted language $\mathcal{L}^u$ with a single sort $\tau$, for each sort $\sigma_k$ of $\mathcal{L}$, a unary predicate $Q_k$, and, for each non-identity relation symbol $R:\bar \sigma$, a relation symbol $R^u:\tau^{|\bar\sigma|}$. Moreover, there is a natural translation  $r_\mathcal{L}:\mathcal{L}\to \mathcal{L}^u$ with one piece $0$ for each sort, sending every sort $\sigma_k$ to the single sort $\tau$ with domain formula $Q_k(x)$, every non-identity relation symbol $R:\bar \sigma$ with $|\bar \sigma|=m$ to $R^u(\bar y)$ with $\bar y:\tau^m$ duplicate-free, and identity to identity. Using this, we can define:   
 \begin{definition}[Unification of a Theory]
        Let $T$ be a many-sorted first-order theory in a language $\mathcal{L}$ with finitely many sorts $\mathcal{S}_T=\{\sigma_1,\ldots,\sigma_n\}$. The \emph{unification} $T^u$ of $T$ is the single-sorted first-order theory defined as follows. The axioms of $T^u$ are\footnote{Strictly speaking, the second condition of \eqref{TU_cover} is superfluous, representing translation of $\exists_{\sigma_i} x (x=_{\sigma_i}x)$ for each $i$. We include it so that one may restrict \eqref{TU_rel} to axioms of $T$.}:
        \begin{align}
            &\delta_{r_\mathcal{L}} &&\text{ for each }\delta\in T\tag{Relativization}\label{TU_rel} \\ 
            &\forall_\tau x\; (\bigvee_{i=1}^n Q_i(x)) \quad \text{and}\quad \exists_\tau x\;( Q_i(x))  &&(1\leq i\leq n)\tag{Cover}\label{TU_cover}\\
            &\forall_\tau x\; \neg(Q_i(x)\land Q_j(x))  &&(i\neq j)\tag{Disjoint}\label{TU_dis}\\
            &\forall_\tau x_1\ldots\forall_\tau x_m\; (R^u(\bar x)\to \bigwedge_{r=1}^m Q_{j_r}(x_r)) &&(R:\sigma_{j_1}\times\cdots\times \sigma_{j_m})\tag{Sorting}\label{TU_sort}
        \end{align}
        \end{definition}
    Barrett and Halvorson show that $T$ and $T^u$ are Morita equivalent for any theory $T$ with finitely many sorts, taking this to be a partial vindication of Quine's idea about many-sorted logic. We will show that $T$ and $T^u$ are piecewise bi-interpretable, thereby demonstrating that they are equivalent in another sense, and also address a problem for unification raised by \citet{hook1985note}: the passage from $T$ to $T^u$ need not preserve interpretability. For instance, let $T_{\mathrm{uno}}$ say that its single sort has exactly one element, and let $T_{\mathrm{dos}}$ have two sorts, each with exactly one element. Then $T_{\mathrm{dos}}$ is interpretable in $T_{\mathrm{uno}}$, but $T_{\mathrm{dos}}^u$ is not interpretable in $T_{\mathrm{uno}}^u$ by an ordinary interpretation, since every structure ordinarily interpretable in a one-element structure has at most one element.\footnote{Hook's example is different and somewhat more interesting, but he only considers one-dimensional interpretations. Compare the example chosen here with \citet[Example~3.13]{visser2021extension}. Relatedly, \citet[Theorem~6.4]{meadows2024beyond} would also be interesting to examine from the perspective of piecewise interpretations.} So unification of domains does not even preserve the interpretability relation, let alone act functorially on interpretations. With piecewise interpretations, both defects disappear (\cref{th:biequivalence_ms_and_ss} and \cref{cor:uni_pres} below). In the following, note that one really only needs the piecewise structure and relativization, and not other features of interpretation, such as modification of identity or many-dimensionality, to achieve the desired goal.
    
    \medskip\noindent
    Let us proceed by constructing the relevant interpretations: By construction, $r_\mathcal{L}$ induces an interpretation $r_T:T\to T^u$. We can also furnish a piecewise interpretation $d_T:T^u\to T$, with set of pieces $[n]$ for $\tau$, with $\tau_{d_T}^k:=\sigma_k$ and $D_\tau^{d_T,k}:=\top_{\sigma_k}$. For the identity and the sortal predicates, we set
    \begin{align}
        (=_\tau)_{d_T}^{i,j}(x,y)&:=
    \begin{cases}
        x=_{\sigma_i}y & \text{if }i=j\\
        \bot_{x\times y} & \text{if }i\neq j
    \end{cases} \quad \text{and}\quad  
    (Q_{k})_{d_T}^{i}(x):=
    \begin{cases}
        \top_x & \text{if }i=k\\
        \bot_x & \text{if }i\neq k
    \end{cases}\\
    (R^u)_{d_T}^{i_1,\ldots,i_m}(\bar x):&=
    \begin{cases}
        R(\bar x)
            & \text{if }(i_1,\ldots,i_m)
                        =(j_1,\ldots,j_m)\\
        \bot_{\bar x}
            & \text{otherwise}
    \end{cases} (R:\sigma_{j_1}\times\cdots\times\sigma_{j_m})
    \end{align}
One can easily verify that $d_T$ respects identity with respect to $T$, noting that \eqref{eq: f-identity-elim} for $d_T$ reduces to \eqref{id_elim} in $\mathcal{L}$. Let us now describe the semantic duals of both translations:

\medskip\noindent
Let $N\vDash T^u$. Then $r_T^*N(\sigma_k)=\{(0,a)\mid a\in N(Q_k)\}/\sim_k$, where $(0,a)\sim_k(0,b)$ iff $a=_\tau b$. For a relation symbol $R:\sigma_{i_1}\times \cdots\times \sigma_{i_m}$,
\begin{equation}
   r_T^*N(R)
=
\left\{
\bigl([0,a_1]_{\sim_{i_1}},\ldots,
      [0,a_m]_{\sim_{i_m}}\bigr)
\;\middle|\;
\begin{array}{l}
a_r\in N(Q_{i_r})\quad(1\le r\le m),\\
(a_1,\ldots,a_m)\in N(R^u)
\end{array}
\right\}.
\end{equation}

\medskip\noindent
Turning to $d_T$ and a model $M\vDash T$, we get $d_T^*M(\tau)
=
\bigl\{(k,a)\mid k\in[n],\ a\in M(\sigma_k)\bigr\}/\!\sim$, where $ (i,a)\sim (j,b)$ iff $i=j$ and $a=_{\sigma_i}b$. The new sortal predicates $Q_k$ receive extensions $d_T^*M(Q_k) = \bigl\{[k,a]_\sim\mid a\in M(\sigma_k)\bigr\}$. Finally, a relation symbol $R^u:\sigma_{i_1}\times\cdots\times\sigma_{i_m}$ receives the extension 
\begin{equation}
d_T^*M(R^u)
= \left\{ \bigl([i_1,a_1]_\sim,\ldots, [i_m,a_m]_\sim\bigr) \;\middle|\;(a_1,\ldots,a_m)\in M(R)
\right\}.
\end{equation}

\begin{proposition}[Piecewise Bi-interpretability of $T$ and $T^u$]\label{prop:T_and_TU}
        Let $T$ be a first-order theory in a language $\mathcal{L}$ with finitely many sorts. Then $T$ and $T^u$ are piecewise bi-interpretable via the interpretations $r_T:T\to T^u$ and $d_T:T^u\to T$.
\end{proposition}
\begin{proof}
We have yet to verify that $d_T$ really is an interpretation. To that end, note that $d_T\circ r_T$ respects identity and can hence be seen to be an interpretation $L(\Sigma_T)\to T$, as can $\id_T$. We now provide an $i$-isomorphism $\theta:d_T\circ r_T \Rightarrow\id_T$ (with respect to $T$) and its semantic dual, where one shows that it is indeed an $i$-isomorphism with said dual via \cref{prop:semantic-recognition}:
\begin{equation}
   \theta_{\sigma_k}^{(0,i),0}(x,y):=
\begin{cases}
x=_{\sigma_k}y & i=k\\
\bot_{x\times y} & i\ne k
\end{cases}, \quad 
\theta^*_{M,\sigma_k}:
   \begin{aligned}
(d_T\circ r_T)^*M(\sigma_k)
&\longrightarrow (\id_T^*M)(\sigma_k)\\
[(0,k),a]&\longmapsto[0,a]
\end{aligned}
\end{equation}
From Propositions \ref{prop:criterion_interpretationhood} and \ref{composition_complex}, we get that $T\vdash (\delta_{r_T})_{d_T}\leftrightarrow \delta$ for all $\mathcal{L}$-sentences $\delta$, whence $T\vdash (\delta_{r_T})_{d_T}$ for all $\delta$ such that $T\vdash \delta$. It remains to check that $d_T$ preserves the axioms \eqref{TU_cover}, \eqref{TU_dis}, and \eqref{TU_sort}. Up to logical equivalence the translations of these axioms under $d_T$ are, using that its domain formulae are just the respective vera:
\begin{align}
  &\bigwedge_{h\in[n]}\forall_{\sigma_h}x\,\bigvee_{k=1}^n(Q_k)^h_{d_T}(x), \quad \bigvee_{h\in[n]}\exists_{\sigma_h}x\,(Q_k)^h_{d_T}(x)
  \tag{Cover$_{d}$}\label{eq:cover-d}\\
  &\bigwedge_{h\in[n]}\forall_{\sigma_h}x\,\neg\bigl((Q_k)^h_{d_T}(x)\land(Q_l)^h_{d_T}(x)\bigr)\quad(k\neq l)
  \tag{Disjoint$_{d}$}\label{eq:disjoint-d}\\
  &\bigwedge_{\bar h\in[n]^m}\forall\bar x\,\bigl((R^u)^{\bar h}_{d_T}(\bar x)\to\textstyle\bigwedge_{r=1}^m(Q_{j_r})^{h_r}_{d_T}(x_r)\bigr)
  \tag{Sorting$_{d}$}\label{eq:sorting-d}
\end{align}
where in \eqref{eq:sorting-d}, $R:\sigma_{j_1}\times\cdots\times\sigma_{j_m}$ and $\bar x=(x_1,\ldots,x_m):\sigma_{h_1}\times\cdots\times\sigma_{h_m}$ is duplicate-free. Since $(Q_k)^h_{d_T}$ is $\top_{x}$ if $h=k$, the disjunct with $k=h$ makes each conjunct of the first formula in \eqref{eq:cover-d} provable; for the second, the disjunct with $h=k$ is $\exists_{\sigma_k}x\,(x=_{\sigma_k}x)$, which is provable by assumption on the calculus; and in each conjunct of \eqref{eq:disjoint-d} at least one of $(Q_k)^h_{d_T}$, $(Q_l)^h_{d_T}$ is $\bot_{x}$. Finally, in \eqref{eq:sorting-d}, if $\bar h=\bar j$, the consequent is a conjunction of instances of $\top_{x_r}$, and otherwise $(R^u)^{\bar h}_{d_T}=\bot_{\bar x}$. Therefore, $d_T:T^u\to T$ is an interpretation.

\noindent
Now, we produce the required $i$-isomorphism $\zeta:r_T\circ d_T
   \Rightarrow\id_{T^u}$, again using \cref{prop:semantic-recognition}:
   \begin{equation}    
\zeta_\tau^{(i,0),0}(x,y):=
Q_i(x)\land x=_\tau y, \quad 
\zeta^*_{N,\tau}:
\begin{aligned}
(r_T\circ d_T)^*N(\tau)
&\longrightarrow(\id_{T^u}^*N)(\tau)\\
[(i,0),a]&\longmapsto[0,a].
\end{aligned}
\end{equation}
\end{proof}
\noindent
Note that this, in particular, shows that $T_{\mathrm{dos}}$ and $T_{\mathrm{dos}}^u$ are piecewise bi-interpretable, but not bi-interpretable without pieces. The assignment $T\mapsto T^u$ can be seen as part of a biequivalence of bicategories, which we will show in a slightly indirect fashion, to save us some work. To that end, consider two full sub-bicategories of $\Thpw$, namely $\Thpw_{\mathrm{fin}}$, the bicategory of theories with finitely many sorts, and $\Thpw_{\mathrm{ss}}$, the bicategory of theories with a single sort. 
\begin{theorem}[Unification is Part of a Biequivalence]\label{th:biequivalence_ms_and_ss}
    The inclusion pseudofunctor $i:\Thpw_{\mathrm{ss}}\to \Thpw_{\mathrm{fin}}$ is part of a biequivalence, with a quasi-inverse sending a theory $T$ to $T^u$ and an interpretation $f:S\to T$ to $r_T\circ (f\circ d_S):S^u\to T^u$.
\end{theorem}
\begin{proof}
    Using choice, a pseudofunctor is part of a biequivalence iff it is an equivalence locally, that is, on each hom-category, and biessentially surjective, that is, surjective up to equivalence in the target bicategory \citep[see, for instance,][Lemma 3.1]{gurski2011biequivalences}.\footnote{In this particular case, choice is not indispensable, since we have explicit choices for the relevant pieces of data.} It is clear that the inclusion pseudofunctor is a local equivalence, and it is biessentially surjective by \cref{prop:T_and_TU}. That we can choose an inverse as above is part of the proof of the characterization of biequivalences.
\end{proof}
\begin{corollary}[Preservation of Interpretability under Unification]\label{cor:uni_pres}
    Let $S$ and $T$ be first-order theories with finitely many sorts. Then $S$ is piecewise interpretable in $T$ iff $S^u$ is piecewise interpretable in $T^u$.
\end{corollary}
\noindent
Therefore, if one takes piecewise interpretations as one's notion of morphism between theories, there is a (pseudo)functorial reduction of many-sorted to single-sorted logic. This vindicates the idea that sorts are, at least from the point of view of the theory of interpretations, ultimately dispensable, if one allows for only finitely many sorts.\footnote{See also \citet{halvorson2024against} on the issue of the dispensability of sorts.} 
\section{Conclusion \& Future Work}
In this paper, we have shown that piecewise interpretations are a well-behaved notion of interpretation between many-sorted first-order theories and developed their basic theory by organizing them into a bicategory. The developed theory of semantic duals is encapsulated by a locally faithful, contravariant pseudofunctor on that bicategory. As an application, we established a biequivalence between the bicategories of finitely many-sorted and single-sorted theories. In particular, we showed that unification, that is, replacing sorts by sortal predicates on a single sort, preserves piecewise interpretability.

\medskip
\noindent
In future work, we will first investigate connections to categorical logic, showing that $\Thpw$ is biequivalent to a $2$-category of (suitably small) Boolean pretopoi, thereby allowing for the transfer of results from categorical logic. Following that, we will explain the relationship between ordinary and piecewise interpretations with a bicategorical analogue of Visser's Kleisli construction \citeyearpar[Section~3.8.8]{visser2021extension}, and use this construction to show that piecewise bi-interpretability coincides with Morita equivalence, appropriately defined \citep[regarding the relationships between these notions, see][]{mceldowney2020morita, meadows2024beyond, li2025theoretical}.

\paragraph{Acknowledgements}
Funded by the European Union (ERC, FORMALISM, 101044114). Views and opinions expressed are, however, those of the author(s) only and do not necessarily reflect those of the European Union or the European Research Council Executive Agency. Neither the European Union nor the granting authority can be held responsible for them. The author would like to thank Albert Visser and Mirko Engler for helpful feedback on various parts of the project. 
\paragraph{AI Statement}
The author used Claude Opus and ChatGPT Astra to help find relevant literature on
the topic and to assist with developing proofs, drafting text, and identifying typos. When parts of an argument or text were suggested by an LLM, great care was taken by the author to check them, render them understandable, and present them in accordance with the overall structure and conventions of the remainder of the document. The author takes full responsibility for the content of the paper.

\bibliography{references}
\begin{appendices}
\makeatletter

\renewcommand{\@seccntformat}[1]{\appendixname~\csname the#1\endcsname: }
\makeatother
\crefname{section}{Appendix}{appendices}
\Crefname{section}{Appendix}{Appendices}
\section{Canonical Context and Substitution}\label{app:cc_and_subst}

\begin{definition}[Canonical Context and Arity]
For a sequence $\bar{s}$, write $\bar{s}^\circ$ for the subsequence obtained by deleting entries occurring more than once from right to left, retaining only their first occurrences, which we will also call its deduplication. Note that $(\bar{a}^\circ \times \bar{b}^\circ)^\circ = (\bar{a}\times\bar{b})^\circ$ for any sequences $\bar{a},\bar{b}$. The \emph{canonical context} $\bar{x}(\varphi)$ of a formula $\varphi$ is defined inductively:
\begin{enumerate}[label=(\roman*)]
\item \emph{Truth Values:} $\bar x(\top)=\bar x(\bot)=\varepsilon$.
    \item \emph{Atomic:} $\bar{x}(R(z_1,\ldots,z_n)) = (z_1,\ldots,z_n)^\circ$.
    \item \emph{Binary:} $\bar{x}(\varphi \star \psi) = (\bar{x}(\varphi)\times\bar{x}(\psi))^\circ$ for $\star \in \{\land, \lor, \to\}$.
    \item \emph{Negation:} $\bar{x}(\neg\varphi) = \bar{x}(\varphi)$.
    \item \emph{Quantifier:} $\bar{x}(\forall z\,\varphi)$ is $\bar{x}(\varphi)$ with $z$ deleted.
\end{enumerate}
The \emph{arity} of $\varphi$ is the sequence of arities of its canonical context, written $\bar{\sigma}(\varphi)$. We also write $\varphi : \bar{\sigma}$. 

\medskip\noindent We will make frequent use of the following map induced by deduplication: If $\bar x$ is a sequence, denote by $d_x:[|\bar x|]\to [|\bar x^\circ|]$ the unique surjection such that $(\bar x^\circ)_{d_x(i)}=x_i$, sending an index $i$ to the index of the first occurrence of $x_i$ in $\bar x^\circ$. We will omit the subscript where the context allows it. The deduplication map $d_x$ induces a function on sequences of the same length as $\bar x^\circ$ to sequences of the length of $\bar x$, via $(\bar y^{d_x})_i:=y_{d_x(i)}$. This function satisfies $(\bar x^\circ)^{d_x}=\bar x$.
\end{definition}
\medskip
\noindent
Recall that, given a context $\bar{y} = (y_1,\ldots,y_m) : \bar{\tau}$ for a formula $\varphi$ with canonical context $\bar{x}(\varphi) = (x_1,\ldots,x_n) : \bar{\sigma}$, there is a unique injection $c : [n] \to [m]$ such that $x_j = y_{c(j)}$ for $1 \leq j \leq n$. We call $c$ the \emph{context injection} of $\bar{y}$ relative to $\varphi$.

\begin{definition}[Simultaneous Substitution]\label{def:sim_subst}
For a formula $\varphi$ with canonical context $\bar{x} : \bar{\sigma}$ and a sequence $\bar{y} : \bar{\sigma}$ of the same arity, the \emph{simultaneous substitution} $\varphi(\bar{y})$ is defined inductively:
\begin{enumerate}[label=(\roman*)]
\item \emph{Truth Values:} $\top(\bar y)=\top$ and $\bot(\bar y)=\bot$.
    \item \emph{Atomic:} $R(\bar z)(\bar{y}) := R(\bar y ^{d_{z}})$.
    \item \emph{Binary:} $(\varphi_1 \star \varphi_2)(\bar{y}) := \varphi_1(\bar{y}^{c_1}) \star \varphi_2(\bar{y}^{c_2})$, where $c_i : [|\bar{x}(\varphi_i)|] \to [|\bar{x}(\varphi_1 \star \varphi_2)|]$ are the context injections for $\bar{x}(\varphi_1\star\varphi_2)$ relative to $\varphi_i$, and $\bar{y}^{c_i}_j := y_{c_i(j)}$.
    \item \emph{Negation:} $(\neg\varphi)(\bar{y}) := \neg(\varphi(\bar{y}))$.
    \item \emph{Quantifier:} Let $\forall_\tau z\,\varphi:\bar\sigma$ and $\bar y:\bar\sigma$.
    Since $\bar x(\forall z\,\varphi)\times z$ is a context for $\varphi$, let
    $c:[|\bar x(\varphi)|]\hookrightarrow[|\bar\sigma|+1]$ be the context injection of
    $\bar x(\forall z\,\varphi)\times z$ relative to $\varphi$. Choose a variable $z':\tau$, fresh for $\bar y$ and $\varphi$, and set
    \begin{equation}
        (\forall z\,\varphi)(\bar y):=\forall z'\,\varphi\bigl((\bar y\times z')^{c}\bigr).
    \end{equation}
\end{enumerate}
\end{definition}

In the quantifier clause of the last definition, different choices of fresh $z'$ will produce $\alpha$-variants of the resulting formulae. This has the consequence that translations, as defined in the paper, are, strictly speaking, only functional up to renaming of bound variables, which will be enough for our purposes. While the following lemmata are essentially bookkeeping and require no special insight, we record them in detail, because the details are a bit finicky. We first need some auxiliary combinatorial lemmata. 

\begin{lemma}[Deduplication under Surjection]\label{lem:dedup-reindex}
Let $\bar{a}_1,\ldots,\bar{a}_n$ be finite sequences and let $s:[m]\twoheadrightarrow[n]$ be a surjection such that $\min s^{-1}(i) < \min s^{-1}(j)$ whenever $i < j$. Then
\begin{equation}
(\bar{a}_{s(1)}\times\cdots\times\bar{a}_{s(m)})^\circ = (\bar{a}_1\times\cdots\times\bar{a}_n)^\circ.
\end{equation}
In particular, if the sequences $\bar{a}_i$ are pairwise disjoint and individually duplicate-free, then $(\bar{a}_{s(1)}\times\cdots\times\bar{a}_{s(m)})^\circ = \bar{a}_1\times\cdots\times\bar{a}_n$. 
\end{lemma}

\medskip
\noindent
In applications, the surjection $s$ will typically be the deduplication map for a sequence other than $\bar a$. One should also note that this covers the case where the concatenated sequences each have length one.
\begin{proof}
Induction on $m\geq 1$. For $m,n\leq1$, the claim is obvious. Now, let $m>1$ and assume the claim holds for all $1\leq m'<m$. If $s(m)=s(j)$ for some $1\leq j\leq m-1$, then the restriction of $s$ to $[m-1]$ satisfies the hypotheses for $[n]$ and $(\bar a_{s(1)}\times\cdots\times \bar a_{s(m-1)}\times \bar a_{s(m)})^\circ= (\bar a_{s(1)}\times\cdots\times \bar a_{s(m-1)})^\circ=(\bar a_1\times\cdots\times\bar a_n)^\circ$ by the induction hypothesis. Otherwise, the hypothesis on $s$ forces $s(m)=n$, whence the restriction of $s$ to $[m-1]$ satisfies the hypothesis for $[n-1]$. Therefore: \begin{align}
    (\bar{a}_{s(1)}\times\cdots\times\bar{a}_{s(m)})^\circ &=((\bar{a}_{s(1)}\times\cdots\times\bar a_{s(m-1)})^\circ\times\bar{a}_{s(m)})^\circ\\
    &=((\bar{a}_{1}\times\cdots\times\bar a_{n-1})^\circ \times\bar{a}_{n})^\circ\\
    &=(\bar a_1\times\cdots\times \bar a_n)^\circ
\end{align}
\end{proof}
\begin{corollary}[Deduplication Map under Surjection]\label{cor:dedup_map_surjection}
    Let $\bar a$ be a sequence of length $n$ and let $s:[m]\twoheadrightarrow[n]$ be a surjection with $\min s^{-1}(i) < \min s^{-1}(j)$ whenever $i < j$. Then 
    \begin{equation}
        d_{a^s}=d_{ a}\circ s
    \end{equation}
\end{corollary}
\noindent
We now record a lemma about the interaction of context injections and deduplications.
\begin{lemma}[Deduplication Map under Context Injection]\label{dedup_map_injection}
Let $\bar y$ be a sequence of length $n$ and let $c:[m]\hookrightarrow[n]$ be an injection, with $\bar y^c$ the induced sequence $(\bar y^c)_i:=y_{c(i)}$. The entries of $(\bar y^c)^\circ$ occur among those of the duplicate-free sequence $\bar y^\circ$, so there is a unique injection $\hat c:[\,|(\bar y^c)^\circ|\,]\hookrightarrow[\,|\bar y^\circ|\,]$ with $(\bar y^\circ)_{\hat c(i)}=((\bar y^c)^\circ)_i$. Then
    \begin{equation}
        d_{y}\circ c = \hat c\circ d_{\bar y^c}.
    \end{equation}
\end{lemma}
\begin{proof}
    Both sides are functions $[m]\to[\,|\bar y^\circ|\,]$. For $i\in[m]$, $(d_{y}\circ c)(i)=d_{y}(c(i))$ is the position of $y_{c(i)}$ in $\bar y^\circ$. On the other side $d_{\bar y^c}(i)$ is the position of $(\bar y^c)_i=y_{c(i)}$ in $(\bar y^c)^\circ$, and $\hat c$ maps it to the position of that same element $y_{c(i)}$ in $\bar y^\circ$. As $\bar y^\circ$ is duplicate-free, this position is unique, so the two maps agree.
\end{proof}
\ccsubst*
\begin{proof}
By induction on $\varphi$. The truth-value and negation cases are immediate, as both substitution and canonical context commute with negation.

\medskip
\noindent
\emph{Atomic:} For any sequence $\bar z$, the deduplication map $d_z$ is surjective and monotone on first occurrences, i.e., it satisfies $\min d_z^{-1}(i)<\min d_z^{-1}(j)$ if $i<j$. Using \cref{lem:dedup-reindex} and the definition of substitution for atomic formulae, we get:
\begin{equation}
    \bar x(R(\bar z)(\bar y))=\bar x(R(\bar y^{d_z}))=(\bar y ^{d_z})^\circ=\bar y^\circ
\end{equation}

\medskip
\noindent
\emph{Binary:} The substitution clause for the binary connectives yields $(\varphi_1\star\varphi_2)(\bar y)=\varphi_1(\bar y^{c_1})\star\varphi_2(\bar y^{c_2})$. By induction $\bar x(\varphi_i(\bar y^{c_i}))=(\bar y^{c_i})^\circ$, so
\begin{equation}
\bar x\big((\varphi_1\star\varphi_2)(\bar y)\big)
=\big((\bar y^{c_1})^\circ\times(\bar y^{c_2})^\circ\big)^\circ
=(\bar y^{c_1}\times\bar y^{c_2})^\circ
=(y_{s(1)},\ldots,y_{s(n_1+n_2)})^\circ,
\end{equation}
where $s$ is the concatenation of $c_1$ and $c_2$ and $n_i=|\bar x(\varphi_i)|$. By definition, $\bar x(\varphi_1\star\varphi_2)=(\bar x(\varphi_1)\times\bar x(\varphi_2))^\circ$, the map $s$ is surjective with $\min s^{-1}$ strictly increasing, so \cref{lem:dedup-reindex} yields $(y_{s(1)},\ldots,y_{s(n_1+n_2)})^\circ=\bar y^\circ$.

\medskip
\noindent
\emph{Quantifier:} Let $\varphi=\forall z\,\psi$ and let $c$ be the context injection of $\bar x(\varphi)\times z$ relative to $\psi$. By the quantifier clause for substitution,
$\varphi(\bar y)=\forall z'\,\psi((\bar y\times z')^c)$ with $z'$ fresh,
and by the induction hypothesis $\bar x(\psi((\bar y\times z')^c))=((\bar y\times z')^c)^\circ$. Hence $\bar x(\varphi(\bar y))$ is $((\bar y\times z')^c)^\circ$ with $z'$ deleted. Now, the position of $z'$ in $(\bar y\times z')^c$, if it occurs, is that of $z$ in $\bar x(\psi)$, whence deleting $z'$ in the former and deduplicating yields $\bar y^\circ$. 
\end{proof}

\cctrans*
\begin{proof}
By induction on $\varphi$, strengthening the inductive hypothesis to: for every context $\bar y:\bar\tau$ for $\varphi$ with context injection $c$ and every $\bar k\in P_{\bar\tau}^f$, the canonical context of $(\bar y\mid\varphi)_f^{\bar k}$ is $\bar x(\varphi)_f^{\bar k^c}$. By the context clause $(\bar y\mid\varphi)_f^{\bar k}=\varphi_f^{\bar k^c}$, so it suffices to verify, for each $\varphi$ in canonical context $\bar x:\bar\sigma$ and each $\bar k\in P_{\bar\sigma}^f$, that $\bar x(\varphi_f^{\bar k})=\bar x_f^{\bar k}$.  By the sortal and no-overlap conditions on the variable function of a reconstrual, the sequences $((\bar x_i)_f^{k_i}$ are individually duplicate-free and pairwise disjoint. The truth-value and negation cases are immediate.

\medskip
\noindent
\emph{Atomic:} For $R(\bar z):\bar \sigma$ with $| \bar z|=n$, the atomic clause yields $R(\bar{z})_f^{\bar k}=R_f^{\bar{k}^{d}}(\bar{z}_f^{\bar{k}^d})$. By \cref{lem:canon-subst} its canonical context is 
\begin{equation}
    (\bar{z}_f^{\bar{k}^d})^\circ=((z_1)_f^{k_{d(1)}}\times\cdots\times (z_n)_f^{k_{d(n)}})^\circ= (((\bar z^\circ)_{d(1)})_f^{k_{d(1)}}\times\cdots\times ((\bar z^\circ)_{d(n)})_f^{k_{d(n)}})^\circ= (\bar z^\circ)_f^{\bar k},
\end{equation}
where $d:=d_z$ and the first equality expands, the second uses $(\bar z^\circ)^{d_z}=\bar z$ for concatenated sequences, and the third applies \cref{lem:dedup-reindex}.    

\medskip
\noindent
\emph{Binary:} Let $\varphi=\varphi_1\star\varphi_2$. By the clause for binary connectives $\varphi_f^{\bar k}=(\bar x\mid\varphi_1)_f^{\bar k}\star(\bar x\mid\varphi_2)_f^{\bar k}$. By induction, the canonical context of $(\bar x\mid\varphi_i)_f^{\bar k}$ is $\bar x(\varphi_i)_f^{\bar k^{c_i}}$. \cref{lem:dedup-reindex} with $s$ the concatenation of $c_1$ and $c_2$ implies $(\bar x(\varphi_1)_f^{\bar k^{c_1}}\times \bar x(\varphi_2)_f^{\bar k^{c_2}})^\circ=\bar x_f^{\bar k}$.

\medskip
\noindent
\emph{Quantifier:} Let $\varphi=\forall z \, \psi$ and let $\bar x$ be the canonical context of $\varphi$. By the quantifier clause for translations,
\begin{equation}
\varphi_f^{\bar k}=\bigwedge_{h\in P_\tau^f}\forall z_f^h\big(D_\tau^h(z_f^h)\to(\bar x\times z\mid\psi)_f^{\bar k\times h}\big).
\end{equation}
By induction, the canonical context of $(\bar x\times z\mid\psi)_f^{\bar k\times h}$ is $\bar x(\psi)_f^{(\bar k\times h)^c}$, and by \cref{lem:canon-subst}, $D_\tau^h(z_f^h)$ has canonical context $z_f^h$. Since $\bar x$ results from $\bar x(\psi)$ by deleting $z$, the sequence $\bar x(\psi)_f^{(\bar k\times h)^c}$ is $\bar x_f^{\bar k}$ with the subsequence $z_f^h$ possibly inserted at the position of $z$. As these subsequences are duplicate-free and pairwise disjoint, the formula in the scope of the quantifier has canonical context $\bigl(z_f^h\times\bar x(\psi)_f^{(\bar k\times h)^c}\bigr)^\circ=z_f^h\times\bar x_f^{\bar k}$. Binding $z_f^h$ deletes this subsequence, whence each conjunct, and therefore the whole conjunction, has canonical context $\bar x_f^{\bar k}$.
\end{proof}
\medskip
\noindent
For the next lemma, we need a preparatory lemma on iterated substitutions. 
\begin{lemma}[Iterated Substitution]\label{lem:iterated-subst}
Let $\varphi:\bar \sigma$ be an $\mathcal{L}$-formula with canonical context $\bar x:\bar \sigma$, and let $\bar y:\bar \sigma$ be another variable sequence of the same arity. By \cref{lem:canon-subst}, the canonical context of $\varphi(\bar y)$ is $\bar y^\circ:\bar \tau$.
Then, for every variable sequence $\bar z:\bar \tau$,
\begin{equation}
    (\varphi(\bar y))(\bar z)= \varphi(\bar z^{d_y})
\end{equation}
up to renaming of bound variables. 
\end{lemma}
\begin{proof}
By induction on $\varphi$. The cases for truth values and negation are clear.

\medskip
\noindent
\emph{Atomic:} 
Let $\varphi=R(\bar w):\bar \sigma$ be an atomic formula. Then $(\varphi(\bar y))(\bar z)=R(\bar y ^{d_w})(\bar z)$. Applying \cref{cor:dedup_map_surjection}, $d_{\bar y^{d_w}}=d_y\circ d_w$. Therefore,
\begin{equation}
R(\bar y ^{d_w})(\bar z)=R(\bar z^{d_{\bar y^{d_w}}})=R(\bar z^{d_y\circ d_w})=R((\bar z^{d_y})^{d_w})=\varphi(\bar z^{d_y}).    
\end{equation}

\medskip
\noindent
\emph{Binary:} Let $\varphi=\varphi_1\star\varphi_2$. Let $c_i$ be the context injection of $\bar x(\varphi)$ relative to $\varphi_i$. Then
\begin{align}
     \big(\varphi(\bar y)\big)(\bar z)&=
\big(\varphi_1(\bar y^{c_1})\star\varphi_2(\bar y^{c_2})\big)(\bar z)= \big(\varphi_1(\bar y^{c_1})\big)(\bar z^{\hat c_1})
\star
\big(\varphi_2(\bar y^{c_2})\big)(\bar z^{\hat c_2})\\
&=\varphi_1\big((\bar z^{\hat c_1})^{d_{\bar y^{c_1}}}\big)\star \varphi_2\big((\bar z^{\hat c_2})^{d_{\bar y^{c_2}}}\big)=\varphi_1\big((\bar z^{d_y})^{c_1}\big)\star \varphi_2\big((\bar z^{d_y})^{c_2}\big)\\
&=\varphi(\bar z^{d_y}),
\end{align}
where $\hat{c}_i$ is the context injection of $\bar y^\circ$ relative to $\varphi_i(\bar y^{c_i})$. The first equation uses the binary clause for substitution for $\bar y$, the second the same clause for $\bar z$, the third applies the induction hypothesis, the fourth applies \cref{dedup_map_injection}, and the fifth again uses the definition of substitution for $\bar z^{d_y}$.

\medskip\noindent
\emph{Quantifier:} Let $\varphi=\forall_\tau u\,\psi$ and
let $c$ be the context injection of $\bar x(\varphi)\times u$ relative to $\psi$. Choose $v:\tau$, fresh for $\bar y$ and $\psi$, so that $\varphi(\bar y)=\forall v\,\psi((\bar y\times v)^c)$ by the quantifier clause for substitution, and $\bar x(\psi((\bar y\times v)^c))=((\bar y\times v)^c)^\circ$ by \cref{lem:canon-subst}.
Let $\hat c$ be the context injection of $\bar y^\circ\times v$ relative to $\psi((\bar y\times v)^c)$.
Now choose $v':\tau$, fresh for $\bar z$, $\psi$ and $\psi((\bar y\times v)^c)$. Then, using \cref{dedup_map_injection} for the second and the fact that $v$ is fresh for $\bar y$ for the third equality, we get
\begin{equation}
    ((\bar z\times v')^{\hat c})^{d_{(\bar y\times v)^c}}=(\bar z\times v')^{\hat c\circ d_{(\bar y\times v)^c}}
    =(\bar z\times v')^{d_{\bar y\times v}\circ c}
    =(\bar z^{d_{y}}\times v')^{c}.
\end{equation}
Hence, up to renaming of bound variables,
\begin{align}
    \bigl(\varphi(\bar y)\bigr)(\bar z)
    &=\forall v'\,\bigl(\psi((\bar y\times v)^c)\bigr)((\bar z \times v')^{\hat{c}})
    =\forall v'\,\psi\bigl(((\bar z \times v')^{\hat{c}})^{d_{(\bar y\times v)^c}}\bigr)\\
    &=\forall v'\,\psi\bigl((\bar z^{d_{y}}\times v')^{c}\bigr)
    =\varphi(\bar z^{d_{y}}).
\end{align}
The first and last equalities apply the quantifier clause of the definition of substitution to $\forall v\,\psi((\bar y\times v)^c)$ with context injection $\hat c$ and to $\forall u\,\psi$ with context injection $c$, respectively; the second uses the induction hypothesis for $\psi$, and the third is the result of the preceding computation. 
\end{proof}

\begin{lemma}[Substitution and Translation]\label{lem:subst-translation}
     Let $f:\Sigma_1\to\Sigma_2$ be a piecewise reconstrual, $\varphi:\bar \sigma$ an $\mathcal{L}_1$-formula with canonical context $\bar x=\bar x(\varphi)$, and $\bar{y}:\bar{\sigma}$ a sequence of variables. Let $\bar \tau:=\bar\sigma(\varphi(\bar y))$, so that $\bar x(\varphi(\bar y))=\bar y^\circ:\bar\tau$ by \cref{lem:canon-subst}. Then, for any piece $\bar{k}\in P^f_{\bar{\tau}}$,
     \begin{equation}
         \varphi(\bar{y})_f^{\bar{k}} \;=\; \varphi_f^{\bar{k}^{d_{ y}}}(\bar{y}_f^{\bar{k}^{d_{ y}}}),
     \end{equation}
     up to renaming of bound variables, where $(\bar k^{d_{y}})_i=k_{d_{y}(i)}$.
\end{lemma}
\begin{proof}
By induction on $\varphi$; the truth-value and negation cases are immediate.

\medskip
\noindent
\emph{Atomic:} 
Let $\varphi=R(\bar z):\bar \sigma$ be atomic and let $\bar y:\bar\sigma$. For a piece $\bar k\in P_{\bar\sigma(\varphi(\bar y))}$, \cref{cor:dedup_map_surjection} implies $d_{\bar y^{d_z}}=d_{y}\circ d_z$. Hence
\begin{equation}
\varphi(\bar y)_f^{\bar k}
=
R(\bar y^{d_z})_f^{\bar k}=
R_f^{\bar k^{d_{\bar y^{d_z}}}}
\big((\bar y^{d_z})_f^{\bar k^{d_{\bar y^{d_z}}}}\big)=
R_f^{(\bar k^{d_y})^{d_z}}
\big((\bar y^{d_z})_f^{(\bar k^{d_y})^{d_z}}\big).
\end{equation}
On the other hand,
\begin{equation}
\varphi_f^{\bar k^{d_y}}(\bar y_f^{\bar k^{d_y}})=
R(\bar z)_f^{\bar k^{d_y}}(\bar y_f^{\bar k^{d_y}})=
\Big(
R_f^{(\bar k^{d_y})^{d_z}}
\big(\bar z_f^{(\bar k^{d_y})^{d_z}}\big)
\Big)(\bar y_f^{\bar k^{d_y}})
=
R_f^{(\bar k^{d_y})^{d_z}}
\big((\bar y^{d_z})_f^{(\bar k^{d_y})^{d_z}}\big),
\end{equation}
where the last equality uses \cref{lem:iterated-subst}. Therefore $\varphi(\bar y)_f^{\bar k}=\varphi_f^{\bar k^{d_y}}(\bar y_f^{\bar k^{d_y}})$, up to renaming of bound variables.

\medskip
\noindent
\emph{Binary:} 
Let $\varphi=\varphi_1\star\varphi_2$, and let $c_i$ be the context injection of
$\bar x(\varphi)$ relative to $\varphi_i$. For a piece
$\bar k\in P_{\bar\sigma(\varphi(\bar y))}$, let $\hat c_i$ be the context injection of $\bar y^\circ$ relative to $\varphi_i(\bar y^{c_i})$. By
\cref{dedup_map_injection},
\begin{equation}
    d_{y}\circ c_i=\hat c_i\circ d_{\bar y^{c_i}} \text{, and hence } (\bar k^{\hat c_i})^{d_{\bar y^{c_i}}}
    =
    (\bar k^{d_y})^{c_i}.
\end{equation}
Thus, by the binary clause for translations and the induction hypothesis,
\begin{align}
\varphi(\bar y)_f^{\bar k}
&=
\big(\varphi_1(\bar y^{c_1})\star\varphi_2(\bar y^{c_2})\big)_f^{\bar k} 
=
\varphi_1(\bar y^{c_1})_f^{\bar k^{\hat c_1}}
\star
\varphi_2(\bar y^{c_2})_f^{\bar k^{\hat c_2}} \\
&=
\varphi_{1,f}^{(\bar k^{\hat c_1})^{d_{\bar y^{c_1}}}}
\big((\bar y^{c_1})_f^{(\bar k^{\hat c_1})^{d_{\bar y^{c_1}}}}\big)
\star
\varphi_{2,f}^{(\bar k^{\hat c_2})^{d_{\bar y^{c_2}}}}
\big((\bar y^{c_2})_f^{(\bar k^{\hat c_2})^{d_{\bar y^{c_2}}}}\big) \\
&=
\varphi_{1,f}^{(\bar k^{d_y})^{c_1}}
\big((\bar y^{c_1})_f^{(\bar k^{d_y})^{c_1}}\big)
\star
\varphi_{2,f}^{(\bar k^{d_y})^{c_2}}
\big((\bar y^{c_2})_f^{(\bar k^{d_y})^{c_2}}\big)\\
&=\Big(
\varphi_{1,f}^{(\bar k^{d_y})^{c_1}}
\star
\varphi_{2,f}^{(\bar k^{d_y})^{c_2}}
\Big)(\bar y_f^{\bar k^{d_y}}) = \Big(
(\bar x(\varphi)\mid\varphi_1)_f^{\bar k^{d_y}}
\star
(\bar x(\varphi)\mid\varphi_2)_f^{\bar k^{d_y}}
\Big)(\bar y_f^{\bar k^{d_y}})\\
&= \varphi_f^{\bar k^{d_y}}(\bar y_f^{\bar k^{d_y}}),
\end{align}
up to renaming of bound variables. Here, the first and second equalities use the binary clauses of substitution and translation, the third the induction hypothesis, and the fourth \cref{dedup_map_injection}. The fifth uses the binary substitution clause in reverse, together with \cref{lem:variables_translations}, while the sixth uses the context clause, with the final equality using the binary clause for translations again.

\medskip\noindent
\emph{Quantifier:} Let $\varphi=\forall_\tau z\,\psi$,
and let $c$ be the context injection of $\bar x(\varphi)\times z$ relative to $\psi$. Let $\bar y:\bar\sigma$ and $\bar l\in P_{\bar\sigma(\varphi(\bar y))}$.
Choose $v:\tau$, fresh for $\bar y$, $\psi$, and such that, for every $h\in P_\tau^f$, $v_f^h$ is fresh for $\varphi^{\bar l^{d_{y}}}_f$. Then $\varphi(\bar y)=\forall v\,\psi((\bar y\times v)^c)$, and let $\hat c$ be the context injection of $\bar y^\circ\times v$ relative to $\psi((\bar y\times v)^c)$. \cref{dedup_map_injection} then implies $d_{\bar y\times v}\circ c=\hat c\circ d_{(\bar y\times v)^c}$. Given a piece $h \in P_\tau^f$, this implies, together with freshness of $v$ for $\bar y$,
\begin{equation}
    ((\bar l\times h)^{\hat c})^{d_{(\bar y\times v)^c}}=(\bar l\times h)^{\hat c\circ d_{(\bar y\times v)^c}}
    =(\bar l\times h)^{d_{\bar y\times v}\circ c}
    =(\bar l^{d_{y}}\times h)^c.
\end{equation}

\smallskip\noindent
For the left-hand side of the claim: By the quantifier clause for translations, the context clause, and the
induction hypothesis for $\psi$,
\begin{align}
    \varphi(\bar y)^{\bar l}_f
    &=\bigwedge_{h\in P^f_\tau}\forall v^h_f
      \bigl(D^h_\tau(v^h_f)\to
      \psi((\bar y\times v)^c)^{(\bar l\times h)^{\hat c}}_f\bigr)\\
    &=\bigwedge_{h\in P^f_\tau}\forall v^h_f
      \bigl(D^h_\tau(v^h_f)\to
      \psi^{(\bar l^{d_{y}}\times h)^c}_f
      \bigl(((\bar y\times v)^c)^{(\bar l^{d_{y}}\times h)^c}_f\bigr)\bigr),
\end{align} up to renaming of bound variables.

\smallskip\noindent
For the right-hand side: Applying the quantifier clause for translations and the context clause for the first equality,
\begin{align}
    \varphi^{\bar l^{d_{y}}}_f
    (\bar y_f^{\bar l^{d_{y}}})
    &=\Bigl(\bigwedge_{h\in P^f_\tau}\forall z^h_f
      \bigl(D^h_\tau(z^h_f)\to
      \psi^{(\bar l^{d_{y}}\times h)^c}_f\bigr)\Bigr)
      (\bar y_f^{\bar l^{d_{y}}})\\
    &=\bigwedge_{h\in P^f_\tau}\forall v^h_f
      \bigl(D^h_\tau(v^h_f)\to
      \psi^{(\bar l^{d_{y}}\times h)^c}_f
      \bigl(((\bar y\times v)^c)^{(\bar l^{d_{y}}\times h)^c}_f\bigr)\bigr),
\end{align}
up to renaming of bound variables, where the second equality first distributes the substituta over the conjunction, and then, using the freshness condition on $v$, over the implication, applying \cref{lem:iterated-subst}. 
\end{proof}

\compositioncomplex*
\begin{proof}
We work in the theory of identity $L(\Sigma_3)$ in the third signature throughout, abbreviating $\eqv_{L(\Sigma_3)}$ to $\eqv$ in this proof. The claim about canonical contexts follows from the definition of the variable function of the composition and \cref{lem:variables_translations} with the computation
\begin{equation}
\bar x(\varphi_{g\circ f}^{\bar m})
= \bar x(\varphi)_{g\circ f}^{\bar m}
= \big(\bar x(\varphi)_f^{\bar m_f}\big)_g^{\bar m_g}
= \bar x\big((\varphi_f^{\bar m_f})_g^{\bar m_g}\big).    
\end{equation}
\noindent
The remainder proceeds by induction on $\varphi$. Write $\bar m_f=(k_1,\ldots,k_n)$ and $\bar m_g=\bar l_1\times\cdots\times\bar l_n$ as in the statement. The truth-value and negation cases are immediate.

\medskip
\noindent
\emph{Atomic:} Let $\varphi=R(\bar z)$ be atomic, with canonical context $\bar x=\bar z^\circ$. Set
\begin{equation}
    \bar v:=\bar z_f^{\bar m_f^{d_z}},
    \qquad \text{ and }\qquad
    \bar q:=\bar l_{d_z(1)}\times\cdots
        \times\bar l_{d_z(|\bar z|)}.
\end{equation}
Then $(\bar m^{d_z})_f=\bar m_f^{d_z}$ and
$(\bar m^{d_z})_g=\bar q$.
Moreover, $\bar v^\circ=\bar x_f^{\bar m_f}$ and
$\bar m_g^{d_{\bar v}}=\bar q$, by the no-overlap
and duplicate-freeness conditions on the variable function.
Since $R_{g\circ f}^{\bar m}=(R_f^{\bar m_f})_g^{\bar m_g}$,
\begin{align}
R(\bar z)_{g\circ f}^{\bar m}
&= R_{g\circ f}^{\bar m^{d_z}}\big(\bar z_{g\circ f}^{\bar m^{d_z}}\big)
= \big(R_f^{\bar m_f^{d_z}}\big)_g^{\bar q}\big(\bar z_{g\circ f}^{\bar m^{d_z}}\big).
\end{align}
On the other hand, by \cref{lem:subst-translation} for $g$ and the definition of the composite variable function,
\begin{align}
\big(R(\bar z)_f^{\bar m_f}\big)_g^{\bar m_g}
&= \Big(R_f^{\bar m_f^{d_z}}(\bar v)\Big)_g^{\bar m_g}
= \big(R_f^{\bar m_f^{d_z}}\big)_g^{\bar q}(\bar v_g^{\bar q}) \\
&= \big(R_f^{\bar m_f^{d_z}}\big)_g^{\bar q}\big(\bar z_{g\circ f}^{\bar m^{d_z}}\big).
\end{align}
The two right-hand sides coincide, so $R(\bar z)_{g\circ f}^{\bar m}=\big(R(\bar z)_f^{\bar m_f}\big)_g^{\bar m_g}$.

\medskip
\noindent
\emph{Binary:} Let $\varphi=\varphi_1\star\varphi_2$, and let $c_i$ be the context injection of $\bar x(\varphi)$ relative to $\varphi_i$. Set 
\begin{equation}
    \bar q_i:=\bar l_{c_i(1)}\times\cdots
        \times\bar l_{c_i(| \bar x(\varphi_i)|)}
    \qquad(i=1,2).
\end{equation}
Then $(\bar m^{c_i})_f=\bar m_f^{c_i}$ and
$(\bar m^{c_i})_g=\bar q_i$.
Using the induction hypothesis with these pieces,
\begin{align}
\varphi_{g\circ f}^{\bar m}
&= (\bar x(\varphi)\mid\varphi_1)_{g\circ f}^{\bar m}\star(\bar x(\varphi)\mid\varphi_2)_{g\circ f}^{\bar m}
= (\varphi_1)_{g\circ f}^{\bar m^{c_1}}\star(\varphi_2)_{g\circ f}^{\bar m^{c_2}} \\
&\eqv \big((\varphi_1)_f^{\bar m_f^{c_1}}\big)_g^{\bar q_1}\star\big((\varphi_2)_f^{\bar m_f^{c_2}}\big)_g^{\bar q_2}.
\end{align}
Now, for the right-hand side, let $\tilde c_i$ be the context injection of $\bar x(\varphi)_f^{\bar m_f}$ relative to $(\varphi_i)_f^{\bar m_f^{c_i}}$.
By \cref{lem:variables_translations}, $\bar m_g^{\tilde c_i}=\bar q_i$.
Then, by the context clause for $f$, the binary and context translation clauses for $g$, and \cref{lem:variables_translations} to identify the translated contexts,
\begin{align}
\big(\varphi_f^{\bar m_f}\big)_g^{\bar m_g}
&= \Big((\bar x(\varphi)\mid\varphi_1)_f^{\bar m_f}\star(\bar x(\varphi)\mid\varphi_2)_f^{\bar m_f}\Big)_g^{\bar m_g}
= \Big((\varphi_1)_f^{\bar m_f^{c_1}}\star(\varphi_2)_f^{\bar m_f^{c_2}}\Big)_g^{\bar m_g} \\
&= \big((\varphi_1)_f^{\bar m_f^{c_1}}\big)_g^{\bar m_g^{\tilde c_1}}
\star
\big((\varphi_2)_f^{\bar m_f^{c_2}}\big)_g^{\bar m_g^{\tilde c_2}}= \big((\varphi_1)_f^{\bar m_f^{c_1}}\big)_g^{\bar q_1}\star\big((\varphi_2)_f^{\bar m_f^{c_2}}\big)_g^{\bar q_2}.
\end{align}
The two right-hand sides coincide, so $\varphi_{g\circ f}^{\bar m}\eqv\big(\varphi_f^{\bar m_f}\big)_g^{\bar m_g}$.

\medskip
\noindent
\emph{Quantifier:} Let $\varphi=\forall_\tau u\,\psi$ and let $\bar x=\bar x(\varphi)$. By the translation clause for the universal quantifier for $f$,
\begin{equation}
    \varphi_f^{\bar m_f}
= \bigwedge_{h\in P_\tau^f}\forall u_f^h\Big(D_\tau^{f,h}(u_f^h)\to(\bar x\times u\mid\psi)_f^{\bar m_f\times h}\Big).
\end{equation}
Applying $g$ at piece $\bar m_g$ and renaming the bound sequence $u_f^h$ as needed, the conjunct belonging to piece $h$ becomes
\begin{equation}
    \bigwedge_{\bar h'\in P_{\tau_f^h}^g}\forall (u_f^h)_g^{\bar h'}\Big(D_{\tau_f^h}^{g,\bar h'}((u_f^h)_g^{\bar h'})\to\bigl(
        \bar x_f^{\bar m_f}\times u_f^h
        \mid
        D_\tau^{f,h}(u_f^h)
        \to(\bar x\times u\mid\psi)_f^{\bar m_f\times h}
    \bigr)_g^{\bar m_g\times\bar h'}
    \Bigr).
\end{equation}
By the binary and context clauses for $g$ and
\cref{lem:subst-translation} applied to
$D_\tau^{f,h}$, the formula in the scope of the
quantifier is equivalent to
\begin{equation}
    \Bigl(
        D_{\tau_f^h}^{g,\bar h'}((u_f^h)_g^{\bar h'})
        \land
        (D_\tau^{f,h})_g^{\bar h'}((u_f^h)_g^{\bar h'})
    \Bigr)
    \to
    \bigl(
        \bar x_f^{\bar m_f}\times u_f^h
        \mid
        (\bar x\times u\mid\psi)_f^{\bar m_f\times h}
    \bigr)_g^{\bar m_g\times\bar h'}.
\end{equation}
Since its antecedent is equivalent to the composite domain formula $ D_\tau^{g\circ f,(h,\bar h')}((u_f^h)_g^{\bar h'})$, it remains to show that its consequent is equivalent to the consequent in the clause for the quantifier. To that end, let $c$ be the context injection of $\bar x\times u$ relative to $\psi$, and set $\bar p:=\bar m\times(h,\bar h')$. Let $\tilde c$ be the context injection of $\bar x_f^{\bar m_f}\times u_f^h$ relative to
$\psi_f^{\bar p_f^c}$. By \cref{lem:variables_translations}, as in the
binary case, $\bar p_g^{\tilde c}=(\bar p^c)_g$ and $\bar p_f^c=(\bar p^c)_f$. Hence, by the context clauses and the induction hypothesis,
\begin{equation}
    \bigl(
        \bar x_f^{\bar m_f}\times u_f^h
        \mid
        (\bar x\times u\mid\psi)_f^{\bar p_f}
    \bigr)_g^{\bar p_g}=
    \bigl(\psi_f^{\bar p_f^c}\bigr)_g
        ^{\bar p_g^{\tilde c}}
    =
    \bigl(\psi_f^{(\bar p^c)_f}\bigr)_g
        ^{(\bar p^c)_g}
    \eqv
    \psi_{g\circ f}^{\bar p^c}
    =
    (\bar x\times u\mid\psi)_{g\circ f}^{\bar p},
\end{equation}
delivering the desired consequent. Collecting over all $(h,\bar h')\in P_\tau^{g\circ f}$
and using
$u_{g\circ f}^{(h,\bar h')}=(u_f^h)_g^{\bar h'}$,
we hence obtain a formula appropriately equivalent to the result of the quantifier clause for $g\circ f$;  that is, $\big(\varphi_f^{\bar m_f}\big)_g^{\bar m_g}\eqv\varphi_{g\circ f}^{\bar m}$.
\end{proof}

\section{Reflection of Bicategorical Structure}
\label{app:bicategorical_results}

To show that $\Thpw$ is a bicategory, we prove and apply a completely general bicategorical theorem, which should be seen as a generalization of the $1$-categorical fact that faithful functors reflect commutative diagrams. This general theorem is hinted at in \citet[p.~9]{leinster1998basic}, but we were not able to find a precise statement in the literature, which is why we state and prove it here. Essentially, it says that bicategorical structure is reflected back along locally faithful (contravariant) pseudofunctors. While this result could have been omitted in favour of a more concrete approach, we think there is some independent value in having it recorded for future use in similar situations, which are likely to occur in the study of notions of interpretations.

\begin{definition}[Pre-bicategory]\label{def:pre-bicategory}
A \emph{pre-bicategory} $\mathcal{B}$ consists of \emph{a class of objects} $\mathrm{ob}(\mathcal{B})$; for all $X,Y\in\mathrm{ob}(\mathcal{B})$ a class $\mathcal{B}(X,Y)$ of \emph{$1$-morphisms} $f: X\to Y$, and for all $f,f'\in\mathcal{B}(X,Y)$ a class $\mathcal{B}(f,f')$ of \emph{$2$-morphisms} $\theta: f\Rightarrow f': X\to Y$; for each $1$-morphism $f$ an \emph{identity $2$-morphism} $\id_f: f\Rightarrow f$, and for $2$-morphisms $\zeta: f\Rightarrow f'$ and $\theta: f'\Rightarrow f''$ a \emph{vertical composite} $\theta\bullet\zeta: f\Rightarrow f''$; for each $X$ an \emph{identity $1$-morphism} $\id_X: X\to X$; for $1$-morphisms $f: X\to Y$ and $g: Y\to Z$ a \emph{composite} $g\circ f: X\to Z$; and for $2$-morphisms $\zeta: f\Rightarrow f'$ between $1$-morphisms $X\to Y$ and $\theta: g\Rightarrow g'$ between $1$-morphisms $Y\to Z$ a \emph{horizontal composite} $\theta \star \zeta: g\circ f\Rightarrow g'\circ f'$; and $2$-morphisms
        \begin{equation}
          a_{h,g,f}: (h\circ g)\circ f \Rightarrow h\circ (g\circ f): W\to Z
          \end{equation}
          \begin{equation}
              l^{f}: \id_X\circ f \Rightarrow f: W\to X \qquad
          r^{f}: f\circ \id_W \Rightarrow f: W\to X
          \end{equation}
for all $f: W\to X$, $g: X\to Y$ and $h: Y\to Z$.
\end{definition}

Hence a pre-bicategory has the same underlying data as a bicategory, but with no axioms imposed. 
\begin{remark}[Properties of a Bicategory]\label{rem:prebicategory_to_bicategory}
A pre-bicategory is a bicategory in the sense of
\cref{def:bicategory} iff each $\mathcal{B}(X,Y)$ is a category
with compositions given by $\bullet$ and identities by $\id$; horizontal composition is functorial, that is, $\id_g\star \id_f=\id_{g\circ f}$ and the  \emph{interchange law} holds: for
composable $2$-morphisms $\eta_1,\eta_2,\zeta_1,\zeta_2$,
\begin{equation}\label{eq:interchange-abstract}
  (\eta_2\bullet\zeta_2) \star (\eta_1\bullet\zeta_1)
  \;=\;
  (\eta_2 \star \eta_1)\bullet(\zeta_2 \star \zeta_1);
\end{equation}
the families $a$, $l$, $r$ are each componentwise invertible and natural: for all $2$-morphisms $\xi: f\Rightarrow f': W\to X$, $\zeta: g\Rightarrow g': X\to Y$ and $\theta: h\Rightarrow h': Y\to Z$,
        \begin{align}
          a_{h',g',f'}\bullet\bigl((\theta\star\zeta)\star\xi\bigr)
            &= \bigl(\theta\star(\zeta\star\xi)\bigr)\bullet a_{h,g,f},
          \tag{Nat-$a$}\label{ax:nat-a}\\
          l^{f'}\bullet(\id_X\,\xi) &= \xi\bullet l^{f},
          \tag{Nat-$l$}\label{ax:nat-l}\\
          r^{f'}\bullet(\xi\,\id_W) &= \xi\bullet r^{f},
          \tag{Nat-$r$}\label{ax:nat-r}
        \end{align}
and the pentagon and triangle identities hold\footnote{Writing these and similar equations out in diagrammatic form makes it much easier to see their significance; we display them here in this equational form merely to save space.}: 
\begin{align}
  (ka_{h,g,f})\bullet a_{k,h\circ g,f}\bullet(a_{k,h,g}f)
    &= a_{k,h,g\circ f}\bullet a_{k\circ h,g,f},
  \tag{Pentagon}\label{ax:pentagon}\\
  (gl^{f})\bullet a_{g,\id_X,f} &= r^{g}f.
  \tag{Triangle}\label{ax:triangle}
\end{align}
\end{remark}

 Next, we widen the definition of (contravariant) pseudofunctor to allow for pseudofunctors with pre-bicategories as domains and $2$-categories as targets. This is all the generality we need.
\begin{definition}[Contravariant Pseudofunctor from a Pre-Bicategory]\label{def:pseudofunctor_from_prebicategory}
Let $\mathcal{B}$ be a pre-bicategory and $\mathcal{C}$ a strict $2$-category. A \emph{contravariant pseudofunctor} $F:\mathcal{B}\to\mathcal{C}$ comprises a function $F:\mathrm{ob}(\mathcal{B})\to \mathrm{ob}(\mathcal{C})$, functions $F_{X,Y}:\mathcal{B}(X,Y)\to \mathcal{C}(FY,FX)$ and $F_{f,g}:\mathcal{B}(f,g)\to\mathcal{C}(Ff,Fg)$ for all $X,Y\in \mathrm{ob}(\mathcal{B})$ and $f,g\in \mathcal{B}(X,Y)$, together with invertible $2$-morphisms of $\mathcal{C}$
\begin{equation}
    \alpha_{g,f}: F(g\circ f)\Rightarrow Ff\circ Fg
    \qquad\text{and}\qquad
    \alpha_{X}: F(\id_X)\Rightarrow \id_{FX},
\end{equation} such that
\begin{align}
    &F(\theta\bullet\zeta)=F\theta\bullet F\zeta,
    \qquad
    F(\id_f)=\id_{Ff},
    \tag{CF-HomFunct}\label{ax:CF_HomFunct}\\
    &\alpha_{g_2,g_1}\bullet F(\theta_2\star\theta_1)
      =\bigl(F\theta_1\star F\theta_2\bigr)\bullet\alpha_{f_2,f_1},
    \tag{CF-CompNat}\label{ax:CF_comp_nat}\\
    &\bigl(\alpha_{g,f}\,Fh\bigr)\bullet\alpha_{h,g\circ f}\bullet F(a_{h,g,f})
      =\bigl(Ff\,\alpha_{h,g}\bigr)\bullet\alpha_{h\circ g,f},
    \tag{CF-Assoc}\label{ax:CF_assoc}\\
    &F(l^{f})=\bigl(Ff\,\alpha_X\bigr)\bullet\alpha_{\id_X,f}, \tag{CF-Left}\label{ax:CF_left}\\ 
    &F(r^{f})=\bigl(\alpha_W\,Ff\bigr)\bullet\alpha_{f,\id_W}.
    \tag{CF-Right}\label{ax:CF_right}
\end{align}
hold for all $\zeta:u\Rightarrow v, \theta:v\Rightarrow w$, $\theta_1:f_1\Rightarrow g_1, \theta_2:f_2\Rightarrow g_2$ with $f_1,f_2$ and $g_2,g_1$ composable, and $f: W\to X$, $g: X\to Y$ and $h: Y\to Z$. It is \emph{locally faithful} if it is injective on $2$-morphisms, that is, iff the following holds for all parallel $2$-morphisms $\eta,\zeta:f\Rightarrow g$:
\begin{align}
    &F\eta=F\zeta\;\Rightarrow\;\eta=\zeta
    \tag{CF-LF}\label{CF_locally_faithful}
\end{align}
\end{definition} 
It is clear that the preceding conditions ensure that any such contravariant pseudofunctor is a contravariant pseudofunctor between bicategories if $\mathcal{B}$ is a bicategory, as the
conditions are exactly the axioms of \cref{def:pseudofunctor} for a
pseudofunctor $\mathcal{B}^{\mathrm{op}}\to\mathcal{C}$, with $\alpha$ replaced by its inverses, as fixed after \cref{def:pseudofunctor}. We can now precisely formulate the following useful theorem: 
\reflection*
\medskip\noindent
Before the proof, note that we could just as easily have formulated a covariant version of the theorem, but the contravariant version is the one we need and avoids the need to define opposite pre-bicategories and other dual notions.  Moreover, the strictness assumption on $\mathcal{C}$ can be relaxed using the coherence theorem for bicategories \citep[see, e.g.,][Section 2.3]{leinster1998basic}.
\begin{proof}
We rely on local faithfulness and strictness of $\mathcal{C}$ throughout: to show that an equation between $2$-morphisms holds in $\mathcal{B}$, it is enough to apply $F$ to it and show that it holds in $\mathcal{C}$.

\medskip\noindent
For objects $X,Y\in \mathrm{ob}(\mathcal{B})$, the assumption \eqref{ax:CF_HomFunct} and associativity of vertical composition in $\mathcal{C}$ imply that $(\mathcal{B}(X,Y),\bullet, \id)$ is a category. We next show that composition is functorial, that is, that horizontal and vertical composition obey the interchange law, and that $\id_g\star \id_f=\id_{g\circ f}$. The latter follows easily from \eqref{ax:CF_comp_nat} and \eqref{ax:CF_HomFunct}, since $F(\id_g\star\id_f)
 =\alpha_{g,f}^{-1}\bullet(\id_{Ff}\star\id_{ Fg})\bullet\alpha_{g,f}
 =\id_{F(g\circ f)}=F(\id_{g\circ f})$. For the former, let $\zeta_i:f_i\Rightarrow g_i$ and $\theta_i:g_i\Rightarrow h_i$ ($i=1,2$) be $2$-morphisms in $\mathcal{B}$. We compute:   
\begin{align}
    F\bigl((\theta_2\bullet\zeta_2)\star(\theta_1\bullet\zeta_1)\bigr)
    &=\alpha_{h_2,h_1}^{-1}\bullet
      \bigl((F\theta_1\bullet F\zeta_1)\star(F\theta_2\bullet F\zeta_2)\bigr)
      \bullet\alpha_{f_2,f_1}\\
      &=\alpha_{h_2,h_1}^{-1}\bullet(F\theta_1\star F\theta_2)
      \bullet(F\zeta_1\star F\zeta_2)\bullet\alpha_{f_2,f_1}\\
    &=F(\theta_2\star\theta_1)\bullet F(\zeta_2\star\zeta_1)
     =F\bigl((\theta_2\star\theta_1)\bullet(\zeta_2\star\zeta_1)\bigr)
\end{align}
Here, the first equality applies \eqref{ax:CF_comp_nat} and \eqref{ax:CF_HomFunct}. The second uses the interchange law of $\mathcal{C}$, the third applies \eqref{ax:CF_comp_nat} twice, and the fourth concludes the computation with \eqref{ax:CF_HomFunct}. 

\medskip\noindent
Moving on to naturality of the associator, let $f,f'\in \mathcal{B}(W,X), g,g'\in \mathcal{B}(X,Y)$, and $h,h'\in \mathcal{B}(Y,Z)$ be $1$-morphisms, and let $\theta_1:f\Rightarrow f', \theta_2:g\Rightarrow g'$, and $\theta_3:h\Rightarrow h'$ be $2$-morphisms. We have to show that $a_{h',g',f'}\bullet\bigl((\theta_3\star \theta_2) \star \theta_1\bigr)
   \bullet a_{h,g,f}^{-1}=\bigl(\theta_3\star(\theta_2\star\theta_1)\bigr)$. Together with \eqref{ax:CF_HomFunct}, two applications of \eqref{ax:CF_comp_nat} and the interchange law yield that
\begin{align}
    &F\bigl(a_{h',g',f'}\bullet\bigl((\theta_3\star \theta_2) \star \theta_1\bigr)
   \bullet a_{h,g,f}^{-1}\bigr)=\\
   &F(a_{h',g',f'})\bullet\alpha_{h'g',f'}^{-1}
   \bullet\bigl(Ff'\,\alpha_{h',g'}^{-1}\bigr)\bullet\Theta
   \bullet\bigl(Ff\,\alpha_{h,g}\bigr)\bullet\alpha_{hg,f}
   \bullet F(a_{h,g,f})^{-1},
\end{align} where $\Theta=F(\theta_1)\star F(\theta_2)\star F(\theta_3)$. Likewise, one sees that \begin{equation}\label{eq:assoc-nat-rhs}
    F\bigl(\theta_3\star(\theta_2\star\theta_1)\bigr)
    =\alpha_{h',\,g'\circ f'}^{-1}
       \bullet\bigl(\alpha_{g',f'}^{-1}\, Fh'\bigr)
     \bullet\Theta\bullet
     \bigl(\alpha_{g,f}\, Fh\bigr)\bullet\alpha_{h,\,g\circ f}
\end{equation} Using that whiskering preserves inverses and vertical composites, and \eqref{ax:CF_assoc}, one sees that the respective factors on the left and right agree.
   
\medskip\noindent
Turning to naturality of the unitors, let $\theta:f\Rightarrow f'$ with $f,f'\in\mathcal{B}(X,Y)$. We need to show that $l^{f'}\bullet(\id_{\id_Y}\star\theta)=\theta\bullet l^{f}$. We compute:
\begin{align}
    &F\bigl(l^{f'}\bullet(\id_{\id_Y}\star\theta)\bigr)
    =\bigl(F f'\,\alpha_Y\bigr)\bullet\bigl(F\theta\star\id_{F\id_Y}\bigr)
      \bullet\alpha_{\id_Y,f}\\
      &=\bigl(F\theta\star\alpha_Y\bigr)\bullet\alpha_{\id_Y,f}
    =\bigl(F\theta\star\id_{\id_{FY}}\bigr)
      \bullet\bigl(F f\,\alpha_Y\bigr)\bullet\alpha_{\id_Y,f}=F\bigl(\theta\bullet l^{f}\bigr)
\end{align}
The first equality uses \eqref{ax:CF_HomFunct} and expands
$F(l^{f'})$ via \eqref{ax:CF_left} and $F(\id_{\id_Y}\star\theta)$ via
\eqref{ax:CF_comp_nat}, cancelling resulting inverse compositors. The second and third result from interchange in $\mathcal{C}$. The final equality uses that $F\theta\star\id_{\id_{FY}}=F\theta$ because of strictness of $\mathcal{C}$ and concludes with \eqref{ax:CF_left} and
\eqref{ax:CF_HomFunct}. The computation for the right unitor is analogous.

\medskip\noindent
We now show that the pentagon identities hold in $\mathcal{B}$. Let $f: W\to X$, $g: X\to Y$, $h: Y\to Z$ and
$k: Z\to V$ be $1$-morphisms of $\mathcal{B}$. The pentagon identity
\eqref{ax:pentagon} says that the composites along
\begin{equation}\label{eq:pent-long}
    ((kh)g)f
    \xRightarrow{\;a_{k,h,g}f\;} (k(hg))f
    \xRightarrow{\;a_{k,hg,f}\;} k((hg)f)
    \xRightarrow{\;k\,a_{h,g,f}\;} k(h(gf))
\end{equation}
and
\begin{equation}\label{eq:pent-short}
    ((kh)g)f
    \xRightarrow{\;a_{kh,g,f}\;} (kh)(gf)
    \xRightarrow{\;a_{k,h,gf}\;} k(h(gf))
\end{equation}
agree. By local faithfulness it suffices to show that their images under $F$ agree.

\medskip\noindent
Set $N:=Ff\circ Fg\circ Fh\circ Fk$, a $1$-morphism $FV\to FW$ of
$\mathcal{C}$, using strictness to omit the brackets. The image of each of the five $1$-morphisms occurring in
\eqref{eq:pent-long} and \eqref{eq:pent-short} under $F$ is
isomorphic to $N$ via the following $2$-morphisms, which one gets by using compositors and their whiskerings to remove the brackets from the outside inwards, and which are hence clearly invertible. 
\begin{align}
    A_1&:=\bigl(Ff Fg\,\alpha_{k,h}\bigr)
          \bullet\bigl(Ff\,\alpha_{kh,g}\bigr)\bullet\alpha_{(kh)g,f}&& : F\bigl(((kh)g)f\bigr)\Rightarrow N,\label{eq:A1}\\
    A_2&:=\bigl(Ff\,\alpha_{h,g}\,Fk\bigr)
          \bullet\bigl(Ff\,\alpha_{k,hg}\bigr)\bullet\alpha_{k(hg),f}
    &&: F\bigl((k(hg))f\bigr)\Rightarrow N,\label{eq:A2}\\
    A_3&:=\bigl(Ff\,\alpha_{h,g}\,Fk\bigr)
          \bullet\bigl(\alpha_{hg,f}\,Fk\bigr)\bullet\alpha_{k,(hg)f}
    &&: F\bigl(k((hg)f)\bigr)\Rightarrow N,\label{eq:A3}\\
    A_4&:=\bigl(\alpha_{g,f}\, FhFk\bigr)
          \bullet\bigl(\alpha_{h,gf}\,Fk\bigr)\bullet\alpha_{k,h(gf)}
    &&: F\bigl(k(h(gf))\bigr)\Rightarrow N,\label{eq:A4}\\
    A_5&:=\bigl(Ff Fg\,\alpha_{k,h}\bigr)
          \bullet\bigl(\alpha_{g,f}\,F(kh)\bigr)\bullet\alpha_{kh,gf}
    &&: F\bigl((kh)(gf)\bigr)\Rightarrow N.\label{eq:A5}
\end{align}

\medskip\noindent
The following equations, proved below, express that the $A_i$ together yield a cocone for the image of the pentagon. Equivalently, the first three give a cocone for the long composite and the last two a cocone for the short composite above. 
\begin{align}
    A_2\bullet F\bigl(a_{k,h,g}f\bigr)&=A_1,
    &
    A_3\bullet F\bigl(a_{k,hg,f}\bigr)&=A_2,
    &
    A_4\bullet F\bigl(k\,a_{h,g,f}\bigr)&=A_3,
    \tag{Cocone-L}\label{eq:cocone-long}\\
    A_5\bullet F\bigl(a_{kh,g,f}\bigr)&=A_1,
    &
    A_4\bullet F\bigl(a_{k,h,gf}\bigr)&=A_5.
    &&
    \tag{Cocone-S}\label{eq:cocone-short}
\end{align}
Using this and invertibility of the $A_i$,
\eqref{eq:cocone-long} and \eqref{eq:cocone-short}, one sees that the images of the composites along
\eqref{eq:pent-long} and \eqref{eq:pent-short} are telescopic, yielding the desired equality:
\begin{align}
    F\Bigl((k\,a_{h,g,f})\bullet a_{k,hg,f}\bullet(a_{k,h,g}f)\Bigr)
    &=F\bigl(k\,a_{h,g,f}\bigr)\bullet F\bigl(a_{k,hg,f}\bigr)
      \bullet F\bigl(a_{k,h,g}f\bigr)\\
    &=\bigl(A_4^{-1}\bullet A_3\bigr)\bullet\bigl(A_3^{-1}\bullet A_2\bigr)
      \bullet\bigl(A_2^{-1}\bullet A_1\bigr)\\
    &=A_4^{-1}\bullet A_1
     =\bigl(A_4^{-1}\bullet A_5\bigr)\bullet\bigl(A_5^{-1}\bullet A_1\bigr)\\
    &=F\bigl(a_{k,h,gf}\bigr)\bullet F\bigl(a_{kh,g,f}\bigr)
     =F\bigl(a_{k,h,gf}\bullet a_{kh,g,f}\bigr),
\end{align}
where the first and last equalities follow by \eqref{ax:CF_HomFunct}. 

\medskip\noindent
It remains to prove \eqref{eq:cocone-long} and \eqref{eq:cocone-short}. For the identity for $(A_2,A_1)$, we see that
\begin{align}
    A_2\bullet F\bigl(a_{k,h,g}f\bigr)
    &=\Bigl(Ff\,\bigl[(\alpha_{h,g}\,Fk)\bullet\alpha_{k,hg}
      \bullet Fa_{k,h,g}\bigr]\Bigr)\bullet\alpha_{(kh)g,f}\\
    &=\Bigl(Ff\,\bigl[( Fg\,\alpha_{k,h})\bullet\alpha_{kh,g}\bigr]\Bigr)
      \bullet\alpha_{(kh)g,f}
    =A_1,
\end{align}
where we use first the instance of \eqref{ax:CF_comp_nat} for $F(a_{k,h,g}\star\id_f)$, \eqref{ax:CF_HomFunct}, and the fact that whiskering commutes with vertical composition, and then \eqref{ax:CF_assoc} for $a_{k,h,g}$ to conclude. An analogous argument for $(A_4,A_3)$, combined with a right whiskering by $Fk$ in the second step, yields
\begin{align}
    A_4\bullet F\bigl(k\,a_{h,g,f}\bigr)
    &=\Bigl(\bigl[(\alpha_{g,f}\,Fh)\bullet\alpha_{h,gf}\bullet Fa_{h,g,f}\bigr]
      \,Fk\Bigr)\bullet\alpha_{k,(hg)f}\\
    &=\Bigl(\bigl[(Ff\,\alpha_{h,g})\bullet\alpha_{hg,f}\bigr]\,Fk\Bigr)
      \bullet\alpha_{k,(hg)f}
    =A_3.
\end{align}
The remaining identities, for $(A_3,A_2)$, $(A_5,A_1)$, and $(A_4,A_5)$, follow from \eqref{ax:CF_assoc} for the triples $(k,hg,f)$, $(kh,g,f)$, and $(k,h,gf)$, respectively, by postcomposing with the leftmost factor in the definition of $A_3$, $A_5$, and $A_4$, respectively, using interchange in the last case.

\medskip\noindent
Finally and more directly, we prove the triangle identities: Let $f: W\to X$ and $g: X\to Y$ be $1$-morphisms in $\mathcal{B}$. By
\eqref{ax:CF_comp_nat} for $F(\id_g\star l^{f})$ and
$F(r^{g}\star\id_f)$, then \eqref{ax:CF_left} for $f$ and
\eqref{ax:CF_right} for $g$, and finally by the fact that whiskering commutes with vertical composition,
\begin{align}
    \alpha_{g,f}\bullet F\bigl(g\,l^{f}\bigr)
    &=\bigl(Fl^{f}\, Fg\bigr)\bullet\alpha_{g,\,\id_X\circ f}
     =\bigl( Ff\,\alpha_X\, Fg\bigr)
      \bullet\bigl(\alpha_{\id_X,f}\, Fg\bigr)
      \bullet\alpha_{g,\,\id_X\circ f}=:B_2,\label{eq:tri-B2}\\
    \alpha_{g,f}\bullet F\bigl(r^{g}f\bigr)
    &=\bigl( Ff\,Fr^{g}\bigr)\bullet\alpha_{g\circ\id_X,\,f}
     =\bigl( Ff\,\alpha_X\, Fg\bigr)
      \bullet\bigl( Ff\,\alpha_{g,\id_X}\bigr)
      \bullet\alpha_{g\circ\id_X,\,f}=:B_1,\label{eq:tri-B1}
\end{align}
both $2$-morphisms in $\mathcal{C}$ to $ Ff\circ Fg$, noting that the identity in $\mathcal{C}$ is strict. Post-composing the instance of \eqref{ax:CF_assoc} for $F(a_{g,\id_X,f})$ with $\bigl(F f\,\alpha_X\,F g\bigr)$ yields
$B_2\bullet F(a_{g,\id_X,f})=B_1$, whence
\begin{equation}
    \alpha_{g,f}\bullet F\Bigl((g\,l^{f})\bullet a_{g,\id_X,f}\Bigr)
    =B_2\bullet F\bigl(a_{g,\id_X,f}\bigr)
    =B_1
    =\alpha_{g,f}\bullet F\bigl(r^{g}f\bigr),
\end{equation}
using \eqref{ax:CF_HomFunct} for the first equality. One concludes with invertibility of $\alpha_{g,f}$.

\end{proof}

   \end{appendices}
\end{document}